\documentclass[11pt,authoryear]{elsarticle}
\usepackage[ansinew]{inputenc}
\usepackage{a4wide}
\usepackage{calc}
\usepackage{times}
\usepackage{amssymb,amsmath,amstext,amsthm,amsfonts, ams fonts}
\usepackage{colortbl,color,graphics,graphicx,epsfig}
\definecolor{gray}{rgb}{0.8,0.8,0.8}
\usepackage{arydshln}
\usepackage{dsfont}
\usepackage{bbm}
\usepackage{booktabs}
\usepackage{lineno}
\usepackage{arydshln}
\usepackage{multirow}
\usepackage{natbib}
\usepackage[onehalfspacing]{setspace}
\usepackage[pdflinkmargin=5pt,pdfstartview={FitBH -32768},plainpages=false]{hyperref}
\renewcommand{\P}{\mathbb{P}}

\newcommand{\E}{\mathbb{E}}
\newcommand{\N}{\mathds{N}}
\newcommand{\R}{\mathds{R}}
\newcommand{\1}{\mathbbm{1}}
\newcommand{\KLEINO}{{\scriptstyle{\mathcal{O}}}}
\DeclareMathAccent{\verywidehat}{\mathord}{largesymbols}{'144}
\newcommand{\var}{\mathbb{V}\hspace*{-0.05cm}\textnormal{a\hspace*{0.02cm}r}}

\newcommand{\cov}{\mathbb{C}\textnormal{o\hspace*{0.02cm}v}}
\newcommand{\bk}{\|\Phi_{jk}\|_n^{-2}\frac{\hat\eta_{kh_n}}{n}}
\newcommand{\ind}{\1_{\{h_n|\zeta_k^{ad}(Y)|\le u_n\}}}

\renewcommand{\:}{\mathrel{\mathop{:}}}
\allowdisplaybreaks[3]
\newtheorem{remark}{Remark}
\newtheorem{theo}{Theorem}
\newtheorem{assump}{Assumption}
\newtheorem{prop}{Proposition}[section]
\newtheorem{lem}{Lemma}
\newtheorem{cor}[prop]{Corollary}

\newtheorem{strassump}{Assumption}

\newtheorem{sstrassump}{Assumption}

\newtheorem{nassump}{Assumption}

\newcounter{assumpl}
\makeatletter
\newcommand*{\rom}[1]{\expandafter\@slowromancap\romannumeral #1@}
\makeatother
\begin{document}
\renewcommand*{\thefootnote}{\fnsymbol{footnote}}
%\title{Common jumps of price and volatility in noisy high-frequency data}
\title{\Large{Common price and volatility jumps in noisy high-frequency data}}
\normalsize
\author{Markus Bibinger$^a$ \\ Lars Winkelmann$^b$\\%\footnote[1]{Financial support from the Deutsche Forschungsgemeinschaft via CRC 649 `\"Okonomisches Risiko', Humboldt-Universität zu Berlin, is gratefully acknowledged.} \\
\emph{\small $^a$Faculty of Mathematics and Computer Science, Philipps-Universität Marburg}\\
\emph{\small $^b$Department of Economics, Freie Universit\"at Berlin}}
%\author{Markus Bibinger, Lars Winkelmann}
%\author[2]{Lars Winkelmann}
%\address[1]{Institut f\"ur Mathematik, Humboldt-Universit\"at zu Berlin, Unter den Linden 6, 10099 Berlin, Germany}
%\address[2]{Institut für Statistik und \"Okonometrie, Freie Universit\"at Berlin, Boltzmannstraße 20, 14195 Berlin, Germany}
\normalsize
\begin{frontmatter}
\date{\today}

%\maketitle

\begin{abstract}
We introduce a statistical test for simultaneous jumps in the price of a financial asset and its volatility process. The proposed test is based on high-frequency data
and is robust to market microstructure frictions. For the test, local estimators of volatility jumps at price jump arrival times are designed using a nonparametric spectral estimator of the spot volatility process. A simulation study and an empirical example with NASDAQ order book data demonstrate the practicability of the proposed methods and highlight the important role played by price volatility co-jumps.
\end{abstract}

%\vspace{.25cm}
\begin{keyword}
%% keywords here, in the form: keyword \sep keyword
High-frequency data \sep microstructure noise \sep nonparametric volatility estimation\sep volatility jumps\\[.25cm]
%% MSC codes here, in the form: \MSC code \sep code
%% or \MSC[2008] code \sep code (2000 is the default)
%C 13, C 58
\end{keyword}

\end{frontmatter}
\thispagestyle{plain}

\section{Introduction\label{sec:1}}
\renewcommand*{\thefootnote}{\arabic{footnote}}
\setcounter{footnote}{0}
%Testing for the presence of jumps from noisy observations has been addressed, for instance, by \cite{lietal}.
In recent years the broad availability of high-frequency intra-day financial data has spurred a considerable collection of works dedicated to statistical modeling and inference for such data. Semimartingales are a general class of time-continuous stochastic processes to model dynamics of intra-day log-prices in accordance with standard no arbitrage conditions. 
%such that the `no free lunch with vanishing risk' property is satisfied. 
We consider a general It\^{o} semimartingale log-price model allowing for stochastic volatility, price and volatility jumps as well as leverage. %Due to the market microstructure of financial data recorded at high frequencies, as effects of transaction costs and bid-ask bounce, log-prices are not directly well fitted by semimartingales. Instead, a noisy observation model turned out to be more suitable. Taking microstructure frictions into account substantially changes statistical properties and involved mathematical concepts of estimators.\\[.1cm]
%At discrete observation times the noise is interfering the latent efficient log-price evolution described by an It\^{o} semimartingale.\\[.1cm]
%One core research topic in statistics, finance and econometrics of high-frequency data is inference on the (integrated) volatility, bringing forth the seminal contributions by \cite{andersenbollerslev98}, \cite{abdl01}, \cite{bn3}, \cite{howoften} and much more literature devoted to this aspect. 
%Reliable estimates of volatility are of key importance in the decision making process of portfolio and risk managers, see, for instance, \cite{modeling}, as well as policy makers, see \cite{dewachter}. 
%As the volatility takes a leading role in the model, it is important to set up accurate stochastic volatility models, see \cite{eraker} among others. 
Uncertainty and risk in these models are usually ascribed to two distinct sources: First, the volatility process of the continuous semimartingale part that permanently influences observed returns and, second, occasional jumps in prices. Modeling and inference on the two components constitutes a core research topic in statistics, finance and econometrics bringing forth the seminal contributions by \cite{andersenbollerslev98}, \cite{abdl01}, \cite{bn3}, \cite{howoften} and much more literature devoted to this aspect. 
%The latter reflect updates of markets' expectations in response to firm specific news, macro or monetary policy events. An important question, often left unaddressed in the literature, is if one should incorporate jumps also in the volatility process. First studies by \cite{eraker2} and \cite{voljumps} suggest to do so and highlight the important implications especially for asset pricing. A natural question arises, if prices and their volatilities jump at common times stimulated by the same events, or not. 
For asset pricing (\cite{duffie}, \cite{premium}), macro and monetary economics (\cite{ecb}) and risk management 
 (\cite{liufinance}) information about jumps is of key importance. While the literature on price jumps is well developed from both a statistical and empirical
 point of view, methods and evidence about volatility jumps are lagging behind. Empirical evidence about volatility jumps is usually based on methods for price jumps applied to an observable volatility measure like the index of implied volatility of S\&P 500 index options (VIX), see \cite{bloom} and \cite{voljumps}. Such modeling strategies inevitably restrict the number of target variables and the overall scope of empirical insights. Since price jumps have often been associated with macro announcements or firm specific news, a natural empirical question arises, if prices and their volatilities jump at common times stimulated by the same events, or not. Such common jumps of price and volatility are often excluded in the statistics literature to avoid technical difficulties. Beyond the question if one should include simultaneous jump times in price and volatility in a model, testing locally for volatility jumps opens up new ways to study effects of information processing and volatility persistence. This is also reflected in an increasing interest to separate the leverage effect in a continuous and a jump part in the current literature, see \cite{lev} and \cite{xiu2}. The asset pricing model of \cite{pastor} illustrates economic forces behind contemporaneous price and volatility jumps. In their model, agents learn about the profitability of a firm in a changing political environment. A change in government policy does not only affect the expected profitability of a firm (price jump) but also triggers a simultaneous volatility jump induced by the impact uncertainty of the new policy. \\[.1cm]
This article presents a statistical test to decide whether intra-day log-prices exhibit common price and volatility jumps. Our main contribution is to extend the pioneering works by \cite{jacodtodorov} and \cite{bandireno} and to provide an approach for an observation model that accounts for market microstructure. It is widely acknowledged that due to market microstructure of financial data recorded at high frequencies, as effects of transaction costs and bid-ask bounce, log-prices are not directly adequately modeled by semimartingales. %Instead, a noisy observation model turned out to be more suitable. 
Taking microstructure frictions into account substantially changes statistical properties and involved mathematical concepts of estimators. We introduce a spectral spot volatility estimator for noisy observations. The test generalizes the theory by \cite{jacodtodorov} for non-noisy observations. We obtain a statistical test by a neat combination of a stable central limit theorem at (almost) optimal rate for the spectral spot volatility estimator and a suitable test function. In analogy to \cite{jacodtodorov}, the new test is self-scaling in the volatility and rate-optimal. Those two properties are crucial to obtain an efficient method. %They rely deeply on the features of the spectral spot volatility estimator which, to the best of the authors' knowledge, is the only estimator for which a stable central limit theorem at optimal rate is available, where the asymptotic variance does not have a more complex structure. %In contrast, with a more complex asymptotic variance structure, as for instance given by a pre-average or realized kernel spot estimator, preserving the self-scaling property does not seem to be feasible. 
The development of a test that can cope with noise is of high relevance and importance as \cite{jacodtodorov} already remark in their empirical application: ``presence of microstructure noise in the prices is nonnegligible''. We show in simulations that compared to an application of the method by \cite{jacodtodorov} based on skip-sampled returns, we can significantly improve the power of the test.\\ % and ``an extension of our tests, while building on the theoretical results here, asks for a significantly more involved mathematical approach which goes beyond the scope of the current paper and is thus left for future work''.\\
Jumps in prices and the volatility are of very different nature. Large price jumps become visible through large returns. More precisely, in a high-frequency context truncation techniques as suggested by \cite{mancini}, \cite{leemykland} and \cite{jacodjumps} can be used to identify returns that involve jumps. Up to some subtle changes due to dilution by microstructure, this remains valid also in the noisy setup. However, the localization of jump times becomes less precise and more difficult under noise. A first localization method for price jumps in the noisy semimartingale model has been introduced by \cite{fanwang} using wavelets. Other localization approaches are included in \cite{leemykland2} and in \cite{bibwink2015}. We adopt the methods from \cite{bibwink2015} to estimate the spot volatility in presence of price jumps and also to locate price-jump times by thresholding. Contrarily to price jumps, volatility jumps are latent and not as obvious as price jumps due to the fact that we can not observe the volatility path. %We thus infer on the latent volatility of an efficient log-price process from indirect observations of these efficient log-prices diluted by microstructure noise. 
The key element to determine volatility jumps will be efficient estimates of the instantaneous volatility from observed prices.\\
Our spectral spot volatility estimator relies on the Fourier method promoted by \cite{reiss} and  \cite{BHMR} for estimating quadratic (co-)variation, combined with truncation techniques of \cite{bibwink2015} to deal with price jumps. These methods attain lower variance bounds for integrated volatility estimation from noisy observations and are, compared to simple smoothing methods and especially skip-sampling to lower observation frequencies, more efficient. While we are the first who address the testing problem under noise, consistent spot volatility estimators under noise are available. \cite{zu14} and \cite{manc15} present local two-scales estimators and prove stable central limit theorems. The construction of a rate-optimal pre-average estimator is sketched in Section 8.7 of \cite{sahaliajacod}. An alternative approach considering deterministic volatility is presented in \cite{munk2010b}. For our estimator, we establish rate-optimality and a stable central limit theorem with smaller asymptotic variance compared to the pre-average approach. The asymptotic theory allows for general heteroscedastic, serially correlated and endogenous noise. %Price jumps are recovered using a truncation procedure which can be
%adapted to the local magnitude and intraday shape of volatility. 
With this estimation approach at hand, we design a test, comparing estimated local volatilities and their left limits at the estimated price-jump times. As a special case, this includes a local test for volatility jumps at some fixed time or stopping time. %For instance, one might want to test for a volatility jump at news arrival times. 
A test with fast convergence rate based on second order asymptotics of the estimator is suggested. While the overarching strategy follows \cite{jacodtodorov}, the specific test function and construction in the noisy observation case are different and profit from the spectral estimation methodology. Compared to previous estimation techniques to smooth noise, the asymptotic variance structure of the spectral volatility estimates in Theorem \ref{cltspot} admits a simpler form. This facilitates a test statistic which is self-scaling in the local volatility and thus furnishes an asymptotic distribution free test with the best possible rate. 
%Using different estimation techniques to smooth noise as realized kernels motivated by \cite{bn2}, or pre-averaging by \cite{JLMPV}, such a construction, if possible at all, will be more cumbersome. 
The Monte Carlo study corroborates the high precision of the methods in finite samples. Our data study shows that price volatility co-jumps occur and are practically relevant.\\[.1cm]
The paper is organized as follows. Section \ref{sec:2} introduces the model and the statistical problem. We discuss the main ideas for the construction of the test including a short review of the approach for non-noisy data. Section \ref{sec:2.2} describes the spectral spot volatility estimation. We state and discuss the assumptions imposed on the model for the asymptotic theory in Section \ref{sec:3.1} before presenting the main results in Section \ref{sec:3.2}. Practical guidance for the implementation and a Monte Carlo study are given in Section \ref{sec:4}. In Section \ref{sec:5} the methods are used to analyze price and volatility jumps in NASDAQ high-frequency intra-day trading data, reconstructed from the order book.  Section \ref{sec:6} concludes. All proofs are gathered in Section \ref{sec:7}.

\section{Model, testing problem and statistical approach\label{sec:2}}
Let $(\Omega^X,\mathcal{F}^X,(\mathcal{F}^X_t),\P^X)$ be a filtered probability space satisfying the usual conditions. The latent log-price process $X$ follows an It\^{o} semimartingale
\begin{align}\label{sm}X_t =X_0+\int_0^tb_s\,ds+\int_0^t\sigma_s\,dW_s&+\int_0^t\int_{\mathds{R}}\delta(s,x)\1_{\{|\delta(s,x)|\le 1\}}(\mu-\nu)(ds,dx)\\
& \notag  +\int_0^t\int_{\mathds{R}}\delta(s,x)\1_{\{|\delta(s,x)|> 1\}}\mu(ds,dx)\,,\end{align}
with $W$ an $(\mathcal{F}^X_t)$-adapted standard Brownian motion, $\mu$ a Poisson random measure on $\mathds{R}_+\times \mathds{R}$ with $\mathds{R}_+=[0,\infty)$ and an intensity measure (predictable compensator of $\mu$) $\nu(ds,dx)=\lambda(dx)\otimes ds$ for a given $\sigma$-finite measure $\lambda$. 
We consider discrete observation times $i/n,i=0,\ldots,n$, on the time span $[0,1]$. The prevalent model, capturing market microstructure effects which interfere the evolution of an underlying semimartingale log-price process at high frequencies, is an indirect observation model with noise:
\begin{align}\label{obs}Y_i=X_{i/n}+\epsilon_i\,,i=0,\ldots,n\,,\end{align}
where $(\epsilon_i)_{0\le i\le n}$ is a discretization of the continuous-time noise process $(U_t)_{t\in[0,1]}$. We consider $X$ and $U$ on a common probability space $(\Omega,\mathcal{F},(\mathcal{F}_t),\P)$ with $\mathcal{F}=\sigma\big(U_s,s\le 1\big)\bigvee \mathcal{F}^X$ and $\mathcal{F}_t=\sigma\big(U_s,s\le t\big)\bigvee\mathcal{F}_t^X$. Here, for two $\sigma$-algebras $\mathcal{F}$ and $\mathcal{H}$, we denote $\mathcal{F}\bigvee\mathcal{H}$ the smallest $\sigma$-algebra which contains $\mathcal{F}\cup\mathcal{H}$. %$\mathcal{F}_t=\bigcap_{s>t}\sigma\big(U_\tau,\tau\le s\big)\otimes\mathcal{F}_s^X$. 
$X$ has the same form \eqref{sm} on this space, see Section 16.1 of \cite{JP} for a formal construction of embedding $X$ and $U$ in a joint probability space. Regularity conditions on the characteristics of the efficient price $X$ and the noise, under which we establish asymptotic results, are given in Section \ref{sec:3.1}. In particular, we work with a general smoothness assumption on the volatility $(\sigma_t)_{t\in[0,1]}$. Similar to \cite{jacodtodorov}, resulting convergence rates of the spot volatility estimator and the asymptotic test hinge on this smoothness. First, readers may think of the typical case that $(\sigma_t)_{t\in[0,1]}$ is an It\^{o} semimartingale with a representation as $X$ in \eqref{sm} and with locally bounded characteristics.%\\[.2cm] 
\subsection{Test for common price and volatility jumps\label{sec:2.1}}
In the presence of price jumps, we design a statistical test to decide if contemporaneous price and volatility jumps occur on the considered time interval $[0,1]$. Let $(S_p)_{p\ge 1}$ be a sequence of stopping times exhausting the jumps of $X$. We denote the process of left limits of the volatility $\sigma_{t-}=\lim_{u\rightarrow t,u<t}\sigma_u$. We address the null hypothesis of no common jump of volatility and price on $[0,1]$:
\begin{align}{\mathds{H}}_{[0,1]}:~~\label{hypo1}\sum_{S_p\le 1}|\sigma^2_{S_p}-\sigma^2_{S_p-}|=0\,,\end{align}
against the alternative hypothesis that there is at least one jump in the volatility at a jump time of $X$.\\
Our test for \eqref{hypo1} relies on two main ingredients. First, localization of price jumps using thresholding. Second, a \emph{local} test for volatility jumps. Suppose we want to test $H_0^*:|\sigma^2_{s}-\sigma^2_{s-}|=0$ at a specific time $s\in(0,1)$, against the alternative hypothesis that the volatility exhibits a jump $|\sigma^2_{s}-\sigma^2_{s-}|>0$. For such a test we require estimates of the squared volatility at time $s$, $\hat \sigma^2_{s}$, and before time $s$, $\hat \sigma^2_{s-}$. An intuitive test statistic is the difference $\hat\sigma^2_{s}-\hat\sigma^2_{s-}$. It turns out that a more general class of statistics $T^*(s)=g\big(\hat \sigma^2_{s},\hat \sigma^2_{s-}\big)$ with a test function $g$ facilitates improved asymptotic properties.\\[.2cm]
If discrete observations of the efficient log-price $X_{i/n}\,,i=0,\ldots,n$, were directly available, and if we assume for this motivation that there are no jumps in $X$, $\sigma_{s}^2$ and $\sigma_{s-}^2$ could be estimated by local versions of realized volatility:
\begin{align}\label{mle}\hat\sigma_{s}^2=\frac{n}{k_n} \sum_{j=\lfloor s n\rfloor +1}^{\lfloor s n\rfloor +k_n} ( X_{(j+1)/n}-X_{j/n})^2~,~\hat\sigma_{s-}^2=\frac{n}{k_n} \sum_{j=\lfloor s n\rfloor -k_n}^{\lfloor s n\rfloor -1} ( X_{j/n}-X_{(j-1)/n})^2\,.\end{align}
For an It\^{o} semimartingale $(\sigma_t)_{t\in[0,1]}$, $k_n=c\,\sqrt{n}$ with some constant $c$, $\hat\sigma_s^2$ yields rate-optimal spot volatility estimators, that is, $(\hat\sigma_{s}^2-\sigma_{s}^2)=\mathcal{O}_{\P}\big(n^{-1/4}\big)$. Further, on the null hypothesis that $\sigma_{s-}=\sigma_{s}$, for $k_n=c\,n^{b}$ with $b=1/2-\delta$ and $\delta>0$ arbitrarily small, a stable central limit theorem can be proved
\begin{align*}n^{b/2}\big(\hat\sigma_{s}^2-\hat\sigma_{s-}^2\big)\stackrel{(st)}{\longrightarrow} MN\big(0,4\sigma_{s}^{4}\big)\,.\end{align*}
For stochastic volatility the limit is mixed normal and it is important that the convergence holds stably in law to allow for confidence intervals. This is a stronger mode of weak convergence which is equivalent to joint weak convergence with every $\mathcal{F}^X$-measurable bounded random variable, see \cite{JP} for an overview on stable limit theorems. This limit theorem readily supplies an asymptotic test for a volatility jump at time $s$ with a rate of convergence $n^{b/2}$. However, the convergence rate is rather slow and not optimal for this testing problem. For the test statistic
\begin{align}\label{JT}\mathcal{T}(s)=2\log\big(\tfrac12\big(\hat\sigma_{s}^2+\hat\sigma_{s-}^2\big)\big)-\log\big(\hat\sigma_{s}^2\big)-\log\big(\hat\sigma_{s-}^2\big)\end{align}
one derives instead $n^{b}\mathcal{T}(s)\stackrel{(st)}{\longrightarrow}\chi_1^2$ with a $\chi_1^2$ limit distribution and a much faster rate. This improves the (asymptotic) power significantly. A key property is that the test is pivotal, since $\mathcal{T}(s)$ is \emph{self-scaling} in the volatility. This means that it does not require some estimated asymptotic variance, since the limit does not depend on any unknown parameter. Such a local test is not separately highlighted in \cite{jacodtodorov}, but is contained as one ingredient of their general method. The final test statistic of \cite{jacodtodorov} is a sum of these local test statistics over all estimated jump times.\\[.2cm]
It is not obvious how to construct a generalization of the local test for a volatility jump to the noisy observations setup \eqref{obs}. Spot volatility estimators, which are local versions of integrated volatility estimators under noise, are available, see for instance \cite{zu14} and \cite{manc15}. For an It\^{o} semimartingale $(\sigma_t)_{t\in[0,1]}$ and i.i.d.\,noise with some moment assumption, stable central limit theorems
\begin{align*}n^{\beta/2}\big(\hat\sigma_{s}^2-\sigma_{s}^2\big)\stackrel{(st)}{\longrightarrow} MN\big(0,\textbf{AVAR}_{s}\big)\,\end{align*}
with optimal $\beta=1/4-\delta$, $\delta>0$, can be proved. Based on $\hat\sigma_{s}^2-\hat\sigma_{s-}^2$, a test with rate $n^{\beta/2}$ could be constructed. Asymptotic variances $\textbf{AVAR}_{s}$ of such estimators are usually sums of at least three addends: one depending on the noise variance, one including the quarticity $\sigma_{s}^4$ and a cross term depending on both. %The form of the variances is thus similar to the ones for integrated volatility estimators, see, for instance, \cite{bn2}, \cite{zhang} and \cite{JLMPV}. 
This applies to the asymptotic variances of the spot volatility estimators in \cite{zu14} and \cite{manc15}, which, however, have sub-optimal slower convergence rates localizing a sub-optimal two-scales integrated volatility estimator. The construction of a rate-optimal pre-average spot volatility estimator with an asymptotic variance of the type above is sketched in Section 8.7 of \cite{sahaliajacod}. Due to this structure of the asymptotic variance, it appears difficult to find a suitable test function that facilitates an asymptotic distribution free test with improved convergence rate.\\
Apart from attaining asymptotic efficiency, our main motivation to construct a method based on spectral spot volatility estimation is that we will be able to prove a stable central limit theorem
\begin{align*}n^{\beta/2}\big(\hat\sigma_{s}^2-\sigma_{s}^2\big)\stackrel{(st)}{\longrightarrow} MN\big(0,8\sigma_{s}^3\eta^{1/2}\big)\,\end{align*}
under mild assumptions for semimartingale volatility. Here, $\eta=\E[\epsilon_i^2]$ is the variance of i.i.d.\,noise, while we consider more general heteroscedastic and serially correlated noise in Section \ref{sec:3}. This enables us to find a suitable test function $g\big(\hat \sigma^2_{s},\hat \sigma^2_{s-}\big)$, such that
\begin{align}n^{\beta}\,T_0(s)~\stackrel{(st)}{\longrightarrow}\chi_1^2\,,\end{align}
for a test statistic $T_0(s)$ which is self-scaling in the volatility. The self-scaling property and the much faster convergence rate are key features to derive a reliable testing procedure.\\[.2cm]
To test the null hypothesis \eqref{hypo1}, local tests are performed at the estimated price-jump times which can be detected by truncation methods. Our asymptotic analysis provides results for the local test at some time $s$ as a special case.\\[.2cm]
%The method is of potential interest not only to test for contemporaneous price and volatility jumps. To evaluate the impact of news arrivals and information processing, for instance, economists might be interested to study volatility adjustments in response to certain firm specific or macroeconomic news announcements.\\[.2cm]
The tests for common price and volatility jumps of \cite{jacodtodorov} for direct observations and our generalization for noisy observations both restrict to finitely many large price adjustments at whose arrival times local tests are performed. Testing for volatility jumps over an interval instead would require a sequence of tests for volatility jumps at infinitely many points and is rather connected to a high-dimensional testing problem. A theory without noise recently has been presented in \cite{BJV} and a generalization of the techniques, which are quite different to \cite{jacodtodorov}, to the model with noise is a challenging topic for future research. It is clear that detecting volatility jumps from noisy observations of the price is especially difficult if we do not specify where to look for potential volatility jumps and the finite-sample performance of a global test is limited, see Section 6 of \cite{BJV}. Restricting to local tests for volatility jumps as in this work facilitates a larger power in finite-sample applications. 
\subsection{Spectral spot volatility estimators\label{sec:2.2}}
Consider a sequence of equispaced partitions of the considered time span $[0,1]$ into bins $[kh_n,(k+1)h_n),k=0,\ldots,h_n^{-1}-1$. For a simple notation suppose $nh_n\in\mathds{N}$, such that on each bin we enclose $nh_n$ noisy observations. A main idea of spectral volatility estimation, constructed in \cite{bibingerreiss}, is to perform optimal parametric estimation procedures localized on the bins. Based on these local estimates, one can build estimators for the spot and the integrated squared volatility. 
We utilize $L^2$-orthogonal functions $(\Phi_{jk})_{1\le j\le J_n}$ for spectral frequencies $1\le j\le J_n$ in the Fourier domain up to a spectral cut-off $J_n\le nh_n$. For $1\le j\le J_n,0\le k\le h_n^{-1}-1$ and $0\leq t\leq 1\,,$ we define
\begin{align}\label{Phi}
\Phi_{j0}(t)=\left(\sqrt{2 h_n}n\sin{\left(\frac{j\pi}{2nh_n}\right)}\right)^{-1}\sin\left(j\pi{h_{n}^{-1}}t\right)\1_{\left[0,h_{n}\right]}(t)\,,\,\Phi_{jk}\left(t\right)=\Phi_{j0}(t-kh_{n})\,.
\end{align}
The indicator functions localize the sine functions to the bins. For the spectral volatility estimation, local linear combinations of the noisy data are used with local weights obtained by evaluating the functions \eqref{Phi} on the discrete grid of observation times $i/n,i=0,\ldots,n$. %This strategy corresponds to performing a discrete sine transformation on the observed returns, similarly as proposed in \cite{corsi}, but localized over the bins. 
We use the notion of empirical scalar products and norms for functions $f,g$ as follows:
\begin{align}\label{sc}\langle f,g\rangle_n\:=\frac{1}{n}\sum_{l=1}^n f\left(\frac{l}{n}\right)g\left(\frac{l}{n}\right)~\mbox{and}~\|f\|_n^2\:=\frac{1}{n}\sum_{l=1}^nf^2\left(\frac{l}{n}\right)=\langle f,f\rangle_n\,.\end{align}
The empirical norms of the sine functions above give for all bins $ k=0,\ldots,h_n^{-1}-1$:
\begin{align}\label{empnorm}\|\Phi_{jk}\|_n^2=\left(4n^2\sin^2{(j\pi/(2nh_n))}\right)^{-1}\,,\end{align}
and we have the discrete orthogonality relations
\begin{align}\label{o}\langle \Phi_{jk},\Phi_{rk}\rangle_n=\|\Phi_{jk}\|_n^2\,\delta_{jr}\,,~j,r\in\{1,\ldots,J_n\}\,,k=0,\ldots,h_n^{-1}-1\,,\end{align}
where $\delta_{jr}=\1_{\{j=r\}}$ is Kronecker's delta. The latter rely on basic discrete Fourier analysis, a detailed proof is given in \cite{bibingerreiss}. The central building blocks of spectral volatility estimation are the spectral statistics
\begin{align}\label{spectralstatistics}S_{jk}=\|\Phi_{jk}\|_n^{-1}\sum_{i=1}^n\Delta_i^n Y\Phi_{jk}\Big(\frac{i}{n}\Big)~,j=1,\ldots,J_n,k=0,\ldots,h_n^{-1}-1\,,\end{align}
in which observed returns $\Delta_i^n Y=Y_{i/n}-Y_{(i-1)/n},i=1,\ldots,n$, are smoothed by bin-wise linear combinations. %with weights from the local discrete sine transformations.
Since the weight functions $\Phi_{jk}(t)$ %_{1\le j\le J_n,0\le k\le h_n^{-1}-1}$ 
are non-zero only on the $k$th bin, the spectral statistics $(S_{jk})$ include returns $(\Delta_i^n Y),i= k nh_n+1 ,\ldots,  (k+1) n h_n $ only over the bin under consideration. In absence of price jumps, bin-wise estimates for the squared volatility $\sigma^2_{kh_n},k=0,\ldots,h_n^{-1}-1$, are provided by weighted sums of bias-corrected squared spectral statistics:
\begin{align}\label{zeta}\zeta_k(Y)=\sum_{j=1}^{J_n} w_{jk}\Big(S_{jk}^2-\bk\Big)\,.\end{align}
For the moment, readers can interpret $(\eta_t)_{t\in[0,1]}$ as time varying variance function of the observation errors in \eqref{obs} and $\hat\eta_{kh_n}$ some consistent estimator. In Section \ref{sec:3.1}, this is further generalized. The oracle optimal weights
\begin{align}\label{orweights}w_{jk}=I_k^{-1}I_{jk}=\frac{\Big(\sigma^2_{kh_n}+\|\Phi_{jk}\|_n^{-2}\tfrac{\eta_{kh_n}}{n}\Big)^{-2}}{\sum_{m=1}^{J_n}\Big(\sigma^2_{kh_n}+\|\Phi_{mk}\|_n^{-2}\tfrac{\eta_{kh_n}}{n}\Big)^{-2}}\,,\end{align}
with $I_k=\sum_{j=1}^{J_n}I_{jk}, I_{jk}=\tfrac12 (\sigma^2_{kh_n}+\|\Phi_{jk}\|_n^{-2}\eta_{kh_n}/n)^{-2}$, follow from minimization of the variance under the constraint of unbiasedness. For a fully adaptive approach we apply a two-stage method and obtain adaptive local estimates $\zeta_k^{ad}(Y)$ by plugging in estimated optimal weights $\hat w_{jk}$ in \eqref{zeta}. 
%The integrated volatility estimator of \cite{BHMR} is simply the average \(\sum_{k=0}^{h_n^{-1}-1}h_n\) \(\zeta_k^{ad}(Y)\).
\begin{remark}
Spectral statistics are related to pre-averages used by \cite{JLMPV}, but the two estimators can not be transformed into one another, see Remark 5.2 in \cite{jacodmykland} for a discussion of their connection. One difference is that for the spectral method we start with a histogram structure and not a rolling kernel and then smooth bin-wise noisy observations in the Fourier domain. The statistics \eqref{spectralstatistics} de-correlate the data for different frequencies and form their local principal components. This is key to the asymptotic efficiency attained by the spectral estimators as shown in \cite{reiss} and \cite{BHMR}. The latter shows that the estimator's asymptotic variance coincides with the minimum asymptotic variance among all asymptotically unbiased estimators. We refer to Remark 3.1 of \cite{jacodmykland} for a recent discussion about efficient volatility estimation under noise.
%The efficiency theory is so far restricted to models with deterministic volatility, without drift and Gaussian noise and it is conjectured that the analogous asymptotic variance, which is determined in \cite{randolf} for the general model, constitutes the general lower bound.
\end{remark}
The spectral volatility estimation provides local estimates \eqref{zeta} for the squared volatility $\sigma^2_{kh_n},k=0,\ldots,h_n^{-1}-1$. 
%However, we should not rely on $\zeta_k^{ad}(Y)\1_{[kh_n,(k+1)h_n]}(t)$ on $[0,1]$ directly as an estimator for $(c_t)_{t\in[0,1]}$ because it will not be consistent as the variance does not decrease when $n$ gets large. Instead, we employ a smoothing method, as typical for function estimation in nonparametric statistics. 
In order to derive an estimate $\sigma^2_s$ at some time $s$, we average the statistics $\zeta_k(Y)$ over a local window around $s$ of length $(r_n^{-1}h_n)\rightarrow 0$ as $n\rightarrow\infty$, $r_n^{-1}\in\N$, slowly enough to ensure $r_n^{-1}\rightarrow\infty$. In the presence of jumps in \eqref{sm}, truncation disentangles bin-wise statistics \eqref{zeta} which include jumps from all others. We use the methods from \cite{bibwink2015} to cope with price jumps for volatility estimation. If $h_n|\zeta_k(Y)|>u_n$ for a threshold sequence $u_n=c\,h_n^{\tau},\tau\in(0,1)$, with some constant $c$, the statistic is too large to be driven by the continuous part and is evoked by a jump of $X$. In order to estimate the volatility, we thus truncate $\zeta_k(Y)$ for these $k$. For estimating the squared volatility and its left limit at a certain time $s$, we use two disjoint windows after and before $s$, respectively.\\ 
When the optimal weights \eqref{orweights} are known, an \emph{oracle spot volatility estimator} $\hat \sigma^2_{s,or}$ for $s\in[r_n^{-1}h_n,1-r_n^{-1}h_n)$ is:
\begin{subequations}
\begin{align}\label{rspotor}%\mbox{right-hand side:} ~~~~
\hat \sigma^2_{s,or}=\sum_{k=\lfloor sh_n^{-1}\rfloor+1}^{\lfloor sh_n^{-1}\rfloor +r_n^{-1}}r_n\sum_{j=1}^{J_n} w_{jk}\big(S_{jk}^2-\|\Phi_{jk}\|_n^{-2}\tfrac{\eta_{kh_n}}{n}\big)\1_{\{h_n|\zeta_k(Y)|\le u_n\}}\,,\end{align}
and the estimator for $\hat \sigma^2_{s-,or}$:
\begin{align}\label{lspotor}%\mbox{left-hand side:} ~~~
\hat \sigma^2_{s-,or}=\sum_{k=\lfloor sh_n^{-1}\rfloor -r_n^{-1}}^{\lfloor sh_n^{-1}\rfloor-1}r_n\sum_{j=1}^{J_n} w_{jk}\big(S_{jk}^2-\|\Phi_{jk}\|_n^{-2}\tfrac{\eta_{kh_n}}{n}\big)\1_{\{h_n|\zeta_k(Y)|\le u_n\}}\,.\end{align}
\end{subequations}
%\begin{subequations}
Close to the boundaries, $s\in[0,r_n^{-1}h_n)\cup (1-r_n^{-1}h_n,1]$, we shrink one window length accordingly. Since the optimal weights \eqref{orweights} hinge on the unknown squared volatility and the noise level $(\eta_t)_{t\in[0,1]}$, we proceed with a two-step estimation approach. First, select a pilot spectral cut-off $J_n^{pi}\ll nh_n$, and build \emph{pilot estimators} for the squared volatility
\begin{align}\label{pilot}\hat \sigma^2_{s,pil}=\sum_{k=\lfloor sh_n^{-1}\rfloor+1}^{\lfloor sh_n^{-1}\rfloor +r_n^{-1}}&r_n\sum_{j=1}^{J_n^{pi}} (J_n^{pi})^{-1}\big(S_{jk}^2-\|\Phi_{jk}\|_n^{-2}\tfrac{\hat\eta_{kh_n}}{n}\big)\\
&~~~~~~~~~~~~~~~~~~~~~\notag \times \1_{\Big\{h_n\big|\sum_{j=1}^{J_n^{pi}}(J_n^{pi})^{-1}\big(S_{jk}^2-\|\Phi_{jk}\|_n^{-2}\tfrac{\hat\eta_{kh_n}}{n}\big)\big|\le u_n\Big\}}\,,\end{align}
and $\hat \sigma^2_{s-,pil}$ analogously. The pilot estimators are hence averages of squared, bias-corrected spectral statistics %, bias-corrected with the estimated noise level $\hat\eta_{kh_n}$, 
over $r_n^{-1}$ bins and $J_n^{pi}$ spectral frequencies. %Estimation of $(\eta_t)_{t\in[0,1]}$ is addressed in Proposition \ref{cornoiseest}. 
In the second step, these pilot estimators are plugged into \eqref{orweights} to determine adaptive weights $\hat w_{jk}$ for the final estimators. We write 
\begin{align}\label{zetaad}\zeta_k^{ad}(Y)=\sum_{j=1}^{J_n} \hat w_{jk}\Big(S_{jk}^2-\|\Phi_{jk}\|_n^{-2}\tfrac{\hat\eta_{kh_n}}{n}\Big)\,.\end{align}
The \emph{spectral estimators of the squared spot volatility} at time $s$ and its left limit are:
\begin{subequations}
\begin{align}\label{rspot}%\mbox{right-hand side:} ~~~~
\hat \sigma^2_{s}=\sum_{k=\lfloor sh_n^{-1}\rfloor+1}^{\lfloor sh_n^{-1}\rfloor +r_n^{-1}}r_n\sum_{j=1}^{J_n} \hat w_{jk}\big(S_{jk}^2-\|\Phi_{jk}\|_n^{-2}\tfrac{\hat\eta_{kh_n}}{n}\big)\1_{\{h_n|\zeta_k^{ad}(Y)|\le u_n\}}\,,\end{align}
\begin{align}\label{lspot}%\mbox{left-hand side:} ~~~
\hat \sigma^2_{s-}=\sum_{k=\lfloor sh_n^{-1}\rfloor -r_n^{-1}}^{\lfloor sh_n^{-1}\rfloor-1}r_n\sum_{j=1}^{J_n} \hat w_{jk}\big(S_{jk}^2-\|\Phi_{jk}\|_n^{-2}\tfrac{\hat\eta_{kh_n}}{n}\big)\1_{\{h_n|\zeta_k^{ad}(Y)|\le u_n\}}\,.\end{align}
\end{subequations}
Estimates \eqref{rspot} and \eqref{lspot} are truncated local averages of the statistics \eqref{zetaad}. 
%Thus, this nonparametric spot volatility estimation is closely related to the usual nonparametric kernel estimation when the statistics \eqref{zetaad} take the role of de-noised observations which are smoothed over local windows. 
%This illuminates the relation to the nonparametric volatility estimator by \cite{kristensen} for the setup without microstructure noise. 
Our approach entails several tuning parameters whose practical choice is discussed in Section \ref{sec:4.2}. 

\section{Asymptotic theory\label{sec:3}}
\subsection{Assumptions with discussion\label{sec:3.1}}
We start with the assumptions on the characteristics of $X$ in \eqref{sm} which are similar to the ones in \cite{jacodtodorov}.
\begin{assump}\label{eff}
For the adapted and locally bounded drift process $(b_s)_{s\ge 0}$, we require a minimal smoothness condition that for $0\le t<s\le 1$, some constant $C$ and some $\iota>0$:
\begin{align}\label{drift}\E[(b_{s}-b_{t})^{2}|\mathcal{F}_t]\le C\,(s-t)^{\iota}~.\end{align}
The volatility process $\sigma_t$ is c\`{a}dl\`{a}g and neither $\sigma_t$ nor $\sigma_{t-}=\lim_{u\rightarrow t,u<t}\sigma_s$ vanish.
\end{assump}
\begin{strassump}\label{H}
We assume that $\sup_{\omega,x}|\delta(t,x)|/\gamma(x)$ is locally bounded for a non-negative deterministic function $\gamma$ satisfying $\int_{\mathds{R}}(\gamma^r(x)\wedge 1)\lambda(dx)<\infty$.
\end{strassump}
We index the assumption in $r\in[0,2]$ to highlight the role of the jump activity index $r$. The larger $r$, the more general jump components are included in our model. In particular for $r=0$ we consider jumps of finite activity. Imposing $r<1$ instead allows for infinite activity jumps which are absolutely summable.
We state the assumptions on characteristics of $X$ with respect to $(\Omega,\mathcal{F},(\mathcal{F}_t),\P)$, with the usual extension from $(\Omega^X,\mathcal{F}^X,(\mathcal{F}_t^X),\P^X)$. Especially, $(W_t)$ in \eqref{sm} is also a standard Brownian motion on this space. For the volatility process, our target of inference, we work with the following general smoothness condition determined by a smoothness parameter $\alpha\in(0,1]$.
\begin{sstrassump}\label{sigma}
The process $\sigma_t$ satisfies $\sigma_t=f_{\sigma}\big(\sigma_t^{(A)},\sigma_t^{(B)}\big)$ with some function $f_{\sigma}:\R^2\rightarrow\R$, continuously differentiable in both coordinates, and two $(\mathcal{F}_t)$-adapted processes $\sigma_t^{(A)},\sigma_t^{(B)}$, where
\begin{itemize}
\item $\sigma^{(A)}$ is an It\^{o} semimartingale
\begin{align}\label{smsigma}&\sigma_t^{(A)}=\sigma_0^{(A)}+\int_0^t\tilde b_s\,ds+\int_0^t\tilde \sigma_s\,dW_s +\int_0^t\tilde \sigma^*_s\,dW'_s\\
& \notag \quad\quad \quad +\int_0^t\int_{\mathds{R}}\tilde \delta(s,x)\1_{\{|\tilde \delta(s,x)|\le 1\}}(\tilde \mu-\tilde \nu)(ds,dx)
+\int_0^t\int_{\mathds{R}}\tilde \delta(s,x)\1_{\{|\tilde \delta(s,x)|> 1\}}\tilde \mu(ds,dx)\,,\end{align}
with an $(\mathcal{F}_t)$-Brownian motion $W'$ independent of $W$, locally bounded characteristics $\tilde\sigma,\tilde \sigma^*,\tilde b,\tilde \mu$ and a random variable $\sigma_0^{(A)}$. $\sigma_t^{(A)}$ satisfies Assumptions \ref{eff} and \hyperref[H]{(H-2)} for $\alpha\le 1/2$. For $\alpha>1/2$, the continuous martingale part of $\sigma^{(A)}$ vanishes and $\sigma^{(A)}$ satisfies Assumptions \ref{eff} and \hyperref[H]{(H-$\alpha^{-1}$)}.
\item $\sigma^{(B)}$ lies in a Hölder ball of order $\alpha$ almost surely, i.e.\,$\big|\sigma_t^{(B)}-\sigma_s^{(B)}\big|\le L |t-s|^{\alpha}$, for all $t,s\in[0,1]$ and a random variable $L$ for which at least fourth moments exist.
\end{itemize}
\end{sstrassump}
The smaller $\alpha$, the less restrictive is Assumption \ref{sigma}. It is natural to develop results for general $\alpha\in(0,1]$ to cover a broad framework and preserve some freedom in the model. This is particularly important, since the precision of nonparametrically estimating a process (or function) foremost hinges on its smoothness $\alpha$.
Therefore, convergence rates in Section \ref{sec:3.2} hinge on $\alpha$. %The composition of the volatility in Assumption \ref{sigma} allows to incorporate recent volatility models and to realistically describe spot volatility dynamics. 
In the composition of the volatility in Assumption \ref{sigma}, $\sigma_t^{(B)}$ can contain a non-Lipschitz seasonality component (Lipschitz continuous seasonalities can as well be modeled by the drift of $\sigma_t^{(A)}$). As pointed out by \cite{jacodtodorov}, $\sigma_t^{(B)}$ can also be a long-memory volatility component as the prominent exponential fractional Ornstein-Uhlenbeck model by \cite{comte}.\\
While an i.i.d.\;assumption on the noise is standard in most works, empirical findings, for instance by \cite{hans06}, motivate to allow for serial correlation and endogeneity in the noise. We develop our theory under the following general assumption.
\begin{nassump}\label{eta}
The noise $(\epsilon_i)_{0\le i\le n}$ process is centered, $\E[\epsilon_i]=0$. For some $p\ge 4$, its $\mathcal{F}^X$-conditional law has finite $p$-th moments, $\E\big[\epsilon_i^p|\mathcal{F}^X\big]<\infty$ almost surely for all $i=0,\ldots,n$. 
The long-run variance process converges
\begin{align}\label{lrvar}\sum_{l=-\lfloor t n\rfloor }^{n-\lfloor t n\rfloor}\cov\big(\epsilon_{\lfloor t n\rfloor },\epsilon_{\lfloor t n\rfloor +l}\big)\rightarrow \eta_t~,\end{align}
for $t\in[0,1]$ uniformly on compacts in probability and we have the mixing behavior 
\begin{align}\label{varrho}\sup_{i=0,\ldots, n}\big|\cov(\epsilon_i,\epsilon_{i+l})\big|=\mathcal{O}\big(|l|^{-1-\varrho}\big)\,,\end{align} 
for some $\varrho>0$, which is specified in the discussion below Theorem \ref{cltspot}. The process $(\eta_t)_{t\in[0,1]}$ is locally bounded and satisfies for all $t,(t+s)\in[0,1]$ the mild smoothness condition:
\begin{align} \label{etasmooth}|\eta_{t+s}-\eta_t|\le K s^{(1/2+\delta) \vee \alpha}\,,\end{align}
with some $\delta>0$. The noise $\epsilon_{i}$ is for all $i$ uncorrelated to $(\Delta_l^n X)_{l=1,\ldots,(i-\tilde Q-1\vee 1)}$ for some $\tilde Q<\infty$ and
\vspace*{-1cm}

\begin{align}\label{endo}\sum_{l={\lfloor t n\rfloor }-\tilde Q}^{{\lfloor t n\rfloor }}\E\big[\epsilon_{\lfloor t n\rfloor }\Delta_l^n X\big]\rightarrow \rho_t~,\end{align}
for some continuous bounded function $(\rho_t)_{t\in[0,1]}$. Furthermore, the noise does not vanish, $\eta_t>0$ for all $t\in[0,1]$. %When $(\eta_t)_{t\in[0,1]}$ is stochastic, for notational convenience, we augment the probability space such that $(\eta_t)_{t\in[0,1]}$ is $\mathcal{F}_t$-adapted.

\end{nassump}
The case that $\cov(\epsilon_i,\epsilon_{i+l})=0$ for all $l\ne 0$ and $\eta=\var(\epsilon_i)$ constant for all $i$ is tantamount
to the classical setup with i.i.d.\;noise. In general the noise is serially correlated, endogenous and heteroscedastic. Different to Assumption (GN) in Section 7.2 of \cite{sahaliajacod}, we do not assume that the noise is conditionally centered to include the correlation to the increments of $X$ in \eqref{endo}. The endogeneity condition \eqref{endo} includes linear models of the form $\epsilon_i=\sum_{l=i-\tilde Q}^{i}c_l\sqrt{n}\Delta_l^n X+U_i$, with $U_i$ exogenous errors and constants $c_l$, similar as in Equation (6) of \citet{koike2016} or considered by \cite{bn2}. If we knew the process $(\eta_t)_{t\in[0,1]}$, Assumption \ref{eta} with a mild lower bound for $\varrho$ would be sufficient for our asymptotic results. For an adaptive method, however, we need to estimate the process $(\eta_{t})_{t\in[0,1]}$. Consistent estimation of the noise long-run variance \eqref{lrvar} requires stronger structural assumptions. For a $Q$-dependent noise process, that is, $\sup_{i=0,\ldots, n}|\cov(\epsilon_i,\epsilon_{i+q})|=0$ for $q>Q$ and some given $Q<\infty$, and if $\eta$ in \eqref{lrvar} is time-invariant, consistent estimation with $\sqrt{n}$-convergence rate of $\eta$ has been established by \cite{haut13}. \cite{BHMR2} show how $Q$ can be found adaptively if it is unknown. \cite{jacodmykland} discuss consistent estimation of the noise variance process under heteroscedasticity, but without serial correlations. %For $\alpha>1/2$ in Assumption \ref{sigma}, we impose a stronger smoothness of $(\eta_t)_{t\ge 0}$ in \eqref{etasmooth}, such that roughness of the long-run noise variance process can not manipulate the resulting convergence rates. We formulate \eqref{etasmooth} as general as possible while from an applied point of view a smoother long-run noise variance process appears realistic. 
For the fully adaptive method, we tighten the assumptions on the noise as follows.
\begin{assump}\label{eta2}Assumption \ref{eta} holds with $p\ge 8$. Moreover, 
\begin{align*}\sup_{i=0,\ldots, n}\big|\cov(\epsilon_i,\epsilon_{i+q})\big|=0\end{align*}
for all $q>Q$ with some $Q<\infty$. 
\end{assump}
Assumption \ref{eta2} is satisfied by a $Q$-dependent noise process. Then, a consistent estimation of the long-run noise variance \eqref{lrvar} process is possible.
\begin{prop}\label{cornoiseest}
Under Assumption \ref{eta2}, for $h_n=\kappa_1 n^{-1/2}\log(n)$, for all $k=0,\ldots,h_n^{-1}-1$, the locally constant approximated noise long-run variance process can be estimated with accuracy
\begin{align}\label{noiseest}\hat\eta_{kh_n}=\eta_{kh_n}+\KLEINO_{\P}\big(n^{-\beta}\big)\,.\end{align}
\end{prop}
Our estimator is given in \eqref{noiseesteq} in the appendix. It is somewhat related to the methods from \cite{haut13} and \cite{BHMR2}, but localized to bins.\\[.2cm]
%When we apply the global method of \cite{BHMR2} localized to bins $[k h_n, (k+1) h_n)$, such that $\eta_{kh_n}$ is estimated in the same way as in Equations (19a)-(19c) of \cite{BHMR2}, but using observations $i= k h_n n +1,\ldots, (k+1) h_n n$ only instead of all observations $i=0,\ldots,n$, the regularity \eqref{etasmooth} renders for $h_n=\kappa_1 n^{-1/2}\log(n)$ under Assumption \ref{eta2} such estimators. For fixed $k\in\{0,\ldots,h_n^{-1}-1\}$, and locally constant $\eta_{kh_n}$, the proof from \cite{BHMR2} can be adopted just using $nh_n$ observations on $[k h_n, (k+1) h_n)$ instead of all $n$ observations. This results in the slower rate $n^{1/4}$ instead of $n^{1/2}$. Regularity \eqref{etasmooth} ensures that the approximation error of setting $\eta_{kh_n}$ locally constant on $[k h_n, (k+1) h_n)$ is asymptotically negligible. This readily gives Proposition \ref{cornoiseest}.\\[.2cm]
The assumptions on the noise are more general than in other works on spectral volatility estimation as in \cite{stable} and in \cite{BHMR2}. In particular, to the best of our knowledge, we consider for the first time heteroscedastic and serially correlated, endogenous noise. %We exclude, however, endogeneity, see \citet{koike2016} for a recent study of endogenous noise.
\begin{remark}[Non-equidistant observations]
For a coherent and simple exposition of the construction of the spectral estimator in \eqref{Phi}-\eqref{spectralstatistics}, we discuss equidistant observations which allows us to rely on discrete-time Fourier identities in \eqref{o}. Considering a heteroscedastic noise-level, our analysis and results are at the same time informative about non-equidistant observations. For general observation schemes $t_i^n,i=0,\ldots,n$, we impose the condition that a differentiable cdf $F$ exists such that observation times $t_i^n=F^{-1}(i/n)$ are obtained by a quantile transformation from the equidistant setting. Moreover, we require that the derivative $F^{\prime}$ is strictly positive and satisfies the same smoothness as $(\eta_t)$ in \eqref{etasmooth}. These assumptions are the same as in Assumption (Obs-d) of \cite{stable}. Then, all our asymptotic results transfer from the equidistant to this general setting when we replace $\eta_s$ by $\eta_s\,(F^{-1})^{\prime}(s)$. This follows directly by the asymptotic equivalence of the respective experiments established in \cite{BHMR}. In particular, having locally less frequent observations is equivalent to having locally an increased noise level. Therefore, under the imposed conditions, $(\eta_t)_{t\in[0,1]}$ and $(F^{-1})^{\prime}(t), t\in[0,1],$ may be pooled. Note that adding the factor $(F^{-1})^{\prime}(s)$ to the noise level $\eta_s$ is the same as generalizing the frequently occurring factor $\|\Phi_{jk}\|_n^{-2}\eta_{kh_n}/n$ to $\|\Phi_{jk}\|_n^{-2}\eta_{kh_n}/(nF^{\prime}(kh_n))$, where $nh_nF^{\prime}(kh_n)$ gives the local sample size. In the equidistant case this is $nh_n$ and we have that $F^{\prime}(s)=1$ is constant.
\end{remark}
%\newpage
\subsection{Asymptotic results\label{sec:3.2}}
Our first main result is on the spot squared volatility estimator and its asymptotic distribution.
\begin{theo}\label{cltspot}Suppose Assumptions \ref{eff}, \ref{eta2} and \ref{H} with some $r<2$ and smoothness Assumption \ref{sigma}, $\alpha\in(0,1]$. Fix some time $s\in(0,1)$, at which we want to estimate $\sigma^2_s$ and $\sigma^2_{s-}$ with \eqref{rspot} and \eqref{lspot}, respectively. Set $h_n=\kappa_1 n^{-1/2}\log(n)$ and $r_n=\kappa_2 n^{-\beta}\log(n)$ with constants $\kappa_1,\kappa_2$ and $J_n\rightarrow\infty$, $J_n=\mathcal{O}(\log(n))$, as $n\rightarrow\infty$. Then, as $n\rightarrow\infty$ and if
\begin{align}\label{beta}0<\beta<\left(\frac{\alpha}{2\alpha+1}\,\wedge\,\tau\Big(1-\frac{r}{2}\Big)\right)\,, \end{align}
and $\tau<1-\beta/(p-2)$ when $p<\infty$ moments of the noise exist, with $\tau$ the truncation exponent in the sequence $u_n$ in \eqref{pilot}, \eqref{rspot} and \eqref{lspot}, the estimators satisfy the $\mathcal{F}$-stable central limit theorem:
%\begin{subequations}
\begin{align}\label{cltspoteq}n^{\beta/2}\left(\begin{matrix}\hat \sigma^2_s-\sigma^2_s\\ \hat \sigma^2_{s-}-\sigma^2_{s-}\end{matrix}\right)\stackrel{(st)}{\longrightarrow} MN\left(0,\left(\begin{matrix}8\sigma_s^{3}\eta_s^{1/2}&0\\0&8\sigma_{s-}^{3}\eta_s^{1/2}\end{matrix}\right)\right)\,.\end{align}
%\begin{align}\label{cltlspot}n^{\nicefrac{\beta}{2}}\big(\hat \sigma^2_{s-}-\sigma^2_{s-}\big)\stackrel{(st)}{\longrightarrow} MN\big(0,8\sigma_{s-}^{3}\eta_s^{\nicefrac12}\big)\,.\end{align}
%\end{subequations}
\end{theo}
For the oracle estimators \eqref{rspotor} and \eqref{lspotor}, the same limit theorem applies under the less restrictive Assumption \ref{eta} with $p=8$, $\varrho>\beta$, and if $\tau<1-\beta/(p-2)$. In fact, we can get arbitrarily close to the optimal rate for estimation which is known to be $n^{\alpha/(4\alpha+2)}$ in this case, see \cite{munk2010}. Balancing the squared bias and the variance guarantees that the estimators \eqref{rspot} and \eqref{lspot} attain the optimal rate. For a central limit theorem we avoid an asymptotic bias by slightly undersmoothing. Most interesting is the case when $\alpha\approx 1/2$, e.g.\ when the volatility is a semimartingale. Then the convergence rate is $n^{1/8}$. In case that $\alpha>1/2$, we obtain faster convergence rates. In case that $\alpha= 1/2$ and if all moments of the noise process exist, for any $r<3/2$ in Assumption \ref{H}, we can choose $\beta=1/4-\varepsilon$ with any $\varepsilon>0$. Under the standard assumption that we only have Assumption \ref{eta} with $p=8$, the condition $\tau<23/24$ results in $r<34/23\approx 1.478$. Hence, restricting to the condition that up to 8th moments of the noise exist leads only to a slightly less general condition on the jump activity. We point out that the restriction $r<3/2$ on the jump activity, to come close to the optimal convergence rate, is less restrictive than the one obtained for integrated squared volatility estimation, $r<1$, in \cite{bibwink2015}. The reason is that for spot volatility estimation we can only obtain slower convergence rates by local smoothing compared to integrated volatility estimation. This, however, works also under more active jumps.\\ The limit variable in \eqref{cltspoteq} is \textit{mixed normal} which we denote by $MN$ and defined on a product space of the original probability space (on which $Y$ is defined) and an orthogonal space independent of $\mathcal{F}$. The convergence is $\mathcal{F}$-stable in law, marked $(st)$. Stability of weak convergence then allows for a so-called \textit{feasible} version of the limit theorem \eqref{cltspoteq} that facilitates confidence sets.
\begin{cor}Under the conditions of Theorem \ref{cltspot}, and also for any $J_n$ fixed as $n\rightarrow\infty$:
\begin{align}\label{feasclt}r_n^{-1/2}\left(\begin{matrix}\hat I_{\lfloor s h_n^{-1}\rfloor+1}^{1/2}(\hat \sigma^2_s-\sigma^2_s)\\[.2cm] \hat I_{\lfloor s h_n^{-1}\rfloor-1}^{1/2}(\hat \sigma^2_{s-}-\sigma^2_{s-})\end{matrix}\right)\stackrel{(st)}{\longrightarrow} N\left(0,\left(\begin{matrix}1&0\\0&1\end{matrix}\right)\right)\,,\end{align}
%\begin{align}\label{feasclt}r_n^{-1/2}\hat I_{\lfloor s h_n^{-1}\rfloor+1}^{1/2}\big(\hat \sigma^2_{s}-\sigma^2_{s}\big)\stackrel{d}{\longrightarrow}N(0,1)\,,\end{align}
with $\hat I_{\lfloor s h_n^{-1}\rfloor +1}$ and $\hat I_{\lfloor s h_n^{-1}\rfloor-1}^{1/2}$ the estimates of $I_{\lfloor s h_n^{-1}\rfloor+1}$ and $I_{\lfloor s h_n^{-1}\rfloor-1}$, as defined in the weights \eqref{orweights}, obtained by inserting the pilot estimates. %This works analogously for $\hat\sigma_{s-}^2$ for which we normalize with $\hat I_{\lfloor s h_n^{-1}\rfloor-1}^{1/2}$ instead.
\end{cor}
The results proved for the spot volatility estimator provide a main building block for our asymptotic test, but are moreover of interest in their own right. They show that the spectral method renders effective spot squared volatility estimators under general noise and in the presence of jumps.\\[.2cm]
In the sequel, let $(S_p)_{p\ge 1}$ be a sequence of stopping times exhausting the jumps of $X$. We address the null hypothesis \eqref{hypo1} that no common jumps of volatility and price occur on $[0,1]$. Under the alternative hypothesis, there is at least one contemporaneous jump in volatility and price.\\
Analogously to \cite{jacodtodorov}, we specify test hypotheses more precisely by focusing on jumps of $X$ with absolute values $|\Delta X_{S_p}|>a$ for $a\ge 0$ and write ${\mathds{H}(a)}_{[0,1]}$. The reason for this is that a suitable test statistic and associated limit theory for ${\mathds{H}(a)}_{[0,1]}$ with $a>0$ works under a much more general setup with jumps of infinite variation while testing ${\mathds{H}(0)}_{[0,1]}$ requires Assumption \hyperref[H]{(H-0)} to hold. In both cases, we concentrate on a finite number of (large) price jumps under the null hypothesis. From an applied point of view this is reasonable, since we are interested in volatility movements at finitely many relevant price adjustments on a fixed time interval.\\[.2cm]
Denote by $g:\mathds{R}_+^2\rightarrow \mathds{R}$ a test function with $g(x,x)=0$ for all $x$. Let us now state the general form of our test statistics:
\begin{align}\label{teststatistic}T_0(h_n,r_n,g)=\sum_{k=r_n^{-1}}^{h_n^{-1}-r_n^{-1}-1}\hat\eta_{kh_n}^{-1/2}g\big(\hat \sigma^2_{kh_n},\hat \sigma^2_{kh_n-}\big)\1_{\big\{h_n|\zeta_k^{ad}(Y)|>\,(u_n\vee a^2)\big\}}\,.\end{align}
Under mild regularity assumptions on $g$ in terms of differentiability in both coordinates, limit theorems for \eqref{teststatistic} can be proved. For testing ${\mathds{H}(a)}_{[0,1]}$, we consider two specific test functions in the following. Adjustments of the test \eqref{hypo1} for sub-intervals of $[0,1]$ are readily obtained by ignoring all jumps elsewhere.
\begin{theo}\label{test}Let $S_1,\ldots,S_{N_1}$ be a finite collection of jump times of $X$ on $[0,1],$ with $|\Delta X_{S_i}|>a$ for all $i$. Consider ${\mathds{H}}(a)_{[0,1]}$, if either $a>0$ and we impose the condition that the L\'{e}vy measure of $X$ does not have an atom in $\{a\}$, or assume $r=0$. On all assumptions of Theorem \ref{cltspot} and if $\tau<3/4$ for $a=0$, when inserting estimates \eqref{rspot} and \eqref{lspot} with $h_n=\kappa_1n^{-1/2}\log{(n)}$, $r_n=\kappa_2n^{-\beta}\log{(n)}$, %, where $0<\beta<\big(\alpha(2\alpha+1)^{-1}\wedge \tau(1-r/2)\big)$, $J_n\rightarrow\infty$, 
$J_n\rightarrow\infty$, $J_n=\mathcal{O}(\log(n))$ in \eqref{teststatistic} with the test function
\begin{align}\label{stang}g(x_1,x_2)=2\sqrt{\frac{x_1+x_2}{2}}-\sqrt{x_1}-\sqrt{x_2}\,,\end{align}
the following asymptotic distribution of the test statistic applies under ${\mathds{H}}(a)_{[0,1]}$:
\begin{align}\label{testlt}n^{\beta}\,T_0(h_n,r_n,g)\stackrel{(st)}{\longrightarrow} \chi^2_{N_1}\,.\end{align}
Under the alternative hypothesis, $n^{\beta}\,T_0(h_n,r_n,g)\rightarrow\infty$ in probability. Therefore, we obtain an asymptotic distribution free test 
%\footnote{That is, the asymptotic distribution does not hinge on any unknown parameters.} 
by the asymptotic $\chi^2$-distribution with $N_1$ degrees of freedom. The test with critical regions
\begin{align}\label{critc}C_n=\big\{n^{\beta}\,T_0(h_n,r_n,g)>q_{1-\alpha}(\chi^2_{\hat N_1})\}\,,\end{align}
where $q_{\alpha}(\chi^2_{\hat N_1})$ denotes the $\alpha$-quantile of the $\chi_{\hat N_1}^2$-distribution, has asymptotic level $\alpha$ and asymptotic power 1.
\end{theo}
In fact, \eqref{critc} contains the estimated number of price jumps $\hat N_1$. Since $\P(\hat N_1-N_1>0)\rightarrow 0$, \eqref{testlt} applies with $N_1$ also.
A naive approach based on the asymptotic normality result \eqref{feasclt} with test function $\tilde g(x_1,x_2)=(x_1-x_2)$ yields as well an asymptotic test. It holds that
\begin{align}\label{naive}r_n^{-1/2}\Bigg(2\sum_{i=1}^{\hat N_1}\hat I_{\lfloor h_n^{-1}S_i\rfloor+1}^{-1}\Bigg)^{-1/2}T_0(h_n,r_n,\tilde g)\stackrel{d}{\longrightarrow} N(0,1)\,,\end{align}
on the null hypothesis $\mathds{H}(a)_{[0,1]}$. Apparently, the rate $r_n^{-1/2}\asymp n^{\beta/2}$,\footnote{We write $a_n\asymp b_n$ for asymptotically equivalent real sequences which means $a_n/b_n\to c$ for some constant $c$.} close to $n^{1/8}$ for $\alpha\le 1/2$, is slower and thus the test in Theorem \ref{test} is preferable. 
\begin{remark}As mentioned by \cite{jacodtodorov}, their test based on \eqref{JT} corresponds to a two-sample likelihood ratio test for equal variances in a Gaussian parametric model with observations $\sqrt{n}\Delta_j^n X\stackrel{iid}{\sim}N(0,\sigma_{s-}^2)$, $j=\lfloor sn\rfloor-k_n,\lfloor sn\rfloor-1$ and $\sqrt{n}\Delta_j^n X\stackrel{iid}{\sim}N(0,\sigma_{s}^2)$, $j=\lfloor sn\rfloor+1,\lfloor sn\rfloor+k_n$. In this simpler model -- closely related to our model in case of no noise -- the likelihood ratio is
\[\Lambda=\frac{(\hat\sigma_s^2\hat\sigma_{s-}^2)^{k_n/2}}{\Big(\frac{\hat\sigma_s^2+\hat\sigma_{s-}^2}{2}\Big)^{k_n}}, ~\mbox{and}~-2\log(\Lambda)=k_n\Big(2\log \tfrac{\hat\sigma_s^2+\hat\sigma_{s-}^2}{2}-\log \hat\sigma_s^2-\log \hat\sigma_{s-}^2\Big)\,,\]
where the estimators \eqref{mle} are the maximum likelihood estimators for this model, and we derive the convergence of $-2\log(\Lambda)$ to a $\chi_1^2$-distribution from the standard asymptotic theory for likelihood ratio tests.\\
The model with noise is more complicated. Our test from Theorem \ref{test} does not directly correspond to a parametric likelihood ratio test and our estimators \eqref{rspot} and \eqref{lspot} do also not agree with the non-explicit maximum likelihood estimators in this model. The choice of $g$ in \eqref{stang} is motivated by studying which properties in \eqref{JT} are important for the asymptotic pivotal distribution under the null.
Any function of the form $g(x,y)=2f(\tfrac{x+y}{2})-f(x)-f(y)$, with some twice differentiable function $f$, is suitable for the construction of tests (in both models) with the fast convergence rate based on second order asymptotics of the estimators, since $\frac{d}{dx}g(\sigma_s^2,\sigma_{s}^2)=\frac{d}{dy}g(\sigma_s^2,\sigma_{s}^2)=0$. On the other hand, that the statistic \eqref{JT} is self-scaling in the volatility leading to the pivotal limit distribution is due to the identity
\[\frac{d^2}{dx^2}g(\sigma_s^2,\sigma_{s}^2)=\frac{d^2}{dy^2}g(\sigma_s^2,\sigma_{s}^2)=-\frac{1}{2}f^{\prime\prime}(\sigma_s^2)=\big(n^{b}\,\var(\hat\sigma_s^2)\big)^{-1}\]
denoting $f^{\prime\prime}$ the second derivative of $f$. With $f(x)=\log(x)$, it holds that $-\frac{1}{2}f^{\prime\prime}(\sigma_s^2)=(2\sigma_s^4)^{-1}$, which guarantees the above identity in the model without noise. In light of the efficient asymptotic variance under noise in Theorem \ref{cltspot}, it is natural to choose $f(x)=\sqrt{x}$, such that 
\[-\frac{1}{2}f^{\prime\prime}(\sigma_s^2)=\frac{1}{8\sigma_s^3}=\frac{\eta_s^{1/2}}{n^{\beta}\var(\hat \sigma_s^2)}\,.\]
Since the noise level $\eta_s$ can be estimated with a much faster rate of convergence than $\sigma_s^2$ -- even under our general assumptions for the noise -- this choice of \eqref{stang} facilitates \eqref{testlt}. \\
The particular choice of the spectral estimators \eqref{rspot} and \eqref{lspot} is not crucial for the construction of the test. Any rate-optimal spot volatility estimator may be used when it is possible to find a function $f$ satisfying the above identities. However, with a more complex asymptotic variance structure, for instance for pre-average or realized kernel estimators, this appears to be difficult. Estimators attaining the same efficient variance as in \eqref{cltspoteq} may be used with the same function $g$ in \eqref{stang}, to derive a test with the same asymptotic properties. A localized QMLE as discussed by \cite{clinet1}, for instance, could allow for analogous results.
\end{remark}
\section{Implementation and numerical study\label{sec:4}}

\subsection{Setup of Monte Carlo simulation study\label{sec:4.1}}
The simulation study examines the finite-sample performance of the proposed methods. We implement a model where observed log-prices are given by
\[Y_{i/n}=\int_0^{\tfrac{i}{n}}\varphi_t\,\sigma_t \,dW_t+\int_0^{\tfrac{i}{n}}\int_{\R}x\mu(dt,dx,dy)+\epsilon_i\,,\]\\
with jump intensity measure $\nu(dt,dx,dy)=\lambda\, dt\,\Pi(dx)\Pi(dy)$ and with Gaussian jump sizes $\Pi\sim N(H,$ $H/100)$ whose magnitude depend on a parameter $H$. The efficient semimartingale log-price process is recorded with additive microstructure noise
\begin{align}\label{NS1}\epsilon_i=\theta\epsilon_{i-1}+u_i\,,~u_i\stackrel{iid}{\sim} N\Big(0,\eta\big(1-\theta\big)^{-2}\Big),~i=1,\ldots,n\,, |\theta|<1\,.\end{align} 
In line with empirical evidence, this model generates serially correlated noise. We further consider two different noise models \eqref{NS2} and \eqref{NS3} below. We set values of $\eta$ according to realistic noise-to-signal ratios. We use the median value of the estimated measure $n\eta\big(\int_0^1\varphi_t^4 \sigma_t^4)^{-1/2}$ found in a comprehensive data study in \cite{BHMR2}. Sample sizes $n=30,000$ and $n=5,000$ in our simulations suggest $\eta^{1/2}\approx 0.005$ and $\eta^{1/2}\approx 0.015$, which we use in the following as two realistic noise levels. According to the data summary in Table \ref{tab4}, 30,000 is a sample size that matches (approximately) the average daily observation numbers of our empirical data. We additionally analyze the methods' performance for smaller samples sizes $n=5,000$, which is realistic for less frequently traded assets. We set $\theta=0.6$ equal to the empirically motivated value in \cite{BHMR2}.\\
\(\varphi_t=1-\tfrac{3}{5}\sqrt{t}+\tfrac{1}{10}t^2\) mimics a deterministic volatility intra-day seasonality pattern
and $\sigma_t^2$ a random stochastic volatility component with leverage:
\[d\sigma^2_t=6(1-\sigma^2_t)\,dt+\sigma^2_t dB_t+dJ_t\,.\]
$B$ is a standard Brownian motion with $d[B,W]_t=\rho\,dt$, where we fix $\rho=0.2$.\\
The jump measure above has a second real argument to incorporate instantaneous arrivals of volatility jumps. The volatility jump component is of the form
\[J_t=\gamma\,\int_0^t\int_{\R}y\mu(dt,dx,dy)+\,\int_0^t\int_{\R}z\tilde \mu(dt,dz)\]
with $\gamma\in \R$ and intensity measure $\tilde \nu(dt,dz)=dt\,\Pi(dz)$. Setting $\gamma=0$ results in no common price and volatility jumps which means the null hypothesis is valid. To simulate the model under the alternative hypothesis, we set $\gamma=1$ instead.

\subsection{Choice of tuning parameters\label{sec:4.2}}
\begin{figure}[t]
\includegraphics[width=7.45cm]{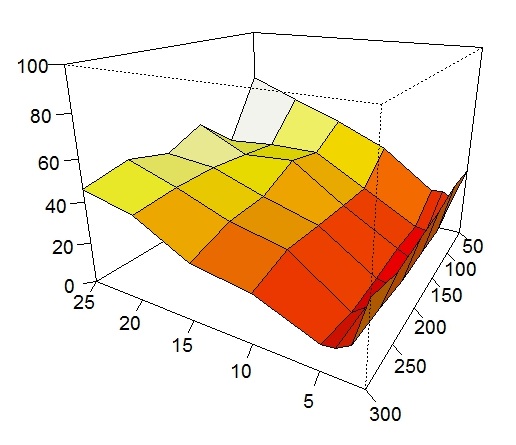}\includegraphics[width=7.45cm]{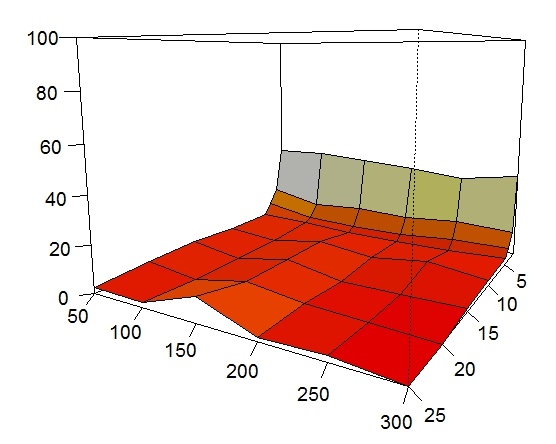}
\caption{\label{Fig:2}Empirical percentage type-II-error rate (right) and empirical percentage global testing error rate (left) for the test of size $\alpha=0.05$, depending on tuning parameters $h_n^{-1}$ and $r_n^{-1}$, with $50\le h_{30,000}^{-1}\le 300$ and $2\le r_{30,000}^{-1}\le 25$. The empirical type-II-error rate measures the empirical amount of realizations  under the alternative hypothesis which are smaller or equal the $.95$-quantile of the $\chi^2_{N_1}$-distribution. The global testing error rate is the sum of the type-II-error rate and the misspecification of the size, that is, the difference between $(1-\alpha)$ and the empirical amount of realizations smaller or equal the $(1-\alpha)$-quantile of the $\chi^2_{N_1}$-distribution, this time on the null hypothesis.}
\end{figure}
%\begin{figure}[t]
%\includegraphics[width=7.45cm]{persp90_5000}\includegraphics[width=7.45cm]{persp95_5000}
%\caption{\label{Fig:3}Empirical global testing error rates for the tests of size $\alpha=0.05$ (right) and $\alpha=0.1$ (left) depending on tuning parameters $h_n^{-1}$ and $r_n^{-1}$, with $20\le h_{5,000}^{-1}\le 100$ and $2\le r_{5,000}^{-1}\le 10$. The global testing error rate is defined as in Figure \ref{Fig:2}.}
%\end{figure}
In the sequel, we provide advice on how to specify the tuning parameters that are involved in the nonparametric procedures. We also conduct a sensitivity analysis for the Monte Carlo study to find suitable values.\\
First, the bin-width $h_n\asymp n^{-1/2}\log{n}$ balances the number of observations on bins $nh_n$, which should be large enough to smooth out noise, and the discretization error by approximating volatility bin-wise constant. %The smoother the underlying volatility process, the smaller the discretization error which allows to consider larger bins. On the other hand, the smaller the noise level, the smaller the bins can be chosen. The bin-width also determines the limits in disentangling small jumps from continuous motion. Smaller bins allow to detect smaller jumps. This applies to price jumps and also to volatility jumps. 
The sensitivity analysis will show that the final test is very robust to modifications of $h_n$. We advise to select $h_n$ such that the number of observations on bins is at least 50 within a range to 250 observations for typical high-frequency financial data. This results in a time resolution of 50-150 bins per trading day.\\
For the spot volatility estimators \eqref{rspot} and \eqref{lspot} and the pilot estimator \eqref{pilot}, we fix spectral cut-offs $J_n$ and $J_n^{pi}$, respectively. The values of the spectral cut-offs do not influence the methods when set sufficiently large. Since the weights \eqref{orweights} decay exponentially for $j\gtrsim \sqrt{n}h_n\asymp\log{n}$, the addends with $j$ large become negligible, such that it suffices to choose $J_n\asymp \log{n}$. The proportionality constant should be larger than 1, we take values between 3 and 12. The pilot estimators \eqref{pilot} instead use averages over frequencies $j=1,\ldots,J_n^{pi}$, such that we fix $J_n^{pi}$ to be smaller. We thus use $J_n^{pi}\asymp \log{n}$ with a proportionality factor smaller than for $J_n$. The threshold sequence $u_n$ determines the bins on which large returns are ascribed to
jumps. We use the practical selection presented in \cite{bibwink2015}.\\ %For the asymptotic theory $u_n\asymp h_n^{\tau}$ works with any $\tau\in(0,1)$. 
%Since in absence of
%jumps local estimates $h_n\zeta_k^{ad}(Y)$ are of order $h_n$ and the maximum over all
%bins at most of order $2\log{(h_n^{-1})}h_n$, a simple global truncation rule is to set $u_n=h_n2\log{(h_n^{-1})}$. This threshold is used for the pre-estimation step of our two-stage method.
%Since the expectation of $\zeta_k^{ad}(Y)$ is $\sigma^2_{kh_n}$ in absence of jumps, we employ in the second step a time-varying adaptive truncation where the first constant threshold is multiplied with the pre-estimated local squared volatility.\\
The most influential tuning parameter for our test is the size of the smoothing window $r_n\asymp n^{-\beta}\log{n}$. If we choose $r_n$ larger, the spot volatility estimates have smaller variance but the bias for rapidly varying volatilities increases. For $\alpha=1/2$, we know the exact order of $r_n$ depending on $n$. There is, however, no simple rule of thumb to fix the constant $\kappa_2$, and we conduct an extensive sensitivity analysis to find the best suitable values. The sensitivity analysis reveals that in order to detect volatility jumps and separate them from a rough continuous semimartingale volatility component, we should use rather small smoothing window sizes.\\[.2cm]
%\begin{figure}[t]
%\includegraphics[width=8cm]{persp95_30000}\includegraphics[width=8cm]{power95_30000}
%\caption{\label{Fig:2}Empirical percentage type-II-error rate (right) and empirical percentage global testing error rate (left) for the test of size $\alpha=0.05$, depending on tuning parameters $h_n$ and $r_n$, with $50\le h_{30,000}\le 300$ and $2\le r_{30,000}\le 25$. The empirical type-II-error rate measures the empirical amount of realizations smaller or equal the $.95$-quantile of the asymptotic $\chi^2_{N_1}$-distribution under the alternative. The global testing error rate is the sum of the type-II-error rate and the misspecification of the size, that is, the difference between $(1-\alpha)$ and the empirical amount of realizations smaller or equal the $(1-\alpha)$-quantile of the $\chi^2_{N_1}$-distribution, this time on the hypothesis.}
%\end{figure}
%\begin{figure}[t]
%\includegraphics[width=8cm]{persp90_5000}\includegraphics[width=8cm]{persp95_5000}
%\caption{\label{Fig:3}Empirical global testing error rates for the tests of size $\alpha=0.05$ (right) and $\alpha=0.1$ (left) depending on tuning parameters $h_n$ and $r_n$, with $20\le h_{5,000}\le 100$ and $2\le r_{5,000}\le 10$. The global testing error rate is defined as in Figure \ref{Fig:2}.}
%\end{figure}
We investigate the performance of the test for common price and volatility jumps depending on the tuning parameters $h_n$ and $r_n$ in the Monte Carlo simulation. We implement the setup from paragraph \ref{sec:4.1} with $\lambda=2,\eta^{1/2}=0.005$ and $H=0.25$ for both sample sizes $n=30,000$ and $n=5,000$. We set $J_n=30$ in all configurations which is large enough to guarantee high efficiency but smaller than $nh_n$ in any configuration. $J_n^{pi}$ is set equal to $25$.\\
Figure \ref{Fig:2} shows the empirical power and a global testing error including misspecification of the size for a typical testing level $\alpha=0.05$ and for $n=30,000$. The power of all configurations is quite high. Starting with values $r_{30,000}^{-1}=2$, that means the smoothing window is two bins in each direction, the power significantly increases by choosing larger values of $r_{30,000}^{-1}$. However, larger values of $r_{30,000}^{-1}$ lead to a misspecification of the size. The global testing error which adds the misspecification of size with equal weight to the power is minimal for $r_{30,000}^{-1}=4$. On the other hand, the performance is remarkably robust across all considered values of $h_{30,000}$.\\
%Figure \ref{Fig:3} gives similar illustrations of the empirical global testing error for $n=5,000$. We plot the graphs for testing level $\alpha=0.05$ and $\alpha=0.1$, to show that the shape does not change much for different testing levels. 
The precise values of empirical power and size for $n=5,000$, depending on $r_{5,000}$ and $h_{5,000}$ are given in Table \ref{tab2} and Table \ref{tab3}. Again, the global error measure becomes minimal when $r_{5,000}^{-1}=4$, not changing much for $r_{5,000}^{-1}=3$ or $5$, and being very robust with respect to $h_{5,000}$.
\begin{table}[t]
\caption{\label{tab2}Empirical power of the $\alpha=0.05$-test for $n=5,000$ depending on tuning parameters $h_{5,000}$ and $r_{5,000}$.}
\begin{center}
\begin{tabular}{|c|ccccccccc|}
\hline
\begin{tabular}[ht]{lr}
& $r_{5000}^{-1}$:\\
$h_{5000}^{-1}$:& \end{tabular}& 2 & 3&4&5&6&7&8&9&10\\
\hline
20&0.498 &0.737 &0.784 &0.852 &0.890 &0.842 &0.869 & 0.843 & 0.831\\
30&0.557 &0.801 &0.852 &0.896 &0.901 &0.898 &0.927 & 0.925 & 0.937\\
40&0.571 &0.831 &0.879 &0.927 &0.934 &0.944 &0.927 & 0.942 & 0.943\\
50&0.601 &0.834 &0.906 &0.922 &0.954 &0.949 &0.948 & 0.950 & 0.957\\
60&0.603 &0.836 &0.914 &0.933 &0.943 &0.945 &0.968 & 0.968 & 0.972\\
70&0.595 &0.879 &0.921 &0.931 &0.950 &0.965 &0.967 & 0.966 & 0.970\\
80&0.611 &0.848 &0.931 &0.949 &0.965 &0.971 &0.970 & 0.972 & 0.983\\
90&0.629 &0.840 &0.926 &0.957 &0.956 &0.977 &0.977 & 0.982 & 0.984\\
100&0.626&0.842 &0.930 &0.956 &0.978 &0.974 &0.983 & 0.973 & 0.991\\
\hline
\end{tabular}
\end{center}
\end{table}
\begin{table}[t]
\caption{\label{tab3}Empirical size of the $\alpha=0.05$-test, that is, the empirical amount of realizations smaller or equal the $0.95$-quantile of the asymptotic $\chi^2_{N_1}$-distribution, for $n=5,000$ depending on tuning parameters $h_{5,000}$ and $r_{5,000}$.}
\begin{center}
\begin{tabular}{|c|ccccccccc|}
\hline
\begin{tabular}[t]{lr}
& $r_{5000}^{-1}$:\\
$h_{5000}^{-1}$:& \end{tabular}& 2 & 3&4&5&6&7&8&9&10\\
\hline
20&0.953 &0.851 &0.747 &0.732 &0.630 &0.603 &0.541 &0.421 &0.459\\
30&0.975 &0.893 &0.794 &0.753 &0.680 &0.614 &0.592 &0.541 &0.491\\
40&0.975 &0.914 &0.856 &0.781 &0.697 &0.684 &0.608 &0.616 &0.528\\
50&0.973 &0.915 &0.845 &0.804 &0.742 &0.669 &0.675 &0.606 &0.535\\
60&0.977 &0.908 &0.855 &0.795 &0.774 &0.737 &0.662 &0.635 &0.614\\
70&0.976 &0.909 &0.868 &0.792 &0.762 &0.711 &0.673 &0.625 & 0.612\\
80& 0.979 &0.911 &0.868 &0.806 &0.787 &0.734 &0.635 &0.666 &0.612\\
90&0.962 &0.924 &0.872 &0.817 &0.771 &0.713 &0.688 &0.667 &0.603\\
100&0.959 &0.906& 0.879 &0.795& 0.778 &0.728& 0.720 &0.644& 0.660\\
\hline
\end{tabular}
\end{center}
\end{table}
%\clearpage
%\newpage
\subsection{Simulation results for spot volatility estimation with a comparison to a multi-scale approach\label{sec:4.3}}
\begin{figure}[t]
\includegraphics[width=7.45cm]{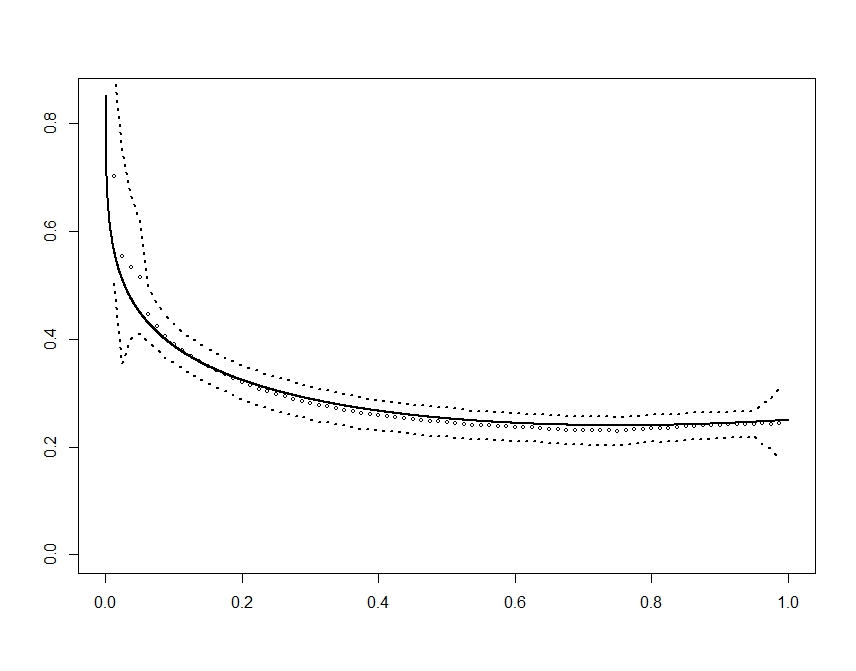}\includegraphics[width=7.45cm]{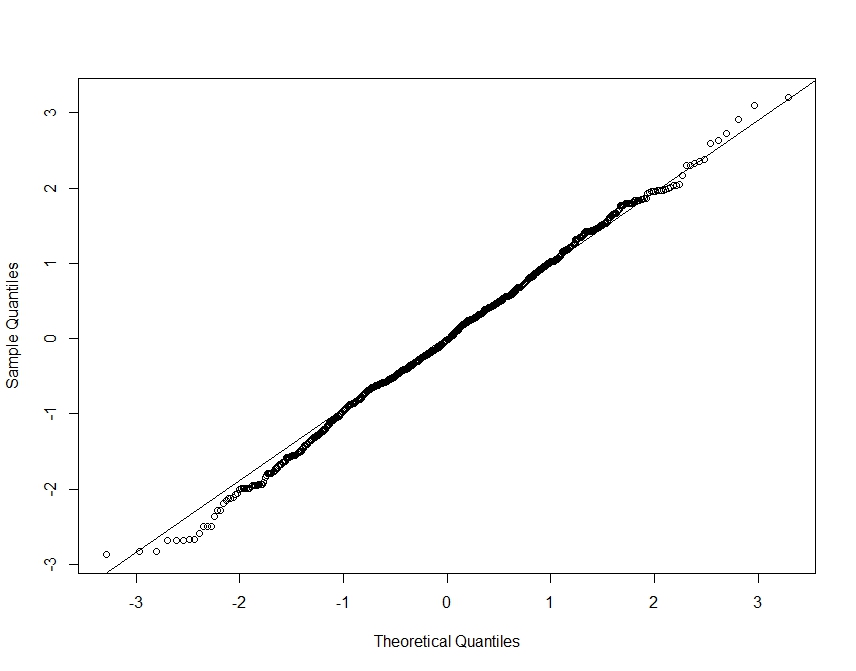}
\caption{\label{Fig:6}Bin-wise averages of spot squared volatility estimates (points) with bin-wise standard deviations (dashed lines) in comparison to the true spot squared volatility (solid line), for $n=5,000$ left. Qq-normal plot for feasible versions of the estimates on bin 40 at $t=1/2$ right.}
\end{figure}

We analyze the accuracy of the spectral spot volatility estimator. First, we illustrate its performance in the model from Section \ref{sec:4.1}, with only a non-random but time-varying volatility component \(\varphi_t=1-\tfrac{3}{5}t^{1/5}+\tfrac{1}{10}t^2\) without volatility jumps. This allows a convenient visualization of the estimation uncertainty. We always use $h_n^{-1}=150$ for $n=30,000$, and $h_n^{-1}=80$ for $n=5,000$, $r_n^{-1}=4$ and $J=30$ and an average of estimators \eqref{rspot} and \eqref{lspot} for the spectral spot volatility estimation of $\sigma_s^2$. For the noise specification \eqref{NS1} with pronounced serial correlations, we use the global version of \eqref{noiseesteq} for the bias-correction terms. Figure \ref{Fig:6} shows the theoretical squared volatility function in comparison to the bin-wise average estimates with standard deviations for $n=5,000$ from 3,000 Monte Carlo runs. The empirical standard deviations on the bins (except the bins close to the boundaries) are quite close to their theoretical values $n^{-1/8}\sqrt{8\sigma_s^3\eta_s^{1/2}}$. For instance, on bin $40$ close to $t=1/2$, we have a ratio of ca.\,1.1 of empirical to theoretical standard deviation. Figure \ref{Fig:6} also depicts the accuracy of the feasible central limit theorem \eqref{feasclt} for bin $40$. 

Next, we compare the performance of our spectral spot volatility estimator to that of a noise-robust multi-scale spot volatility estimator. The multi-scale estimator for integrated volatility is adopted from \cite{zhang}. Applied to all data it estimates $\int_0^1\sigma_t^2\,dt$ and we denote it by $\widehat{\langle X,X\rangle}_1$. In order to obtain an estimator of $\sigma_t^2$ at some $t\in(0,1)$, we use a local difference $\widehat{\langle X,X\rangle}_t-\widehat{\langle X,X\rangle}_{t-\delta}$ with suitable small $\delta$. This extends the methods by \cite{manc15} and \cite{zu14} from two-scale to multi-scale versions. Though no theoretical results are established for this estimator, it is clear that for optimal $\delta$ the approach renders a rate-optimal multi-scale spot volatility estimator. A tuning parameter, the multi-scale frequency, is chosen data-driven in an optimal-way, for which a formula is provided in Section 6 of \cite{bibinger}.\\
The multi-scale estimator gets biased under autocorrelated noise as in \eqref{NS1}. Thus, we focus on noise models without serial correlation to draw a meaningful comparison. First, consider
\begin{align}\label{NS2}\epsilon_i\stackrel{iid}{\sim} N\Big(0,\eta\Big(\int_0^1\varphi_t^4\sigma_t^4\,dt\Big)^{1/2}\Big),~i=0,\ldots,n\,.\end{align}
In this model, the bias-correction of the spectral estimator uses a standard noise variance estimator for i.i.d.\ noise. Further, we examine the estimators in the following noise model with time-varying and endogenous noise:
\begin{align}\label{NS3}\epsilon_i\sim N\Big(0,\eta\Big(\int_0^1\varphi_t^4\sigma_t^4\,dt\Big)^{1/2}+\Big(\sum_{l=1}^5\frac{l}{15}|\Delta_{i-l}^nX|\Big)^2\Big),~i=5,\ldots,n\,,\end{align}
and \eqref{NS2} for $i=0,\ldots,5$. Here, we use locally bin-wise estimated noise levels for the bias-correction terms.\\
Since generated volatility paths in our simulation model are random and thus different in each run, we measure the discrepancy for each path. A suitable global quantity to assess the estimators' qualities from $M$ Monte Carlo iterations is an average normalized mean integrated squared error
\[\operatorname{MISE}=\frac{1}{M}\sum_{m=1}^M\int_0^1\Big(\frac{\hat\sigma_t^2}{\sigma_t^2}-1\Big)^2\,dt\,.\]
\begin{table}[t]
\caption{\label{tabest}Accuracy of spectral and multi-scale spot volatility estimators.}
\begin{center}
\begin{tabular}{ccccc}
\toprule
$n$& noise model & $\eta^{1/2}$&\multicolumn{2}{c}{$\operatorname{MISE}$}\\
  \cline{4-5} \\[-1.05em]
&&&spectral&multi-scale\\ \hline \\[-.6em]
30,000&\eqref{NS2}&0.01&0.0216&0.0713\\
30,000&\eqref{NS2}&0.005&0.0162&0.0421\\
30,000&\eqref{NS2}&0.0025&0.0146&0.0285\\
5,000&\eqref{NS2}&0.01&0.0328&0.0855\\
5,000&\eqref{NS2}&0.005&0.0246&0.0702\\
5,000&\eqref{NS2}&0.0025&0.0227&0.0698\\
\midrule
30,000&\eqref{NS3}&0.01&0.0231&0.0792\\
30,000&\eqref{NS3}&0.005&0.0184&0.0555\\
30,000&\eqref{NS3}&0.0025&0.0170&0.0469\\
5,000&\eqref{NS3}&0.01&0.0597&0.1015\\
5,000&\eqref{NS3}&0.005&0.0540&0.0892\\
5,000&\eqref{NS3}&0.0025&0.0517&0.0875\\
\bottomrule
\end{tabular}
\end{center}
\end{table}
\begin{figure}[ht]
\centering\includegraphics[width=12cm]{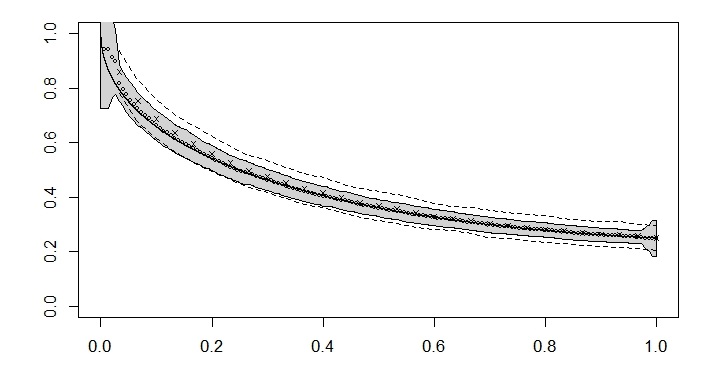}
\caption{\label{Fig:spot}Bin-wise averages of spectral (points) and multi-scale (crosses) spot squared volatility estimates with bin-wise standard deviations (dashed lines) in comparison to the true spot squared volatility (solid line), for $n=30,000$. The area around the spectral estimates determined by their standard deviations is gray colored such that the other dashed lines depict the standard deviations of the multi-scale estimates.}
\end{figure}
The integrals are approximated by sums. For the multi-scale estimator, we set $\delta=K_{MS}^{-1}$ and compute spot volatility estimates on a grid of $K_{MS}$ equidistant time points. An optimization of the $\operatorname{MISE}$ led us to fix $K_{MS}=30$ for $n=30,000$, and $K_{MS}=10$ for $n=5,000$. For the spectral estimator the discretization is given by the $h_n^{-1}$ bins of length $h_n$.\\
An overview of the results for different noise levels, each quantity based on $M=3,000$ Monte Carlo runs, is given in Table \ref{tabest}. The spectral estimation outperforms the ad hoc multi-scale approach in each model specification. The efficiency gains are most relevant for larger noise and more frequent observations. Figure \ref{Fig:spot} visualizes spectral and multi-scale spot volatility estimates with their standard deviations when the true volatility is deterministic and given by \(\varphi_t=1-\tfrac{3}{5}\sqrt{t}+\tfrac{1}{10}t^2\). The confidence regions sketched by the point-wise standard deviations are wider for the multi-scale than for the spectral estimator. We further see a small positive bias of the multi-scale estimates. The discretization, chosen to optimize $\operatorname{MISE}$, is also coarser than the bins of the spectral method which we expect to be the main reason for this bias.\\
Overall, the estimation results for the spectral spot volatility estimator are promising. They confirm that it provides a useful statistical device which is of interest beyond its use as one ingredient for the statistical test for common price and volatility jumps.

\subsection{Simulation results for the test with a comparison to a skip-sampling approach\label{sec:4.4}}
\begin{figure}[t]
\includegraphics[width=7.45cm]{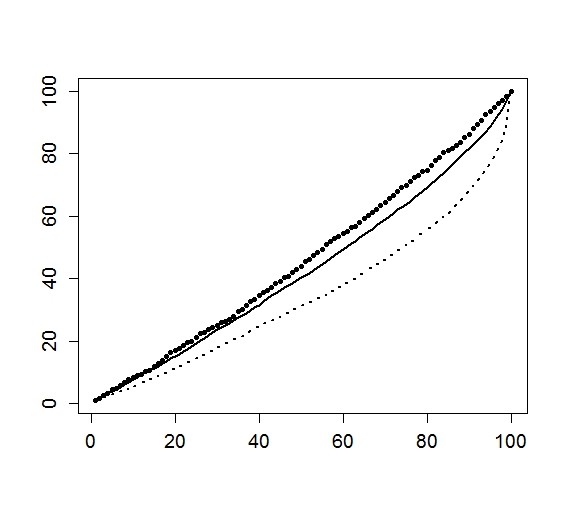}\includegraphics[width=7.45cm]{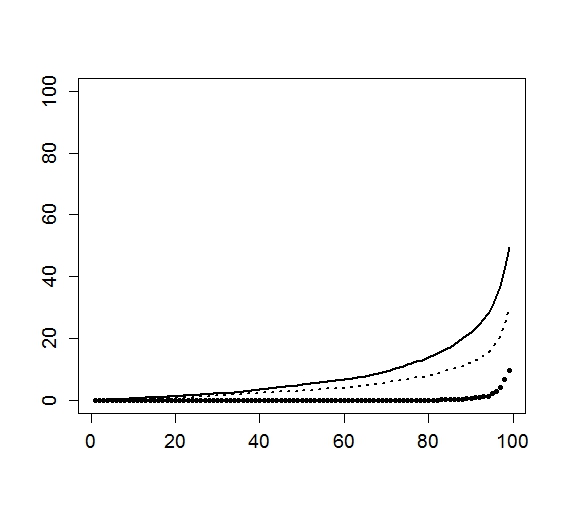}\\[-1cm]
\caption{\label{Fig:4}{\bf{Comparison of the test and the Jacod-Todorov test based on skip-sampled data.}}\newline Empirical size and power of the tests in scenario \rom{1} under the null hypothesis (left) and alternative hypothesis (right). Empirical amount of realizations smaller or equal percentiles of theoretical asymptotic distribution under the null (y-axis) against those percentiles (x-axis). The dotted line shows results for our test and the solid and dashed line two versions of the Jacod-Todorov test using two different tuning parameters. The skip-sampling frequency is optimized to allow for the highest power.} 
\end{figure}

In the sequel, we first study the empirical size and power of our test with respect to different calibrations of volatility jump sizes, noise level and number of observations. To evaluate the improved performance in comparison to the test by \cite{jacodtodorov}, we also implement the latter based on appropriately down sampled discretized simulated paths.\\
The parameter configurations used in the Monte Carlo study for different scenarios are summarized in Table \ref{tab1} together with the chosen tuning parameters according to the values found to be optimal in the sensitivity analysis.
\begin{table}[t]
\caption{\label{tab1}Parameter specification for Monte Carlo.}
\begin{center}
\begin{tabular}{|cccccc:cccc|}
\hline
Scenario& $n$& $\lambda$ & $H$&$\eta^{1/2}$&$\gamma$& $h_n^{-1}$& $J$& $J^{pi}$& $r_n^{-1}$\\
\hline
\rom{1}-Hyp& 30,000& 2 & 0.25 & 0.005 &0& 150 & 30 & 25 & 4\\
\rom{1}-Alt& 30,000& 2&0.25&0.005&1&150&30&25&4\\
\rom{2}-Hyp & 5,000& 2&0.25&0.005&0&80&30&25&4\\
\rom{2}-Alt & 5,000& 2&0.25&0.005&1&80&30&25&4\\
\rom{3}-Alt& 30,000&2 &0.10&0.005&1&200&30&25&5\\
\rom{4}-Hyp& 5,000&2&0.25&0.015&0&80&30&25&4\\
\rom{4}-Alt&5,000&2&0.25&0.015&1&80&30&25&4\\
\hline
\end{tabular}
\end{center}
\end{table}
In scenario \rom{2} (\rom{1}) the average price jump is approx.\,20 (60) times larger than the average absolute return. The identification of price jumps by truncation thus works with only very few errors. Hence, we can use the results from all Monte Carlo iterations to analyze our methods' performance. Examining the ability of thresholding to locate price jumps in different situations has been addressed in \cite{bibwink2015}. Here, the focus is on the test for common price and volatility jumps. The volatility jumps in scenarios \rom{1}, \rom{2} and \rom{4} are a bit smaller than half the size of the average range of the simulated continuous part of the intra-day volatility path. Figure \ref{Fig:7} illustrates that in empirical applications much larger volatility jumps occur. In scenario \rom{3} the jump in the volatility is less than 20\% of the range of the continuous intra-day volatility motion. In scenarios \rom{1}, \rom{2} and \rom{4} we thus have a volatility jump size where the test should attain reasonable power, while scenario \rom{3} investigates the behavior for rather small volatility jumps.\\[.2cm]
We compare the performance of our test based on the statistic \eqref{teststatistic} in scenario \rom{1} for our simulated model with the method by \cite{jacodtodorov}. We cannot apply the latter to the simulated $n=30,000$ high-frequency observations, since the simulated data contains noise. If we apply the test for direct observations to noisy data, the statistics are heavily biased and the performance is very poor. Instead, we skip-sample simulations at a coarser frequency. A heuristic optimization leads us in scenario \rom{1} of our simulation study to an optimal skip-sample frequency resulting in ca.\,500 ``de-noised'' observations on $[0,1]$. For intra-day NASDAQ data this translates in using one observation per 46.8 seconds. \cite{jacodtodorov} employ a one minute frequency for different -- but also very liquid -- data in their application part. Moderate changes of the skip-sampling frequency do not affect the results substantially. Figure \ref{Fig:4} demonstrates a very good performance of our test in scenario \rom{1}. The power is 97.7\% for the $\alpha=0.05$-test and above 90\% even for level $\alpha=0.01$. Similar to our test, the performance of the Jacod-Todorov test applied to the 500 coarse returns is crucially influenced by the length of the smoothing window of local realized volatilities. We visualize two configurations with $k_n=50,100$ in the spot volatility estimators given in \eqref{mle}. The choice $k_n=100$ is in favor of higher power, but the accuracy of the asymptotic quantiles on the null hypothesis is not good. Setting $k_n=50$, we obtain less power but the empirical quantiles on the null hypothesis track the asymptotic ones more closely. In all configurations, the performance of the Jacod-Todorov test applied to skip-sampled data is inferior to the power of our noise-robust approach. This is not surprising, since for our approach we rely on an efficient smoothing technique while skip-sampling can be seen as the simplest method to smooth out noise. The performance of the Jacod-Todorov test is reasonably well also, but in a situation with large available sample sizes and significant noise it is worth to apply the more efficient, noise-robust procedure. If sample sizes are smaller (and the noise not larger), the difference between the two methods becomes smaller.
\begin{figure}[t!]
\hspace*{.715cm}Scenario \rom{2} - Hyp ~~~~~~~~~~Scenario \rom{2} - Alt~~~~~~~~~~~~Scenario \rom{3} - Alt ~~~~~~~~~~Scenario \rom{4} - Alt\\
\includegraphics[width=3.7cm]{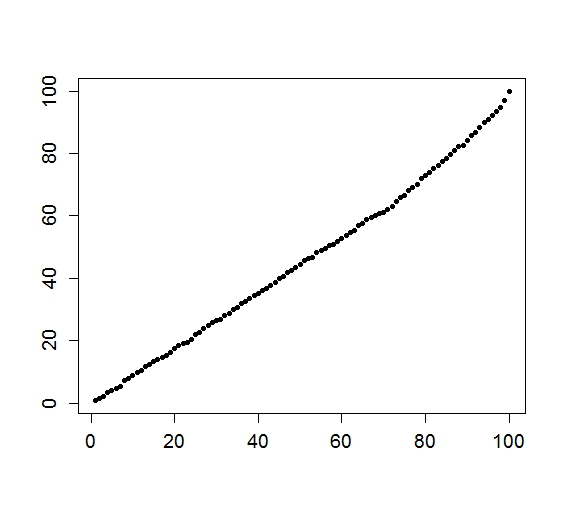}\includegraphics[width=3.7cm]{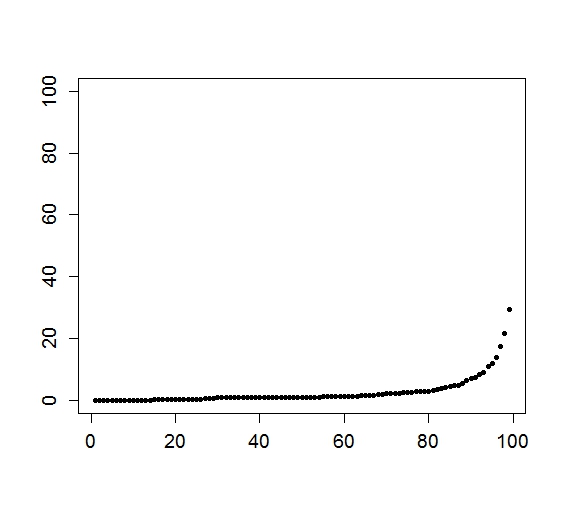}\includegraphics[width=3.7cm]{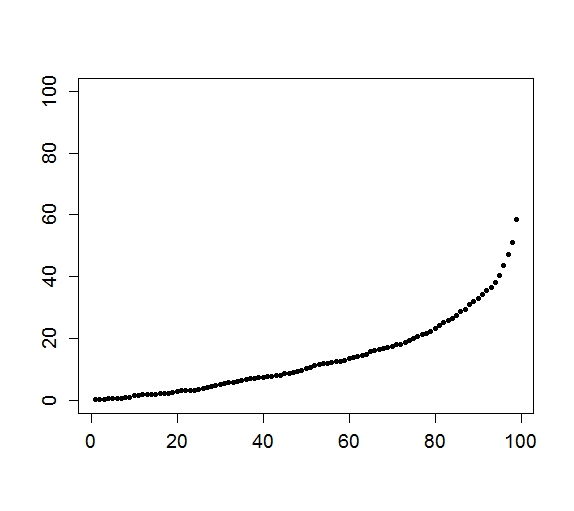}\includegraphics[width=3.7cm]{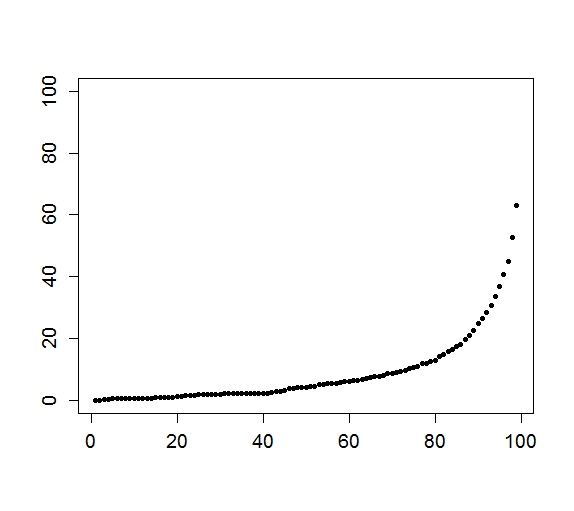}
\caption{\label{Fig:5}Empirical size and power of the tests in scenarios \rom{2}, \rom{3} and \rom{4}. Empirical amount of realizations smaller or equal percentiles of theoretical asymptotic distribution under the null (y-axis) against those percentiles (x-axis).}
\end{figure}
Figure \ref{Fig:5} shows the performance in other scenarios \rom{2}, \rom{3} and \rom{4}. Decreasing the sample size to $n=5,000$ observations in scenario \rom{2}, while all parameters are equal as in scenario \rom{1}, leads to a slightly smaller power and larger misspecification of the size. The power is still higher than for the skip-sample approach, but the difference is less relevant. With the tuning parameters which minimize the global empirical testing error, the misspecification of the size is still acceptable. Larger noise levels result in smaller power as shown for scenario \rom{4} in Figure \ref{Fig:5}, while the fit of the size remains good. In this situation, the Jacod-Todorov method would only work for less frequent skip-sampling resulting as well in smaller power. For the alternative hypothesis with a small volatility jump in scenario \rom{3}, a sensitivity analysis as in Section \ref{sec:4.2} led us to slightly different tuning parameters, $h_{30,000}^{-1}=200$ and $r_{30,000}^{-1}=5$. Since smaller bins give a higher time resolution, it is not surprising that detecting small volatility jumps in a rapidly time-varying spot volatility works better for a finer time resolution. On the other side, choosing $r_n^{-1}$ slightly larger leads to almost the same window length $r_n^{-1}h_n$ for spot volatility estimation as before. The power for such small volatility jumps is less, but still ca.\,60\% for $\alpha=0.05$.
%\newpage\noindent

\section{Data study\label{sec:5}}

To provide evidence about the practical relevance of price-volatility co-jumps and to study the usefulness of our estimators and test in a real-world data environment, we apply our methodology to stocks traded at the exchange platform NASDAQ. The data study is based on limit order book data taken from the online data tool LOBSTER\footnote{LOBSTER academic data- \url{lobsterdata.com}, powered by NASDAQ OMX}. The example refers to stocks of the online and technology companies Amazon.com Inc.\ (AMZN), Apple Inc.\ (AAPL), Facebook Inc.\ (FB), Intel Corp.\ (INTC) and Microsoft Corp.\ (MSFT). We focus on transaction prices of 252 trading days in the year 2013. A trading day spans from 9:30 to 16:00 EDT and includes for a single stock a minimum of 4,267 (AMZN 2013-07-03) up to a maximum of 210,812 (FB 2013-10-31) transactions. One benefit of our estimator and test is that we can directly plug-in traded log-prices, reconstructed from the order book, without considering any skip-sampling or synchronization procedures. Since the method is robust against market microstructure noise, we efficiently take into account all information stored in the data.\\
\begin{table}[t]
\begin{center}
\caption{\label{tab4}Testing for price and volatility jumps in NASDAQ order book data.}
\label{emp}
\vspace{.25cm}
\begin{tabular}{lccccccc}
 \hline\hline \vspace{-.8em} \\
\multirow{2}{*}{Stocks}&\multirow{2}{*}{\centering \parbox{2.3cm}{\centering\# of days with price jumps}} &
\multicolumn{2}{c}{\rotatebox{0}{\parbox{3.58cm}{\centering Rejection rate  (price-volatility jumps)}}} &&
\multicolumn{3}{c}{\rotatebox{0}{\parbox{3.5cm}{\centering Sample Averages (whole year 2013)}}} \vspace{.2em} \\ \cline{3-4} \cline{6-8}
\vspace*{-.85em}\\

&&$\alpha=5\%$&$\alpha=10\%$&&$n$&$\Delta\hat\sigma^2_s$&$\hat{IV} $ \\ \hline \vspace{-.8em} \\
Amazon.com Inc.\         &21 & 52.4\% & 61.9\%&&10,924&$31.2\%$&1.47 \vspace{.2em} \\
Apple Inc.\                   & 22 & 63.6\% & 72.7\%&&36,947&$36.5\%$&1.52 \vspace{.2em} \\
Facebook Inc.\     		   & 37 & 46.0\% & 51.4\%&&41,354&27.8\% &$3.12$  \vspace{.2em} \\
Intel Corp.\     	        	   & 47 & 27.7\% & 36.2\%&&18,535&$23.0\%$&0.93 \vspace{.2em} \\
Microsoft Corp.\     	   & 22 & 50.0\% & 50.0\%&&28,052&$31.2\%$&0.97 \vspace{.2em}
\vspace{.2em} \\ \hline \hline
\multicolumn{8}{p{14cm}}{\footnotesize Notes: Estimation and test executed for each day in the year 2013 separately.  $n$ indicates the average number of observed trades per trading day, $\Delta\hat\sigma^2_s=|\log(\hat\sigma_{s}^2)-\log(\hat\sigma_{s-}^2)|$ is the average estimated relative size of volatility jumps, $\hat{IV}$ the average spectral estimate of the integrated squared volatility times $10^{-4}$. Sample period:  01-02-2013 to 12-31-2013 (252 days).}
\end{tabular}
\end{center}
\end{table}
Guided by our theoretical results and the simulations, estimates and tests are based on spectral statistics calculated for $k=0,1,...,h_n^{-1}-1$ bins over a trading day, with $h_n^{-1}=\lfloor 3\sqrt{n}/\log(n)\rfloor $. 
We set $J=30$ and $J^{pi}=15$. Jumps in prices are detected with the locally adaptive threshold $\hat u_k=2\log(h_n^{-1})h_n\hat\sigma_{kh_n,pil}^2$, with $\hat\sigma_{kh_n,pil}^2$ the pilot estimator \eqref{pilot} of the spot squared volatility. We fix constant window lengths $r_n^{-1}=4$. Surely, $r_n^{-1}$ determines a crucial parameter which can be studied to learn about the persistence or live-time of a break in spot volatility. 
\begin{figure}[t]
\begin{center}
\includegraphics[width=14.8cm]{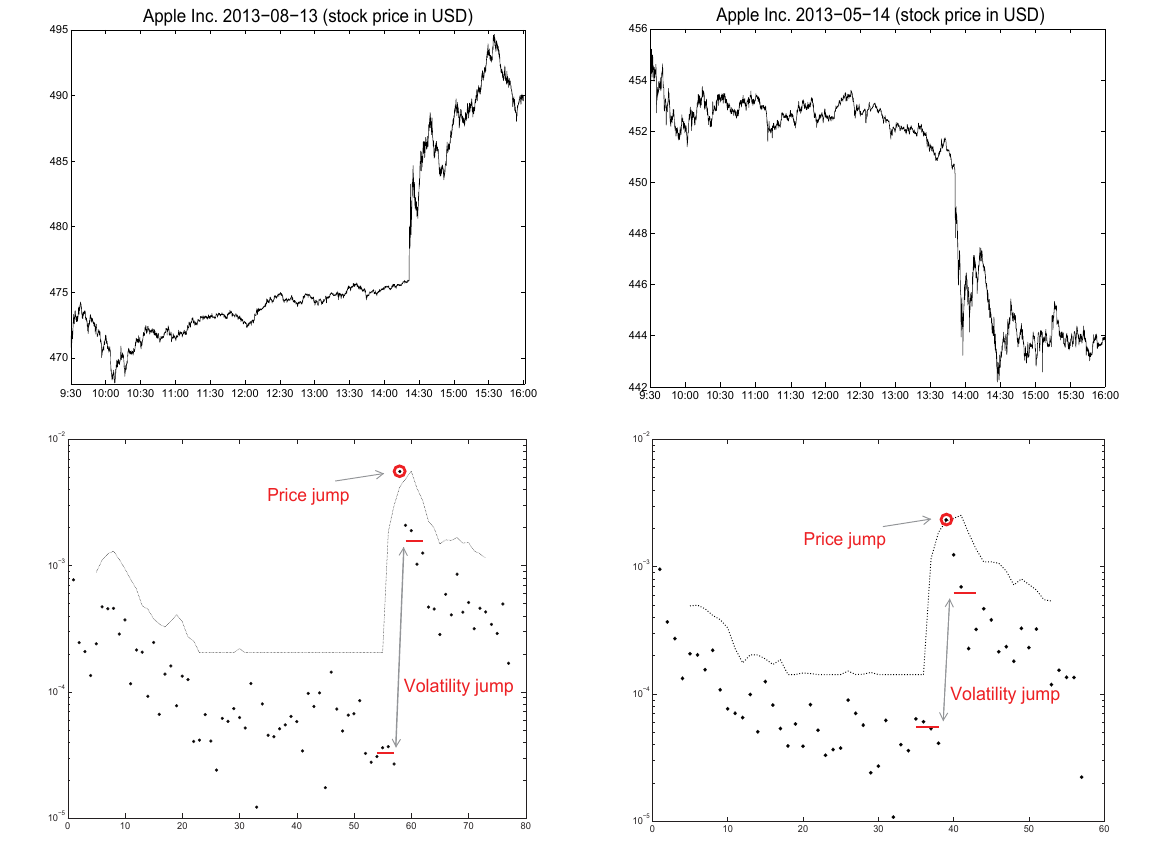}
\caption{\label{Fig:7} {\bf Examples of common price and volatility jumps.} Upper figures indicate price processes of the Apple Inc.\ stock. Lower figures display the related spectral quadratic variation estimates for the bins $k=0,\ldots,h_n^{-1}$. Estimates exceeding the threshold (dotted line) are marked as price jumps. The difference between spot volatility estimates (bars) estimate the volatility jump. 2013-08-13: $n=87,445$, 2013-05-14: $n=40,707$.}
\end{center}
\end{figure}
We apply the test to each day separately.

Table \ref{emp} reports the rejection rates for the 5\% and 10\% significance levels. Results indicate that on a 10\% significance level 36\% (INTC) up to 73\% (AAPL) of jumps in prices are accompanied by jumps in volatility.  It appears that the rejection rate decreases in the number of detected price jumps. This leads to relatively stable frequencies of price-volatility co-jumps over time across the considered stocks. Referring to the 5\% significance level, the Amazon.com stock displays with around 4.4\% of the trading days the lowest frequency of common price and volatility jumps. With around 6.7\% of trading days, Facebook Inc.\ has the largest number of common jumps. Absolute jump sizes of the log squared volatility processes reported in Table \ref{emp} are considerably large.    

Figure \ref{Fig:7} illustrates the mechanisms behind the test for common price and volatility jumps. Left hand plots show an upward jump in prices on bin $k=58$, whereas right hand plots show a downward jump in prices on bin $k=39$. Both price jumps are associated with a significant contemporaneous upward jump in spot volatility. The $p$-value in both examples is 0.00. On the first example date, August, 13th 2013, the investor Carl Icahn has taken a large stake of AAPL stocks. On May 14th, the downward jump example date, figures of mobile phone sales have been reported.\\
We find evidence for frequent occurrences of simultaneous jumps in price and volatility and quite large volatility jump sizes. Yet, by far not all detected price jumps are accompanied by volatility jumps. Understanding the economic sources of different jump events and their consequences for price-volatility co-jumps is of interest for future research.

\section{Conclusion\label{sec:6}}
We present a new test for the presence of contemporaneous jumps of price and volatility based on high-frequency data. The test transfers the methodology of \cite{jacodtodorov} to a setup accounting for microstructure noise by employing a spectral estimation of the spot volatility and an accurate test function. The nonparametric spot volatility estimator shows appealing asymptotic and finite-sample qualities and is of interest beyond the scope of this article. It opens up several new ways for inference in models for high-frequency financial data with noise. Simulations demonstrate that the proposed noise-robust test increases the finite-sample performance considerably compared to an application of the test by \cite{jacodtodorov} to skip-sampled data. Our data study reveals cogent significance of price and volatility co-jumps in NASDAQ high-frequency data. %As a consequence, they should be considered 
%This has consequences for the future modeling of price and volatility. %Investigating why many price jumps are accompanied by volatility adjustments, but others not, appears important through the lenses of economic theory of information processing and surprise elements. 
%Since many price adjustments are due to news releases the topic of simultaneous volatility adjustments is as well important to analyze for communication management.\\
The presented methods can be generalized in various directions. For instance, our methods guide the way how a test for correlation of price and volatility jumps, as presented by \cite{gernot} for a non-noisy observation design, can be constructed. A general global test for volatility jumps under noise generalizing the methods from \cite{BJV} could be addressed with a related high-dimensional testing procedure.
%\section*{Acknowledgement}
%Data has been provided by LOBSTER academic data, powered by NASDAQ OMX.
%-------------------------------------------------------------------------------------------------------------
\section{Proofs\label{sec:7}}
\subsection{Preliminaries}
 On the finite time horizon $[0,1]$, we may augment local boundedness to uniform boundedness in Assumption \ref{H} and Assumption \ref{eff}, such that we can assume that there exists a constant $\Lambda$ with
 \begin{align*}\max{\{|b_s(\omega)|,\sigma^2_s(\omega),|X_s(\omega)|,|\delta_{\omega}(s,x)|/\gamma(x),\eta_s(\omega),\rho_s(\omega)\}}\le \Lambda\,,\end{align*}
for all $(\omega,s,x)\in(\Omega,\mathds{R}_+,\mathds{R})$.
This standard procedure can be found in Section 4.4.1 of \cite{JP}.
Throughout the proofs $K$ is a generic constant and $K_p$ a constant emphasizing dependence on $p$.
We decompose the semimartingale $X$ in its continuous part
\[C_t =X_0+\int_0^tb_s\,ds+\int_0^t\sigma_s\,dW_s\,,\]
and the jump component
\[J_t =\int_0^t\int_{\mathds{R}}\delta(s,x)\1_{\{|\delta(s,x)|\le 1\}}(\mu-\nu)(ds,dx) +\int_0^t\int_{\mathds{R}}\delta(s,x)\1_{\{|\delta(s,x)|> 1\}}\mu(ds,dx)\,.\]
The processes
\begin{align}\tilde C^n_t=\int_0^t \sigma_{\lfloor sh_n^{-1}\rfloor h_n}\,dW_s\end{align}
serve as an approximation of $C_t$ by simplified processes without drift and with locally constant volatility.
We separate jumps with absolute value bounded from above by some $\varepsilon<1$ and larger jumps:
 \begin{align*}J_t=J(\varepsilon)_t+\int_0^t\int_{\mathds{R}\backslash A_{\varepsilon}}\delta(s,x)\1_{\{|\delta(s,x)|\le 1\}}(\mu-\nu)(ds,dx)+\int_0^t\int_{\mathds{R}}\delta(s,x)\1_{\{|\delta(s,x)|> 1\}}\mu(ds,dx)\,,\end{align*}
 with $A_{\varepsilon}=\{z\in\mathds{R}|\gamma(z)\le \varepsilon\}$ and later let $\varepsilon\rightarrow 0$.
 Let us recall some usual estimates on Assumptions \ref{eff}, \ref{H} and \ref{sigma} which are crucial for the following proofs.
 \begin{subequations}
For the continuous semimartingale part, we have
 \begin{align}\label{e1}\forall p\in \mathds{N},s,t\ge 0:~\E\left[|C_{s+t}-C_s|^p\big|\mathcal{F}_s\right]\le K_p t^{\frac{p}{2}}\,.\end{align}
For given $0<\varepsilon<1$, for $J(\varepsilon)$ the estimate
 \begin{align}\notag \forall p\in \mathds{N},\forall s,t\ge 0:~\E\left[|J(\varepsilon)_{s+t}-J(\varepsilon)_s|^p\big|\mathcal{F}_s\right]&\le K_p\,\E\Big[\Big(\int_s^{(s+t)}\int_{A_{\varepsilon}}(\gamma^2(x)\wedge 1)\mu(d\tau,dx)\Big)^{\frac{p}{2}}\Big]\\ \label{e2}&  \le K_p t^{(\frac{p}{2}\wedge 1)}\gamma_{\varepsilon}^{(\frac{p}{2}\wedge 1)}\,,\end{align}
 holds with \(\gamma_{\varepsilon}=\int_{A_{\varepsilon}}\left(\gamma^2(x)\wedge 1\right)\lambda(dx) \le K \varepsilon^{(2-r)}\).\\
The continuous semimartingale increments satisfy local Gaussianity in the sense that
 \begin{align}\notag \forall p\in \mathds{N},s,t\ge 0:~\E\left[|C_{s+t}\hspace*{-.05cm}-\hspace*{-.05cm}C_s\hspace*{-.05cm}-\hspace*{-.05cm}(\sigma_s(W_{s+t}\hspace*{-.05cm}-\hspace*{-.05cm}W_s))|^p\big|\mathcal{F}_s\right]& \le K_p \,\E\Big[\Big(\int_s^{s+t}|\sigma_{\tau}-\sigma_s|^2\,d\tau\Big)^{\frac{p}{2}}\big|\mathcal{F}_s\Big]\\
 &\notag \le K_p \, t^{\frac{p}{2}}\,\E\Big[\sup_{\tau\in[s,s+t]}{\left(|\sigma_{\tau}-\sigma_s|^p\right)}\big|\mathcal{F}_s\Big]\\
&\label{e4}\le K_p t^{\frac{p}{2}(1+2\alpha)}\,,
  \end{align}
 on Assumption \ref{sigma}. The probability of a frequent occurrence of large jumps is small. Precisely, the expectation of jumps with absolute value larger than $\varepsilon$ is bounded:
  \begin{align}\label{e5}\forall s,t\ge 0:~\E\left[|J_{s+t}-J_s-\left(J(\varepsilon)_{s+t}-J(\varepsilon)_s\right)|\big|\mathcal{F}_s\right]\le K t\varepsilon^{-r}\,.\end{align}
Under Assumption \ref{H} with $r\ge 1$, the jumps moreover satisfy
 \begin{align}\notag\forall s,t\ge 0:~\E\left[|J_{t}-J_s|^p\big|\mathcal{F}_s\right]&\le K_p\,\E\Big[\Big(\int_s^{t}\int_{\R}(\gamma^r(x)\wedge 1)\lambda(dx)ds\Big)^{p/r}\Big]\\ \label{boundjumps}&  \le K_p |t-s|^{((p/r)\wedge 1)}\,.\end{align}
Under Assumption \ref{sigma} for $0\le s<t\le 1$, the squared volatility satisfies:
\begin{align}\E[|\sigma^2_t-\sigma^2_s|\,|\mathcal{F}_s]\le|t-s|^{\alpha}\,.\end{align}
\end{subequations}
Proofs of these bounds can be found, for instance, in Chapter 13 of \cite{JP}. \eqref{e2} follows from Equation (54) in \cite{aitjac10}.\\
In the sequel, we gather more properties of the basis functions \eqref{Phi}. We define $(\Phi_{jk})$ in \eqref{Phi} in the same way as \cite{bibingerreiss} in their Equation (4b) to exploit discrete-time Fourier identities under equidistant sampling. The asymptotic properties of the estimator remain the same when we use
\begin{align}\label{Phi2}
\tilde\Phi_{jk}(t)&=\frac{\sqrt{2h_n}}{j\pi} \sin{\left(j\pi h_n^{-1}\left(t-kh_n\right)\right)}\1_{[kh_n,(k+1)h_n)}(t),j\ge 1\,,
\end{align}
instead which equals the definition from Equation (2.2) in \cite{BHMR}. We heavily exploit the summation by parts identity for spectral statistics
\begin{align}\label{sbp}S_{jk}=\|\Phi_{jk}\|_n^{-1} \Bigg(\sum_{i=1}^n\Delta_i^n X\Phi_{jk}\big(\tfrac{i}{n}\big)-\sum_{i=1}^{n-1}\epsilon_i\, \varphi_{jk}\Big(\frac{i+1/2}{n}\Big)\frac{1}{n}\Bigg)\,,\end{align}
with $\varphi_{jk}(t)=\sqrt{2}h_n^{-1/2}\cos{\big(j\pi h_n^{-1}(t-kh_n)\big)}\1_{[kh_n,(k+1)h_n]}(t)$, see Lemma 6.1 of \cite{stable}. For all $(\varphi_{jk})$, it holds that
\begin{align}n^{-1}\sum_{i=1}^{n-1}\varphi_{jk}\Big(\frac{i+1/2}{n}\Big)\varphi_{j'k}\Big(\frac{i+1/2}{n}\Big)=\delta_{jj'}\,.\end{align}
For the asymptotic theory, we shall further use the following identities
\begin{subequations}
\begin{align}\label{squarephi} \int_0^1\tilde\Phi_{jk}^2(t)\,dt=\frac{h_n^2}{\pi^2j^2}~,~\int_0^1\varphi_{jk}^2(t)\,dt=1\,,\\
\label{crossphi}\int_0^1 \tilde\Phi_{jk}(t)\,\varphi_{uk}(t)\,dt=\frac{(1-\cos(\pi j)\cos(\pi u))2 h_n}{\pi^2(j^2-u^2)}\,.\end{align}
\end{subequations}
The latter gives $4h_n/(\pi^2(j^2-u^2))$ whenever $j$ is odd and $u$ even, or the other way round, and vanishes in all other cases. 
Recall the definition of the weights \eqref{orweights}. The magnitude of these weights is
\begin{align}\notag w_{jk}\le I_{jk}=\tfrac{1}{2}\Big(\sigma^2_{kh_n}+\tfrac{\eta_{kh_n}}{n}\|\Phi_{jk}\|^{-2}_{n}\Big)^{-2}&=\mathcal{O}\Big( \Big(1+\frac{j^2}{nh_n^2}\Big)^{-2}\Big)\\
&\label{orderweights}= \begin{cases}\mathcal{O}(1)~~~~~~~~~~~~~~\,\mbox{for}\,j\le \sqrt{n}h_n\\ \mathcal{O}(j^{-4}n^2h_n^{4})~~\mbox{for}\,j>\sqrt{n}h_n\end{cases}\,,\end{align}
with $\|\Phi_{jk}\|_n^{-2}\approx \pi^2j^2h_n^{-2}=\big( \int_0^1\tilde\Phi_{jk}^2(t)\,dt\big)^{-1}$.

In the proofs, we use the notation $\zeta_k^{ad}(Z)$ and $\zeta_k(Z)$ from \eqref{zeta} analogously also for different processes $Z$. This means that we insert in \eqref{zeta} spectral statistics $S_{jk}(Z)$, analogous to \eqref{spectralstatistics}, computed from the sequence $Z_{i/n},i=0,\ldots,n$, especially $\zeta_k(X)$ for the statistics based on the unobserved efficient price.
\subsection{Estimation of the noise long-run variance}
First, consider the standard case where $\alpha\le 1/2$ in Assumption \ref{sigma}, such that $\beta< 1/4$. To estimate $(\eta_{kh_n})$ under \eqref{etasmooth}, we use $nh_n$ observations on the bin $[kh_n,(k+1)h_n]$. For $k=0,\ldots,h_n^{-1}-1$, and $u=1,\ldots,Q$, define the cumulative empirical autocorrelation statistics
\begin{align*}
Z_{kh_n}^{(u)}&=\frac{1}{2nh_n}\sum_{i=knh_n+1}^{(k+1)nh_n}(\Delta_i^n Y)^2+\frac{1}{nh_n}\sum_{l=1}^{u}\sum_{i=knh_n+1}^{(k+1)nh_n-u}\Delta_i^n Y\Delta_{i+l}^n Y\,,\\
\tilde Z_{kh_n}^{(u)}&=\frac{1}{2nh_n}\sum_{i=knh_n+1}^{(k+1)nh_n}(\Delta_i^n Y)^2+\frac{1}{nh_n}\sum_{l=1}^{u}\sum_{i=knh_n+u+1}^{(k+1)nh_n}\Delta_i^n Y\Delta_{i-l}^n Y\,.
\end{align*}
For $u=0$, the rescaled local realized volatilities in the first addend define $Z_{kh_n}^{(0)}$. We estimate $\eta_{kh_n}$ by
\begin{align}\label{noiseesteq}\hat\eta_{kh_n}=\sum_{u=0}^Q(u+1)Z_{kh_n}^{(u)}+\sum_{u=1}^Qu \tilde Z_{kh_n}^{(u)}\,.\end{align}
We assume that $Q$ is known. However, the same result applies if the process is $\tilde Q$-dependent with $\tilde Q<Q$. It thus suffices to take $Q$ sufficiently large. A statistical method to infer $Q$ is provided by \cite{BHMR2}.\\ 
We consider separately the case $\alpha>1/2$ with possible values $1/4\le \beta<1/3$. Then, we exploit the increased smoothness of the noise by \eqref{etasmooth} to estimate $(\eta_{kh_n})$ with an improved convergence rate. We partition $[0,1]$ in $n/M_n$ windows of lengths $M_n/n$, each with $M_n$ observations, where $M_n=c_M\,n^{1-(2\alpha+1)^{-1}}$. For a simple exposition we may suppose $M_n,n/M_n\in\N$ again. Completely analogously as before, we compute the cumulative empirical autocorrelation statistics $Z_{k\frac{M_n}{n}}^{(u)},\tilde Z_{k\frac{M_n}{n}}^{(u)}$ for $k=0,\ldots,n/M_n-1$. The estimates $(\hat\eta_{kh_n})$ are now obtained by 
\begin{align}\label{noiseesteq2}\hat\eta_{kh_n}=\hat\eta_{\tilde k\frac{M_n}{n}}\1_{\{kh_n\in[\tilde k\frac{M_n}{n},(\tilde k+1)\frac{M_n}{n})\}}\,,\end{align}
with $\hat\eta_{\tilde k\frac{M_n}{n}}$ analogous to \eqref{noiseesteq} over the coarser time windows.\\[.1cm]
\newpage\noindent
\textbf{Proof of Proposition \ref{cornoiseest}}\\[.1cm]
We begin with the case $\alpha\le 1/2$ in Assumption \ref{sigma}, such that $\beta< 1/4$. We prove that 
\begin{align}\label{noiseestsp}\hat\eta_{kh_n}=\eta_{kh_n}+\mathcal{O}_{\P}\big(n^{-1/4}\big)\,.\end{align}
Considering the expectation of the cumulative empirical autocorrelation statistics, all terms involving increments $\Delta_i^n X$ are of order $\mathcal{O}_{\P}(n^{-1/2})$, and even smaller under exogenous noise. Thus, we have that
\begin{align*}\E\big[Z_{kh_n}^{(u)}|{\cal{F}}_{kh_n}\big]&=\frac{1}{2nh_n}\sum_{i=knh_n+1}^{(k+1)nh_n}\E\big[\epsilon_i^2+\epsilon_{i-1}^2-2\epsilon_i\epsilon_{i-1}|{\cal{F}}_{kh_n}\big]+\mathcal{O}_{\P}(n^{-1/2})\\
& \quad +\frac{1}{nh_n}\sum_{l=1}^{u}\sum_{i=knh_n+1}^{(k+1)nh_n-u}\E\big[\epsilon_i\epsilon_{i+l}+\epsilon_{i-1}\epsilon_{i+l-1}-\epsilon_{i-1}\epsilon_{i+l}-\epsilon_{i}\epsilon_{i+l-1}|{\cal{F}}_{kh_n}\big]\\
&=\frac{1}{nh_n}\sum_{i=knh_n+1}^{(k+1)nh_n}\E\big[\epsilon_i^2-\epsilon_i\epsilon_{i-1}|{\cal{F}}_{kh_n}\big]+\mathcal{O}_{\P}(n^{-1/2})\\
& \quad +\frac{1}{nh_n}\sum_{i=knh_n+1}^{(k+1)nh_n-u}\E\big[\epsilon_i\epsilon_{i+u}+\epsilon_{i-1}\epsilon_{i}-\epsilon_{i-1}\epsilon_{i+u}-\epsilon_{i}^2|{\cal{F}}_{kh_n}\big]\,,
\end{align*}
where we use $(nh_n)^{-1}=\mathcal{O}(n^{-1/2})$ for the first and the telescoping sum for the second addend. We obtain that
\begin{align*}\E\big[Z_{kh_n}^{(u)}|{\cal{F}}_{kh_n}\big]= \frac{1}{nh_n}\sum_{i=knh_n+1}^{(k+1)nh_n-u}\E\big[\epsilon_i\epsilon_{i+u}-\epsilon_{i-1}\epsilon_{i+u}|{\cal{F}}_{kh_n}\big]+\mathcal{O}_{\P}(n^{-1/2})\,,
\end{align*}
for all $0\le u\le Q$. Summing over $u\in\{0,\ldots,Q\}$, we exploit another telescoping sum:
\begin{align*}\sum_{u=0}^Q(u+1)\E\big[ Z_{kh_n}^{(u)}|{\cal{F}}_{kh_n}\big]&=\frac{1}{nh_n}\sum_{i=knh_n+1}^{(k+1)nh_n-Q}\sum_{u=0}^Q(u+1)\big(\cov\big(\epsilon_i,\epsilon_{i+u}|{\cal{F}}_{kh_n}\big)-\cov\big(\epsilon_{i-1},\epsilon_{i+u}|{\cal{F}}_{kh_n}\big)\big)\\ &\quad+\mathcal{O}_{\P}(n^{-1/2})\\
&=\frac{1}{nh_n}\sum_{i=knh_n+1}^{(k+1)nh_n-Q}\sum_{u=0}^Q\cov\big(\epsilon_i,\epsilon_{i+u}|{\cal{F}}_{kh_n}\big)+\mathcal{O}_{\P}(n^{-1/2})\,,\end{align*}
since $\cov\big(\epsilon_{i-1},\epsilon_{i+Q}|{\cal{F}}_{kh_n}\big)=0$. There are at most $\tilde Q<\infty$ addends $i=knh_n+1,\ldots,knh_n+\tilde Q$, for that $\E[\epsilon_i|{\cal{F}}_{kh_n}]\ne 0$ is possible by endogeneity, which are asymptotically negligible in the above sum. A similar computation for $\tilde Z_{kh_n}^{(u)}$ gives:
\begin{align*}\E\big[\tilde Z_{kh_n}^{(u)}|{\cal{F}}_{kh_n}\big]&=\frac{1}{2nh_n}\sum_{i=knh_n+1}^{(k+1)nh_n}\E\big[\epsilon_i^2+\epsilon_{i-1}^2-2\epsilon_i\epsilon_{i-1}|{\cal{F}}_{kh_n}\big]+\mathcal{O}_{\P}(n^{-1/2})\\
& \quad +\frac{1}{nh_n}\sum_{l=1}^{u}\sum_{i=knh_n+u+1}^{(k+1)nh_n}\E\big[\epsilon_i\epsilon_{i-l}+\epsilon_{i-1}\epsilon_{i-l-1}-\epsilon_{i-1}\epsilon_{i-l}-\epsilon_{i}\epsilon_{i-l-1}|{\cal{F}}_{kh_n}\big]\\
&=\frac{1}{nh_n}\sum_{i=knh_n+1}^{(k+1)nh_n}\E\big[\epsilon_{i-1}^2-\epsilon_i\epsilon_{i-1}|{\cal{F}}_{kh_n}\big]+\mathcal{O}_{\P}(n^{-1/2})\\
& \quad +\frac{1}{nh_n}\sum_{i=knh_n+u+1}^{(k+1)nh_n}\E\big[\epsilon_i\epsilon_{i-1}+\epsilon_{i-1}\epsilon_{i-u-1}-\epsilon_{i-1}^2-\epsilon_{i}\epsilon_{i-u-1}|{\cal{F}}_{kh_n}\big]\,,
\end{align*}
and thus that
\begin{align*}\sum_{u=1}^Qu \E\big[\tilde Z_{kh_n}^{(u)}|{\cal{F}}_{kh_n}\big]&=\frac{1}{nh_n}\sum_{i=knh_n+Q+1}^{(k+1)nh_n}\sum_{u=1}^Qu\big(\cov\big(\epsilon_i,\epsilon_{i-u}|{\cal{F}}_{kh_n}\big)-\cov\big(\epsilon_{i},\epsilon_{i-u-1}|{\cal{F}}_{kh_n}\big)\big)\\ &\quad+\mathcal{O}_{\P}(n^{-1/2})\\
&=\frac{1}{nh_n}\sum_{i=knh_n+Q+1}^{(k+1)nh_n}\sum_{u=1}^Q\cov\big(\epsilon_i,\epsilon_{i-u}|{\cal{F}}_{kh_n}\big)+\mathcal{O}_{\P}(n^{-1/2})\,,\end{align*}
since $\cov\big(\epsilon_{i},\epsilon_{i-Q-1}|{\cal{F}}_{kh_n}\big)=0$ for all, except finitely many, $i$. This yields that for the estimator \eqref{noiseesteq}
\begin{align}\label{expnoisep}\E\big[\hat\eta_{kh_n}|{\cal{F}}_{kh_n}\big]=\frac{1}{nh_n}\sum_{i=knh_n+Q+1}^{(k+1)nh_n-Q}\sum_{u=-Q}^Q\cov\big(\epsilon_i,\epsilon_{i+u}|{\cal{F}}_{kh_n}\big)+\mathcal{O}_{\P}(n^{-1/2})\,,\end{align}
such that $\sup_{t\in[kh_n,(k+1)h_n]}|\eta_{kh_n}-\eta_t|=\mathcal{O}\big(h_n^{(1/2+\delta)\vee \alpha}\big)=\mathcal{O}\big(n^{-1/4}\big)$ and \eqref{lrvar} give that 
\begin{align}\label{expnoise}\E\big[\hat\eta_{kh_n}|{\cal{F}}_{kh_n}\big]=\eta_{kh_n}+\mathcal{O}_{\P}\big(n^{-1/4}\big)\,.\end{align}
The following bound for the conditional variance of the estimator \eqref{noiseesteq} completes the proof of \eqref{noiseestsp}. It holds uniformly in $k$ that
\begin{align*}
\var\big(Z_{kh_n}^{(u)}|{\cal{F}}_{kh_n}\big)&\le \frac{1}{2n^2h_n^2}\sum_{i=knh_n+1}^{(k+1)nh_n}\sum_{u=-Q}^Q\cov\big((\Delta_i^n \epsilon)^2,(\Delta_{i+u}^n \epsilon)^2|{\cal{F}}_{kh_n}\big)\\
&\quad +\frac{2}{n^2h_n^2}\sum_{i,j,l,u}\cov\big(\Delta_i^n \epsilon \Delta_{i+l}^n \epsilon ,\Delta_j^n \epsilon\Delta_{j+u}^n \epsilon|{\cal{F}}_{kh_n}\big)+\KLEINO_{\P}(n^{-1/2})\\
&=\mathcal{O}_{\P}\big((nh_n)^{-1}Q^3\big)=\mathcal{O}_{\P}\big(n^{-1/2}\big)\,,
\end{align*}
since the covariances vanish whenever the difference of two indices exceeds $Q<\infty$. Analogously, we derive that $\var\big(\tilde Z_{kh_n}^{(u)}|{\cal{F}}_{kh_n}\big)=\mathcal{O}_{\P}\big(n^{-1/2}\big)$ for all $k$. This readily implies that $\var(\hat\eta_{kh_n})=\mathcal{O}\big(n^{-1/2}\big)$, and with Chebyshev's inequality and \eqref{expnoise} we conclude that $\hat\eta_{kh_n}=\eta_{kh_n}+\mathcal{O}_{\P}(n^{-1/4})$.\\
It remains to prove \eqref{noiseesteq} for $\alpha\ge 1/2$. Then, $\sup_{t\in[\tilde k M_n/n,(\tilde k+1)M_n/n]}|\eta_{\tilde k M_n/n}-\eta_t|=\mathcal{O}\big((M_n/n)^{\alpha}\big)=\mathcal{O}\big(n^{-\frac{\alpha}{2\alpha+1}}\big)=\KLEINO(n^{-\beta})$ by \eqref{beta}. Repeating the steps for estimates $\hat\eta_{kh_n}$ from $nh_n$ observations, we now obtain with $M_n$ observations, for all $\tilde k=0,\ldots,n/M_n-1$, that 
\begin{align*}\hat\eta_{\tilde k \frac{M_n}{n}}=\eta_{\tilde k \frac{M_n}{n}}+\mathcal{O}_{\P}(M_n^{-1/2})=\eta_{\tilde k \frac{M_n}{n}}+\mathcal{O}_{\P}\Big(n^{-\frac{1}{2}(1-\frac{1}{(2\alpha+1)})}\Big)&=\eta_{\tilde k \frac{M_n}{n}}+\mathcal{O}_{\P}\big(n^{-\frac{\alpha}{2\alpha+1}}\big)\\
&=\eta_{\tilde k \frac{M_n}{n}}+\KLEINO_{\P}(n^{-\beta})\,.\end{align*}
This proves Proposition \ref{cornoiseest}.
\hfill\qed
\subsection{Stable convergence of the spot squared volatility estimators}
We first prove two lemmas, one on moments of the noise terms in the spectral statistics and one on moments of the statistics \eqref{zeta}.
\begin{lem}\label{lemmacornoise}On Assumption \ref{eta} with $p\ge 4$:
\begin{subequations}
\begin{align}\E\Big[\Big(\sum_{i=1}^{n-1}\epsilon_i\,\varphi_{jk}\Big(\frac{i+1/2}{n}\Big)n^{-1}\Big)^2\Big|\mathcal{F}_{kh_n}\Big]=\eta_{kh_n} n^{-1}+\KLEINO_{\P}\big(n^{-1}\big)\,,\end{align}
\begin{align}\E\Big[\Big(\sum_{i=1}^{n-1}\epsilon_i\,\varphi_{jk}\Big(\frac{i+1/2}{n}\Big)n^{-1}\Big)^4\Big|\mathcal{F}_{kh_n}\Big]=3\,\eta_{kh_n}^2 n^{-2}+\KLEINO_{\P}\big(n^{-2}\big)\,.\end{align}
%If we suppose the existence of $2p$th moments $\E[\epsilon_i^{2p}]$ for all $i=0,\ldots,n$ 
Under Assumption \ref{eta} with $p=2p',p'>2$, it holds that
\begin{align}\E\Big[\Big(\sum_{i=1}^{n-1}\epsilon_i\varphi_{jk}\Big(\frac{i+1/2}{n}\Big)n^{-1}\Big)^{2p'}\Big|\mathcal{F}_{kh_n}\Big]\le K_{p'}\,\eta_{kh_n}^{p'} n^{-p'}+\KLEINO_{\P}\big(n^{-p'}\big)\,.\end{align}
\end{subequations}
\end{lem}
\begin{proof}
\begin{align*}\E\Big[\Big(\sum_{i=1}^{n-1}\epsilon_i\varphi_{jk}\Big(\frac{i+1/2}{n}\Big)n^{-1}\Big)^2\Big|\mathcal{F}_{kh_n}\Big]&=\sum_{i=1}^{n-1}\sum_{l=-i}^{n-i-1}\E[\epsilon_i\epsilon_{i+l}|\mathcal{F}_{kh_n}] \frac{\varphi_{jk}\Big(\frac{i+1/2}{n}\Big)\varphi_{jk}\Big(\frac{i+l+1/2}{n}\Big)}{n^2}\\
&=\big(\eta_{kh_n}+\KLEINO_{\P}(1)\big)n^{-1}\sum_{i=1}^{n-1}\varphi_{jk}^2\Big(\frac{i+1/2}{n}\Big)n^{-1}+R_n\\
&=\eta_{kh_n}\,n^{-1}+\KLEINO_{\P}(n^{-1})+R_n\,.\end{align*}
To control the remainder $R_n$, we perform a Taylor expansion
\[\varphi_{jk}\Big(\frac{i+l+1/2}{n}\Big)-\varphi_{jk}\Big(\frac{i+1/2}{n}\Big)=\sum_{r\ge 1}\frac{\varphi_{jk}^{(r)}\Big(\frac{i+1/2}{n}\Big)}{r !}\,\frac{l^r}{n^r}\,,\]
with $\varphi_{jk}^{(r)}$ the existing $r$th derivative of $\varphi_{jk}$.
If $|\E[\epsilon_i\epsilon_{i+l}]|\le |l|^{-1-\varrho}$ for all $i=0,\ldots,n-l$, it follows for any $i= k n h_n +1,\ldots, (k+1) n h_n$ that
\begin{align*}\sum_{l= k n h_n +1-i}^{ (k+1) n h_n-i}\E[\epsilon_i\epsilon_{i+l}|\mathcal{F}_{kh_n}]\frac{l^r}{n^rh_n^r}=\mathcal{O}\Big(\sum_{l=1}^{nh_n}\frac{l^{r-1-\varrho}}{n^rh_n^r}\Big)=\mathcal{O}\big((nh_n)^{-\varrho}\big)\end{align*} 
which tends to zero and is $\KLEINO(n^{-\beta/2})$ when $\varrho>\beta$. Since $\varphi_{jk}^{(r)}\lesssim h_n^{-r}\,\varphi_{jk}$ and $\varphi_{jk}$ is zero outside the interval $[kh_n,(k+1)h_n]$ it follows that $R_n=\mathcal{O}_{\P}(n^{-1}(nh_n)^{-\varrho})$.\\
Considering fourth moments yields
\begin{align*}\E\Big[\Big(\sum_{i=1}^{n-1}\epsilon_i\,\varphi_{jk}\Big(\frac{i+1/2}{n}\Big)n^{-1}\Big)^4\Big|\mathcal{F}_{kh_n}\Big]&=\sum_{i,u,l,v}n^{-4}\E[\epsilon_i\epsilon_{l}\epsilon_u\epsilon_v|\mathcal{F}_{kh_n}] \,\varphi_{jk}\Big(\frac{i+1/2}{n}\Big)\varphi_{jk}\Big(\frac{l+1/2}{n}\Big)\\
&\quad\quad\quad\quad\quad\quad\quad\quad \quad\quad \times \varphi_{jk}\Big(\frac{u+1/2}{n}\Big)\varphi_{jk}\Big(\frac{v+1/2}{n}\Big)\\
&=3n^{-4}\sum_{i=1}^{n-1}\sum_{l=-i}^{n-i+1}\E[\epsilon_i\epsilon_{i+l}|\mathcal{F}_{kh_n}]\,(1+\KLEINO_{\P}(1))\\%\varphi_{jk}\Big(\frac{i+1/2}{n}\Big)\varphi_{jk}\Big(\frac{i+l+1/2}{n}\Big)\\
& \hspace*{-6cm}\times \varphi_{jk}\Big(\frac{i+1/2}{n}\Big)\varphi_{jk}\Big(\frac{i+l+1/2}{n}\Big)\sum_{u=1}^{n-1}\sum_{v=-u}^{n-u+1}\E[\epsilon_u\epsilon_{u+v}|\mathcal{F}_{kh_n}]\varphi_{jk}\Big(\frac{u+1/2}{n}\Big)\varphi_{jk}\Big(\frac{u+v+1/2}{n}\Big)\\%(1+\KLEINO_{\P}(1))\\
&= 3n^{-2}\eta_{kh_n}^2  \Bigg( \sum_{u=1}^{n-1}\frac{\varphi_{jk}^2\Big(\frac{u+1/2}{n}\Big)}{n}\Bigg)^2+ \tilde R_n\,.\end{align*}
%\end{align*}
The conditional expectation $\E[\epsilon_i\epsilon_{l}\epsilon_u\epsilon_v|\mathcal{F}_{kh_n}] $ is negligible unless $|i-l|$ and $|u-v|$ are small, or $|i-u|$ and $|l-v|$ are small, or $|i-v|$ and $|u-l|$ are small. In the first identity, we have neglected the sum over terms where all four indices are close which is of the order
\[(1+\KLEINO(1))\sum_{i=1}^n\frac{\varphi_{jk}^4\Big(\frac{i+1/2}{n}\Big)}{n}\,n^{-3}\cdot \mathcal{O}_{\P}(1)=\mathcal{O}_{\P}(h_n^{-1}n^{-3})=\KLEINO_{\P}(n^{-2})\]
given that $\E[\epsilon_i^4|\mathcal{F}^X]<\infty$ almost surely for all $i$. That no fourth moments of the noise appear in the leading term is natural, as in standard proofs of central limit theorems using a moment method, since there are only $n$ addends with $i=l=u=v$. %The terms where all four indices are close to each other are negligible due to the factor $n^{-4}$ given that fourth moments of the noise process exist. Therefore
%\begin{align*}\E\Big[\Big(\sum_{i=1}^{n-1}\epsilon_i\,\varphi_{jk}\Big(\frac{i+1/2}{n}\Big)n^{-1}\Big)^4\Big|\mathcal{F}_{kh_n}\Big]&=3n^{-4}\sum_{i=1}^{n-1}\sum_{l=-i}^{n-i+1}\E[\epsilon_i\epsilon_{i+l}|\mathcal{F}_{kh_n}]\,(1+\KLEINO_{\P}(1))\\%\varphi_{jk}\Big(\frac{i+1/2}{n}\Big)\varphi_{jk}\Big(\frac{i+l+1/2}{n}\Big)\\
%& \hspace*{-6cm}\times \varphi_{jk}\Big(\frac{i+1/2}{n}\Big)\varphi_{jk}\Big(\frac{i+l+1/2}{n}\Big)\sum_{u=1}^{n-1}\sum_{v=-u}^{n-u+1}\E[\epsilon_u\epsilon_{u+v}|\mathcal{F}_{kh_n}]\varphi_{jk}\Big(\frac{u+1/2}{n}\Big)\varphi_{jk}\Big(\frac{u+v+1/2}{n}\Big)\\%(1+\KLEINO_{\P}(1))\\
%&= 3n^{-2}\eta_{kh_n}^2  \Bigg( \sum_{u=1}^{n-1}\frac{\varphi_{jk}^2\Big(\frac{u+1/2}{n}\Big)}{n}\Bigg)^2+ \tilde R_n\,.\end{align*}
That the remainder term $\tilde R_n$ is asymptotically negligible follows with the Taylor expansion from above.\\
Analogously, given that $2p'$th moments of the noise process exist for some $p'>2$, an analogous computation yields that
\begin{align*}\E\Big[\Big(\sum_{i=1}^{n-1}\epsilon_i\varphi_{jk}\Big(\frac{i+1/2}{n}\Big)n^{-1}\Big)^{2p'}\Big|\mathcal{F}_{kh_n}\Big]&=\big( (2p'-1)\cdot (2p'-3)\cdot \ldots \cdot 1\big)\frac{\eta_{kh_n}^{p'}}{n^{p'}}  \Bigg( \sum_{u=1}^{n-1}\frac{\varphi_{jk}^2\Big(\frac{u+1/2}{n}\Big)}{n}\Bigg)^{p'}\\
&\quad + \KLEINO_{\P}(n^{-p'})\,.\end{align*}
\end{proof}
\begin{lem}\label{moments}
%\begin{subequations}
On Assumptions \ref{eff}, \ref{sigma}, \ref{H} and \ref{eta}, we obtain the moment bounds 
\begin{align}\label{m1}\E\big[|\zeta_k(C+\epsilon)|^p|\mathcal{F}_{kh_n}\big]&\le K_p \big(\log(n)+\mathcal{O}_{\P}(1)\big)\,.
%\begin{align}\label{m2}\E\big[\zeta_k(C+\epsilon)|^p|\mathcal{F}_{kh_n}\big]\le K_p \log(n)
%\label{m3}\E\Big[\Big|\sum_{j=1}^{J_n}w_{jk}\|\Phi_{jk}\|_n^{-2}\big(\sum_{i=1}^n\Delta_i^n J\Phi_{jk}(i/n)\big)^2\Big|^p|\mathcal{F}_{kh_n}\Big]&\le K_p \,.
\end{align}
%\end{subequations}
\end{lem}
\begin{proof}
First, \eqref{zeta} is a convex combination and applying Jensen's inequality (for convex combinations) and Young's inequality, we derive that
\begin{align*}\E\big[|\zeta_k(C+\epsilon)|^p|\mathcal{F}_{kh_n}\big]&\le \sum_{j=1}^{J_n} w_{jk}\E\big[\big|S_{jk}^2(C+\epsilon)-\|\Phi_{jk}\|_n^{-2}\frac{\hat\eta_{kh_n}}{n}\big|^p |\mathcal{F}_{kh_n}\big]\\
&\le \sum_{j=1}^{J_n} w_{jk}\,2^{p-1}\E\big[|S_{jk}(C+\epsilon)|^{2p}+\big|\|\Phi_{jk}\|_n^{-2}\frac{\hat\eta_{kh_n}}{n}\big|^{p} |\mathcal{F}_{kh_n}\big]\,.\end{align*}
For the second addends, we obtain with \eqref{orderweights} and $J_n=\mathcal{O}(\log(n))$ that
\begin{align*}\sum_{j=1}^{J_n}w_{jk}\,2^{p-1}\big|\|\Phi_{jk}\|_n^{-2}\frac{\eta_{kh_n}}{n}\big|^p&\le K_p\sum_{j=1}^{J_n}w_{jk}\Big(\frac{j^2}{nh_n^2}\Big)^p\\
&\le K_p \Big(\sum_{j=1}^{J_n}j^{2p}(\log(n))^{-2p}\Big)\le K_p\,\log(n)\,.\end{align*}
With Proposition \ref{cornoiseest} this bound applies to the conditional expectation with $\hat\eta_{kh_n}$ also.\\
For the term with spectral statistics $S_{jk}(C+\epsilon)$, depending on the process $(C_t)_{t\in[0,1]}$ and the noise, we infer with Young's inequality and since $\E[\Delta_i^n C\Delta_l^n C]=\mathcal{O}(n^{-2})$ for all $i\ne l$, that
\begin{align*}\E\big[|S_{jk}(C+\epsilon)|^{2p} |\mathcal{F}_{kh_n}\big]&\le 2^{2p-1}\Big(\E\Big[\Big(\|\Phi_{jk}\|_n^{-2}\sum_{i=1}^n(\Delta_i^n C)^2 \Phi_{jk}^2\big(\tfrac{i}{n}\big)\Big)^{p} \,\big|\mathcal{F}_{kh_n}\Big]\\
&\quad +\E\Big[\Big|\|\Phi_{jk}\|_n^{-1}\sum_{i=1}^{n-1}\epsilon_i\,\varphi_{jk}\Big(\frac{i+1/2}{n}\Big)\frac{1}{n}\Big|^{2p}\ \big|\mathcal{F}_{kh_n}\Big]\Big)\,.
\end{align*}
Applying Jensen's inequality again yields for the first addends
\begin{align*}\E\Big[ \Big(\|\Phi_{jk}\|_n^{-2}\sum_{i=1}^n(\Delta_i^n C)^2 \Phi_{jk}^2\big(\tfrac{i}{n}\big)\Big)^{p} \,\big|\mathcal{F}_{kh_n}\Big]&\le \|\Phi_{jk}\|_n^{-2}\sum_{i=1}^n\frac{\Phi_{jk}^2\big(\tfrac{i}{n}\big)}{n}\,n^{p}\E\big[(\Delta_i^n C)^{2p}|\mathcal{F}_{kh_n}\Big] \le K_p\,\end{align*}
by \eqref{e1}. For the noise term, Lemma \ref{lemmacornoise} implies that
\begin{align*}\E\Big[\Big|\|\Phi_{jk}\|_n^{-1}\sum_{i=1}^{n-1}\epsilon_i\,\varphi_{jk}\Big(\frac{i+1/2}{n}\Big)\frac{1}{n}\Big|^{2p}\ \big|\mathcal{F}_{kh_n}\Big]&\le K_p\|\Phi_{jk}\|_{n}^{-2p}\,\eta_{kh_n}^{p}(1+\KLEINO_{\P}(1))\,n^{-p}\\
&\le K_p\,\Big(\frac{j^2}{nh_n^2}\Big)^p (1+\KLEINO_{\P}(1))\le K_p (1+\KLEINO_{\P}(1))\,\end{align*}
for all $j=1,\ldots,J_n=\mathcal{O}(\log(n))$. Inserting the bounds above yields \eqref{m1}.%\\
%For the continuous part we directly exploited a bound for the general moments of $\Delta_i^n C$. In order show \eqref{m2} instead, it is important that $\Delta_i^n J$ vanishes with high probability for most $i=\lfloor k h_n\rfloor+1,\ldots,\lfloor (k+1) h_n\rfloor$. Using $\|\Phi_{jk}\|_n^{-1}\Phi_{jk}(t)\le \sqrt{2/ h_n}$ uniformly, we deduce the simple bound
%\begin{align}\E\Big[\Big|\sum_{j=1}^{J_n}w_{jk}\|\Phi_{jk}\|_n^{-2}\big(\sum_{i=1}^n\Delta_i^n J\Phi_{jk}(i/n)\big)^2\Big|^p|\mathcal{F}_{kh_n}\Big]&\le K_p h_n^{-p}\,.\end{align}
%\E\Bigg|\sum_{j=1}^{J_n}w_{jk}\|\Phi_{jk}\|_n^{-2}\Big(\sum_{i=1}^n\Delta_i^n J\Phi_{jk}\big(\tfrac{i}{n}\big)\Big)^2\Bigg|^p \le \sum_{j=1}^{J_n}w_{jk}\E\left[\Big(\sum_{i=1}^n\|\Phi_{jk}\|_n^{-2}(\Delta_i^n J)^2\Phi_{jk}^2\big(\tfrac{i}{n}\big)\Big)^p\right]+\KLEINO(1)\\
%&\quad\quad  \le \sum_{j=1}^{J_n}w_{jk}\sum_{i=1}^n\|\Phi_{jk}\|_n^{-2}\frac{\Phi_{jk}^2\big(\tfrac{i}{n}\big)}{n}\,n^p\E\left[|\Delta_i^n J|^{2p}\right]+\KLEINO(1)=\mathcal{O}(1)\,.\end{align*}
\end{proof}
\noindent
\textbf{Proof of Theorem \ref{cltspot}}\\[.1cm]
The proof is structured in five steps. We establish the marginal stable central limit theorem for the estimator \eqref{rspot}. Since we may consider the continuous martingale part of $X$ time-reversed, the mathematical analysis for the second component follows the same arguments and we restrict ourselves to the right-limit case explicitly. Then, we address the joint convergence in the fifth step of the proof. The Steps 1-4 are structured according to the following decomposition:
\begin{align*}n^{\beta/2}\Big(\hat \sigma^2_{s}-\sigma^2_{s}\Big)&=n^{\beta/2}\Bigg(\Bigg(\sum_{k=\lfloor s h_n^{-1}\rfloor +1}^{\lfloor s h_n^{-1}\rfloor +r_n^{-1}}r_n\zeta_k^{ad}(Y)\ind\Bigg)-\sigma^2_{s}\Bigg)\\
&=n^{\beta/2}\Bigg(\Bigg(\sum_{k=\lfloor s h_n^{-1}\rfloor +1}^{\lfloor s h_n^{-1}\rfloor +r_n^{-1}}r_n\zeta_k(\tilde C^n+\epsilon)\Bigg)-\sigma^2_{s}\Bigg)\\
&\quad\quad+ n^{\beta/2}\Bigg(\sum_{k=\lfloor s h_n^{-1}\rfloor +1}^{\lfloor s h_n^{-1}\rfloor +r_n^{-1}}r_n\big(\zeta_k( C+\epsilon)-\zeta_k(\tilde C^n+\epsilon)\big)\Bigg)\\
&\quad\quad+ n^{\beta/2}\Bigg(\sum_{k=\lfloor s h_n^{-1}\rfloor +1}^{\lfloor s h_n^{-1}\rfloor +r_n^{-1}}r_n\Big(\zeta_k(Y)\1_{\{h_n|\zeta_k(Y)|\le u_n\}}-\zeta_k(C+\epsilon)\Big)\Bigg)\\
&\quad\quad+ n^{\beta/2}\Bigg(\sum_{k=\lfloor s h_n^{-1}\rfloor +1}^{\lfloor s h_n^{-1}\rfloor +r_n^{-1}}r_n\Big(\zeta_k^{ad}(Y)\ind-\zeta_k(Y)\1_{\{h_n|\zeta_k(Y)|\le u_n\}}\Big)\Bigg)\,.\end{align*}
In Step 1, we establish the stable limit theorem for the oracle spectral estimator \eqref{rspotor} built from observations of the process $\tilde C^n$ in the simplified model with noise. 
Working more generally than under Assumption \ref{eta2} with Proposition \ref{cornoiseest}, just suppose that we have some estimator
\begin{align}\label{expnoiseprep}\E\big[\hat\eta_{kh_n}|{\cal{F}}_{kh_n}\big]=\frac{1}{nh_n}\sum_{i=knh_n}^{(k+1)nh_n-1}\sum_{u=knh_n-i}^{(k+1)nh_n-1-i}\cov\big(\epsilon_i,\epsilon_{i+u}|{\cal{F}}_{kh_n}\big)+\KLEINO_{\P}(n^{-\beta/2})\,,\end{align}
as well as
\begin{align}\label{varnoiserep}\var(\hat\eta_{kh_n}|{\cal{F}}_{kh_n})=\KLEINO_{\P}(n^{-\beta})\,.\end{align}
Then, on Assumptions \ref{eff}, \ref{eta} with $p=8,\varrho>\beta$, \ref{H} with $r<2$ and \ref{sigma} and if $0<\beta<\alpha/(2\alpha+1)$, as $n\rightarrow\infty$: 
\begin{align}\label{stableCLT}&{\bf{Step~1:}}&n^{\beta/2}\Bigg(\Bigg(\sum_{k=\lfloor s h_n^{-1}\rfloor +1}^{\lfloor s h_n^{-1}\rfloor+r_n^{-1}} r_n\zeta_k(\tilde C^n+\epsilon)\Bigg)-\sigma^2_{s}\Bigg)\stackrel{(st)}{\longrightarrow}MN(0,8\sigma^3_s\eta_{s}^{1/2})\,.\end{align}
\textit{Proof of Step 1:}
In order to prove a point-wise central limit theorem we verify three conditions: one addressing the conditional bias, one the variance and one Lindeberg-type criterion.
Additionally we have to show that the convergence holds \emph{stably} in law.\\
First, we establish asymptotic unbiasedness of the local estimates \eqref{zeta}: 
\begin{align}\label{bias}\E[\zeta_k(\tilde C^n+\epsilon)|\mathcal{F}_{kh_n}]=\sigma^2_{kh_n}+\KLEINO_{\P}\big(n^{-\beta/2}\big)~~\mbox{for all $k$}\,.\end{align} 
Using the summation by parts identity \eqref{sbp}, we decompose 
\begin{align*}&\E[\zeta_k(\tilde C^n+\epsilon)|\mathcal{F}_{kh_n}]=\E\Big[\sum_{j=1}^{J_n}w_{jk}\Big(S_{jk}^2(\tilde C^n+\epsilon)-\bk\Big)\big|\mathcal{F}_{kh_n}\Big]\\
&\,=\sum_{j=1}^{J_n}w_{jk}\|\Phi_{jk}\|_n^{-2}\bigg(\E\Big[\Big(\sum_{i=1}^n\Delta_i^n\tilde C^{n}\Phi_{jk}\Big(\frac{i}{n}\Big)\Big)^2-2\sum_{i=1}^n\Delta_i^n\tilde C^{n}\Phi_{jk}\Big(\frac{i}{n}\Big)\sum_{l=1}^{n-1}\epsilon_l\varphi_{jk}\Big(\frac{l+1/2}{n}\Big)\frac{1}{n}\big|\mathcal{F}_{kh_n}\Big]\\
&\hspace*{7.35cm}+\E\Big[\Big(\sum_{i=1}^{n-1}\epsilon_i\varphi_{jk}\Big(\frac{i+1/2}{n}\Big)\frac{1}{n}\Big)^2-\frac{\hat\eta_{kh_n}}{n}\big|\mathcal{F}_{kh_n}\Big]\bigg)
\end{align*}
and consider the three terms separately. For the first term we obtain with the martingale property that
\begin{align*}\sum_{j=1}^{J_n}w_{jk}\|\Phi_{jk}\|_n^{-2}\E\Big[\Big(\sum_{i=1}^n\Delta_i^n\tilde C^{n}\Phi_{jk}\Big(\frac{i}{n}\Big)\Big)^2\big|\mathcal{F}_{kh_n}\Big]=\sum_{j=1}^{J_n}w_{jk}\|\Phi_{jk}\|_n^{-2}\sum_{i=1}^n\frac{\sigma^2_{kh_n}}{n}\Phi_{jk}^2\Big(\frac{i}{n}\Big)=\sigma^2_{kh_n}.
\end{align*}
For the noise and bias-correction term, we obtain with the bound for the remainder from Lemma \ref{lemmacornoise} and with \eqref{expnoiseprep} that
\begin{align*}&\sum_{j=1}^{J_n}w_{jk}\|\Phi_{jk}\|_n^{-2}\E\Big[\Big(\sum_{i=1}^{n-1}\epsilon_i\varphi_{jk}\Big(\frac{i+1/2}{n}\Big)\frac{1}{n}\Big)^2-\frac{\hat\eta_{kh_n}}{n}\big|\mathcal{F}_{kh_n}\Big]\\
&=\sum_{j=1}^{J_n}w_{jk}\|\Phi_{jk}\|_n^{-2}\sum_{i=knh_n}^{(k+1)nh_n-1}\sum_{l=knh_n-i}^{(k+1)nh_n-i-1}\bigg(\E[\epsilon_i\epsilon_{i+l}|\mathcal{F}_{kh_n}]n^{-2}\varphi_{jk}\Big(\frac{i+1/2}{n}\Big)\varphi_{jk}\Big(\frac{i+l+1/2}{n}\Big)\\
&\hspace*{8cm}-\frac{\E[\hat\eta_{kh_n}|\mathcal{F}_{kh_n}]}{n}\bigg)\\
&=\sum_{j=1}^{J_n}w_{jk}\|\Phi_{jk}\|_n^{-2}\sum_{i=knh_n}^{(k+1)nh_n-1}\sum_{l= knh_n-i}^{(k+1)nh_n-i-1}\E[\epsilon_i\epsilon_{i+l}|\mathcal{F}_{kh_n}]\Big(\frac{\varphi_{jk}\Big(\frac{i+1/2}{n}\Big)\varphi_{jk}\Big(\frac{i+l+1/2}{n}\Big)}{n^2}-\frac{1}{n^2h_n}\Big)\\
&\hspace*{8cm}+\KLEINO_{\P}(n^{-\beta/2})\\
&=\sum_{j=1}^{J_n}w_{jk}\|\Phi_{jk}\|_n^{-2}n^{-1}\big(\eta_{kh_n}+\KLEINO_{\P}(1)\big)\bigg(\sum_{i=knh_n}^{(k+1)nh_n-1}\frac{\varphi_{jk}^2\Big(\frac{i+1/2}{n}\Big)}{n}-1+\KLEINO(n^{-\beta/2})\bigg)+\KLEINO_{\P}(n^{-\beta/2})\\
&=\KLEINO_{\P}(n^{-\beta/2})\,,\end{align*}
by \eqref{squarephi} since $\|\Phi_{jk}\|_n^{-2}n^{-1}$ is uniformly bounded. %We used the smoothness of $(\varphi_{jk})$ and \eqref{expnoiseprep}.\\
The expectation of cross terms clearly vanishes under independence of noise and $(X_t)$. Under \eqref{endo}, we derive that
\begin{align*}&\sum_{j=1}^{J_n}w_{jk}\|\Phi_{jk}\|_n^{-2}\E\Big[\sum_{i=1}^n\Delta_i^n\tilde C^{n}\Phi_{jk}\Big(\frac{i}{n}\Big)\sum_{l=1}^{n-1}\epsilon_l\varphi_{jk}\Big(\frac{l+1/2}{n}\Big)\frac{1}{n}\big|\mathcal{F}_{kh_n}\Big]\\
&=\sum_{j=1}^{J_n}w_{jk}\|\Phi_{jk}\|_n^{-2}\sum_{i=knh_n}^{(k+1)nh_n-1}\sum_{l=knh_n}^{i}\E[\epsilon_i\Delta_l^n\tilde C^{n}|\mathcal{F}_{kh_n}\Big]\Phi_{jk}\Big(\frac{l}{n}\Big)\varphi_{jk}\Big(\frac{i+1/2}{n}\Big)\frac{1}{n}\\
&=\sum_{j=1}^{J_n}w_{jk}\|\Phi_{jk}\|_n^{-2}\big(\rho_{kh_n}+\KLEINO_{\P}(1)\big)\sum_{i=knh_n}^{(k+1)nh_n-1}\frac{1}{n}\Big(\Phi_{jk}\Big(\frac{i}{n}\Big)+\mathcal{O}((h_nn)^{-1})\Big)\Big(\varphi_{jk}\Big(\frac{i}{n}\Big)+\mathcal{O}(n^{-1})\Big)\\
&=\mathcal{O}_{\P}((nh_n)^{-1}+h_n^{3/2})\Big)=\KLEINO_{\P}(n^{-\beta/2})\,,\end{align*}
since $\sum_{i=knh_n}^{(k+1)nh_n-1}\Phi_{jk}\big(\frac{i}{n}\big)\varphi_{jk}\big(\frac{i}{n}\big)=0$ and using the bound $|\int\Phi_{jk}(t)dt|\le 2\sqrt{2}h_n^{3/2}j^{-2}$, whereas $(\varphi_{jk})$ integrate to zero. To put it simply, that the integrals in \eqref{crossphi} vanish for $j=u$ guarantees that the endogenous noise does not induce any non-negligible bias term. This completes the proof of \eqref{bias}.\\
For the expectation of the left-hand side in \eqref{stableCLT}, we deduce that
\begin{align*}&n^{\beta/2}\Bigg(\sum_{k=\lfloor s h_n^{-1}\rfloor +1}^{\lfloor s h_n^{-1}\rfloor +r_n^{-1}}r_n\E[\zeta_k(\tilde C^n+\epsilon)-\sigma^2_{s}|\mathcal{F}_{kh_n}]\Bigg)=n^{\beta/2}\Bigg(\sum_{k=\lfloor s h_n^{-1}\rfloor +1}^{\lfloor s h_n^{-1}\rfloor +r_n^{-1}}r_n (\sigma^2_{kh_n}-\sigma^2_s)\Bigg)\\
&\quad=\mathcal{O}_{\P}\Bigg(n^{\beta/2}r_n\sum_{k=\lfloor s h_n^{-1}\rfloor +1}^{\lfloor s h_n^{-1}\rfloor +r_n^{-1}} (kh_n)^{\alpha}\Bigg)=\mathcal{O}_{\P}\Big(n^{\beta/2}(h_n/r_n)^{\alpha}\Big)\\
&\quad =\mathcal{O}_{\P}\Big(n^{\beta (\alpha+1/2)}n^{-\alpha/2}\log^{\alpha}(n)\Big)=\KLEINO_{\P}(1)\,,\end{align*}
because $\alpha>0$ and $\beta<\alpha(2\alpha+1)^{-1}$. % for any $\alpha$ implying $\beta<1/4$ for $\alpha<1/2$.\\
By \eqref{varnoiserep} and using that $\|\Phi_{jk}\|_n^{-2}n^{-1}$ is uniformly bounded for all $j$, we obtain that
\begin{align*}
\var\Big(\sum_{j=1}^{J_n}w_{jk}\|\Phi_{jk}\|_n^{-2}\frac{\hat\eta_{kh_n}}{n}\big|\mathcal{F}_{kh_n}\Big)=\Big(\sum_{j=1}^{J_n}w_{jk}n^{-1}\|\Phi_{jk}\|_n^{-2}\Big)^2\var\big(\hat\eta_{kh_n}|\mathcal{F}_{kh_n}\big)=\KLEINO_{\P}(n^{-\beta})\,.
\end{align*}
Thus, the estimation of $\eta_{kh_n}$ in the bias-correction is negligible in the variance of $\hat\sigma_s^2$. In case of exogenous noise, with Lemma \ref{lemmacornoise}, we can readily adopt the identity 
\begin{align}\label{adoptvar}\var\Big(\zeta_k(\tilde C^n+\epsilon)|\mathcal{F}_{kh_n}\Big)=\sum_{j=1}^{J_n}w_{jk}^2I_{jk}^{-1}=I_k^{-1}\end{align}
from Section 6.2.2 of \cite{stable} with $I_k,I_{jk}$ from \eqref{orweights}. We consider additionally the conditional variance terms due to endogenous noise under condition \eqref{endo}. With similar estimates for the remainders as in the bias term above, we obtain that
\begin{align*}&\cov\Big(\Big(\sum_{i=1}^n\Delta_i^n\tilde C^{n}\Phi_{jk}\Big(\frac{i}{n}\Big)\Big)^2\,,\,\Big(\sum_{l=1}^{n-1}\epsilon_l\varphi_{uk}\Big(\frac{l+1/2}{n}\Big)\frac{1}{n}\Big)^2\big|\mathcal{F}_{kh_n}\Big)\\
&=\sum_{i,p=1}^n\sum_{l,q=1}^{n-1}\Big(\E\big[\Delta_i^n\tilde C^{n}\Delta_p^n\tilde C^{n}\epsilon_l\epsilon_q\big|\mathcal{F}_{kh_n}\big]-\E\big[\Delta_i^n\tilde C^{n}\Delta_p^n\tilde C^{n}\big|\mathcal{F}_{kh_n}\big]\E\big[\epsilon_l\epsilon_q\big|\mathcal{F}_{kh_n}\big]\Big)\\
&\hspace*{4cm}\times n^{-2}\Phi_{jk}\Big(\frac{i}{n}\Big)\Phi_{jk}\Big(\frac{p}{n}\Big)\varphi_{uk}\Big(\frac{l+1/2}{n}\Big)\varphi_{uk}\Big(\frac{q+1/2}{n}\Big)\\
&=2\sum_{l=knh_n}^{(k+1)nh_n-1}\sum_{i=knh_n}^{l}\E\big[\epsilon_l\Delta_i^n\tilde C^{n}|\mathcal{F}_{kh_n}\big]\sum_{q=knh_n}^{(k+1)nh_n-1}\sum_{p=knh_n}^{q}\E\big[\epsilon_q\Delta_p^n\tilde C^{n}|\mathcal{F}_{kh_n}\big]\\
&\hspace*{4cm}\times n^{-2}\Phi_{jk}\Big(\frac{i}{n}\Big)\Phi_{jk}\Big(\frac{p}{n}\Big)\varphi_{uk}\Big(\frac{l+1/2}{n}\Big)\varphi_{uk}\Big(\frac{q+1/2}{n}\Big)(1+\KLEINO_{\P}(1))\\
&=2\sum_{l=knh_n}^{(k+1)nh_n-1}\hspace*{-.15cm}\big(\rho_{kh_n}+\KLEINO_{\P}(1)\big)n^{-1}\Phi_{jk}\Big(\frac{l}{n}\Big)\varphi_{uk}\Big(\frac{l}{n}\Big)\sum_{q=knh_n}^{(k+1)nh_n-1}\hspace*{-.15cm}\big(\rho_{kh_n}+\KLEINO_{\P}(1)\big)n^{-1}\Phi_{jk}\Big(\frac{q}{n}\Big)\varphi_{uk}\Big(\frac{q}{n}\Big)\\
&\hspace*{10cm}+\mathcal{O}_{\P}\big(h_n^{3/2}n^{-1}\big)\\
&=2\rho^2_{kh_n}\Big(\int_0^1\tilde \Phi_{jk}(t)\varphi_{uk}(t)\,dt\Big)^2(1+\KLEINO_{\P}(1))+\mathcal{O}_{\P}\big(h_n^{3/2}n^{-1}\big)\,.
\end{align*}
In the first identity the terms for $i=p$ and $i$ not close to $l,q$ cancel. We used the smoothness of $(\Phi_{jk})$ and $(\varphi_{jk})$ again.
Analogously, we obtain that
\begin{align*}&\cov\Big(\Big(\sum_{l=1}^{n-1}\epsilon_l\varphi_{jk}\Big(\frac{l+1/2}{n}\Big)\frac{1}{n}\Big)^2\,,\,\Big(\sum_{i=1}^n\Delta_i^n\tilde C^{n}\Phi_{uk}\Big(\frac{i}{n}\Big)\Big)^2\big|\mathcal{F}_{kh_n}\Big)\\
&=2\rho^2_{kh_n}\Big(\int_0^1\tilde \Phi_{uk}(t)\varphi_{jk}(t)\,dt\Big)^2(1+\KLEINO_{\P}(1))+\mathcal{O}_{\P}\big(h_n^{3/2}n^{-1}\big)\,.
\end{align*}
With similar computations, we obtain that
\begin{align*}&\cov\Big(\sum_{i=1}^n\Delta_i^n\tilde C^{n}\Phi_{jk}\Big(\frac{i}{n}\Big)\sum_{l=1}^{n-1}\epsilon_l\varphi_{jk}\Big(\frac{l+1/2}{n}\Big)\frac{1}{n}\,,\,\Big(\sum_{p=1}^n\Delta_p^n\tilde C^{n}\Phi_{uk}\Big(\frac{p}{n}\Big)\Big)^2\big|\mathcal{F}_{kh_n}\Big)\\
&=2\rho_{kh_n}\sigma^2_{kh_n}\hspace*{-0.05cm}\Big(\int_0^1\hspace*{-0.05cm}\tilde\Phi_{jk}(t)\tilde\Phi_{uk}(t)\,dt\hspace*{-0.05cm}\int_0^1\hspace*{-0.05cm}\tilde\Phi_{uk}(t)\varphi_{jk}(t)\,dt\Big)\hspace*{-0.05cm}(1+\KLEINO_{\P}(1))\hspace*{-0.05cm}+\hspace*{-0.05cm}\mathcal{O}_{\P}\big(h_n^{3/2}n^{-1}\big)\hspace*{-0.05cm}=\hspace*{-0.05cm}\mathcal{O}_{\P}\big(h_n^{3/2}n^{-1}\big),\end{align*}
since $\int_0^1\tilde \Phi_{jk}(t)\tilde\Phi_{uk}(t)\,dt\int_0^1\tilde\Phi_{uk}(t)\varphi_{jk}(t)\,dt=0$ for all $j,u$. Analogously, since  $\int_0^1\varphi_{jk}(t)\varphi_{uk}(t)\,dt$ $\int_0^1\tilde \Phi_{jk}(t)\varphi_{uk}(t)\,dt=0$ for all $j,u$, the conditional covariance of cross terms and noise terms is of the same order $h_n^{3/2}n^{-1}$ in probability.
The only other (at first) non-negligible additional conditional variance term thus comes from
\begin{align*}&\cov\Big(\sum_{i=1}^n\Delta_i^n\tilde C^{n}\Phi_{jk}\Big(\frac{i}{n}\Big)\sum_{l=1}^{n-1}\epsilon_l\varphi_{jk}\Big(\frac{l+1/2}{n}\Big)\frac{1}{n}\,,\,\sum_{p=1}^n\Delta_p^n\tilde C^{n}\Phi_{uk}\Big(\frac{p}{n}\Big)\sum_{q=1}^{n-1}\epsilon_q\varphi_{uk}\Big(\frac{q+1/2}{n}\Big)\frac{1}{n}\big|\mathcal{F}_{kh_n}\Big)\,.
\end{align*}
Using the same approximations as in the previous terms and subtracting the term already contained in $I_{jk}^{-1}$ from the exogenous setup, we obtain the overall \emph{additional} conditional variance
\begin{align*}
&\sum_{j,u=1}^{J_n}w_{jk}w_{uk}\|\Phi_{jk}\|_n^{-2}\|\Phi_{uk}\|_n^{-2}\bigg(\cov\Big(\Big(\sum_{i=1}^n\Delta_i^n\tilde C^{n}\Phi_{jk}\Big(\frac{i}{n}\Big)\Big)^2\,,\,\Big(\sum_{l=1}^{n-1}\epsilon_l\varphi_{uk}\Big(\frac{l+1/2}{n}\Big)\frac{1}{n}\Big)^2\big|\mathcal{F}_{kh_n}\Big)\\
&\quad +\cov\Big(\Big(\sum_{l=1}^{n-1}\epsilon_l\varphi_{jk}\Big(\frac{l+1/2}{n}\Big)\frac{1}{n}\Big)^2\,,\,\Big(\sum_{i=1}^n\Delta_i^n\tilde C^{n}\Phi_{uk}\Big(\frac{i}{n}\Big)\Big)^2\big|\mathcal{F}_{kh_n}\Big)\\
&\quad +4\cov\Big(\sum_{i=1}^n\Delta_i^n\tilde C^{n}\Phi_{jk}\Big(\frac{i}{n}\Big)\sum_{l=1}^{n-1}\epsilon_l\varphi_{jk}\Big(\frac{l+1/2}{n}\Big)\frac{1}{n}\,,\,\sum_{p=1}^n\Delta_p^n\tilde C^{n}\Phi_{uk}\Big(\frac{p}{n}\Big)\sum_{q=1}^{n-1}\epsilon_q\varphi_{uk}\Big(\frac{q+1/2}{n}\Big)\frac{1}{n}\big|\mathcal{F}_{kh_n}\Big)\bigg)\\
&\quad -\sum_{j=1}^{J_n}w_{jk}^24\|\Phi_{jk}\|_n^{-2}\sigma_{kh_n}^2\frac{\eta_{kh_n}}{n}\\
&=\sum_{j,u=1}^{J_n}w_{jk}w_{uk}\|\Phi_{jk}\|_n^{-2}\|\Phi_{uk}\|_n^{-2}\rho^2_{kh_n}(1+\KLEINO_{\P}(1)) \Bigg(\mathcal{O}_{\P}\big(h_n^{3/2}n^{-1}\big)\\
&\quad+ 2\Big(\int \tilde\Phi_{jk}(t)\varphi_{uk}(t)\,dt\Big)^2+2\Big(\int \tilde\Phi_{uk}(t)\varphi_{jk}(t)\,dt\Big)^2+4\int \tilde\Phi_{jk}(t)\varphi_{uk}(t)\,dt\,\int \tilde\Phi_{uk}(t)\varphi_{jk}(t)\,dt\Bigg)\,.
\end{align*}
However, by \eqref{crossphi} the integrals sum up to zero. Since the remainder is $\mathcal{O}_{\P}\big((\log(n))^3n^{-1/4}\big)$, the effect of the endogenous noise becomes negligible at first asymptotic order. We conclude \eqref{adoptvar}.\\
In the sequel we write $w_{jk}, I_{jk}, I_k$ as functions of the squared volatility and $\eta$: $I_j(\sigma^2,\eta)=\tfrac12\big(\sigma^2+\|\Phi_{jk}\|_n^{-2}\tfrac{\eta}{n}\big)^{-2}$, $I(\sigma^2,\eta)=\sum_{j=1}^{J_n}I_j(\sigma^2,\eta)$ and $w_j(\sigma^2,\eta)=(I(\sigma^2,\eta))^{-1}I_j(\sigma^2,\eta)$. Note that $\|\Phi_{jk}\|_n^{-2}$ is equal for all $k$ such that the time-dependence of $I,I_j,w_j$ is only in the squared volatility $\sigma^2$ and $\eta$. 
For the sum of conditional variances of the left-hand side of \eqref{stableCLT}, we obtain that
\begin{align*}n^{\beta}\sum_{k=\lfloor s h_n^{-1}\rfloor +1}^{\lfloor s h_n^{-1}\rfloor +r_n^{-1}}r_n^2\,\var\Big(\zeta_k(\tilde C^n+\epsilon)|\mathcal{F}_{kh_n}\Big)&=n^{\beta}r_n\Bigg(\sum_{k=\lfloor s h_n^{-1}\rfloor +1}^{\lfloor s h_n^{-1}\rfloor +r_n^{-1}}\sum_{j=1}^{J_n}r_n w_{jk}^2I_{jk}^{-1}+\KLEINO_{\P}(1)\Bigg)\\
&=\log{(n)}I^{-1}\Big(\sigma^2_{\lfloor sh_n^{-1}\rfloor h_n},\eta_{\lfloor sh_n^{-1}\rfloor h_n}\Big)+R_n\,.\end{align*}
%The remainder in the first equality is due to fourth moments of the noise when the noise is non-Gaussian and asymptotically negligible. 
We exploit bounds on the derivative of the weights with respect to $\sigma^2$ and $\eta$
\begin{align}\label{derweights}\frac{\partial w_j(\sigma^2,\eta)}{\partial \sigma^2}=\mathcal{O}\big(w_j(\sigma^2,\eta)\log^2{(n)}\big)\,,\end{align}
here and several times below. The bound is proved as Equation (77) in \cite{stable}. $\partial w_j(\sigma^2,\eta)/{(\partial \eta)}$ can be bounded analogously. Observe that by the chain and product differentiation rule
\[\frac{\partial}{\partial\sigma^2}\big(w_j^2(\sigma^2,\eta)(I_j(\sigma^2, \eta))^{-1}\big)=2w_j(\sigma^2,\eta)\frac{\partial w_j}{\partial \sigma^2 }(\sigma^2,\eta)(I_j(\sigma^2,\eta))^{-1}+w_j^2(\sigma^2,\eta)4\big(\sigma^2+\|\Phi_{jk}\|_n^{-2}\tfrac{\eta}{n}\big)\,.\]
Thus, we can find an upper bound for the remainder $R_n$ using
\begin{align*}\sum_{j=1}^{J_n}\big(1\vee \|\Phi_{jk}\|_n^{-2}n^{-1}\big)\big(1\wedge \|\Phi_{jk}\|_n^{8}n^{4}\big)=\mathcal{O}\Bigg(\sum_{j=1}^{\lfloor\sqrt{n}h_n\rfloor}1+\sum_{j=1}^{J_n}\|\Phi_{jk}\|_n^{6}n^3\Bigg)=\mathcal{O}(\log^6{(n)})\\
\Rightarrow R_n=\mathcal{O}_{\P}\Bigg(n^{\beta}r_n^2\sum_{k=\lfloor s h_n^{-1}\rfloor +1}^{\lfloor s h_n^{-1}\rfloor +r_n^{-1}}\log^6{(n)}\Big(\sigma^2_{kh_n}-\sigma^2_{\lfloor sh_n^{-1}\rfloor h_n}\Big)\Bigg)=\mathcal{O}_{\P}\Big(\log^7{(n)}(h_n/r_n)^{\alpha}\Big)
\end{align*}
with \eqref{orderweights}, which tends to zero as $n\rightarrow\infty$ because $\alpha>0$. By \eqref{etasmooth}, the locally constant approximation of the long-run noise variance induces an error of smaller or at most equal order.\\
The Lindeberg condition is proved by the stronger Lyapunov criterion considering fourth moments:
\begin{align*}&n^{2\beta}\sum_{k=\lfloor s h_n^{-1}\rfloor +1}^{\lfloor s h_n^{-1}\rfloor +r_n^{-1}}r_n^4\;\E\left[\zeta_k^4(\tilde C^n+\epsilon)\big|\mathcal{F}_{kh_n}\right]%\le 
%n^{2\beta}\hspace*{-.05cm}\sum_{k=\lfloor s h_n^{-1}\rfloor +1}^{\lfloor s h_n^{-1}\rfloor +r_n^{-1}}\hspace*{-.25cm}r_n^4\Bigg(\sum_{j=1}^{J_n}\hspace*{-.05cm}w_{jk}\Bigg(\hspace*{-.05cm}\E\left[\Big(\tilde S_{jk}^2-\tfrac{\hat\eta}{n}\|\Phi_{jk}\|_n^{-2}\Big)^4\right]\hspace*{-.1cm}\Bigg)^{\tfrac14}\Bigg)^4\\
=\mathcal{O}_{\P}\Big(n^{-\beta}\log(n)\Big)=\KLEINO_{\P}(1)\,,\end{align*}
using Lemma \ref{moments} (replacing $C$ by $\tilde C^n$, the proof of Lemma \ref{moments} applies in the same way). We obtain the variance in \eqref{stableCLT}, since the bin-wise Fisher informations
\[I_k=\frac{1}{2}\sum_{j=1}^{J_n}\Big(\sigma^2_{kh_n}+\|\Phi_{jk}\|_n^{-2}\frac{\eta_{kh_n}}{n}\Big)^{-2}\]
satisfy the following convergences (see Section 6.2.2 of \cite{stable}):
\begin{align}\frac{1}{\log{(n)}}I_k\longrightarrow \int_0^{\infty}\frac{1}{2}\Big(\sigma^2_{kh_n}+\eta_{kh_n}\pi^2 x^2\Big)^{-2}\,dx=\Big(8\,\sigma_{s}^{3}\eta_{s}^{1/2}\Big)^{-1}\,,\end{align}
and the reciprocal of the right-hand side thus constitutes the asymptotic variance of $\hat\sigma_s^2$.\\
Finally, stability of the weak convergence is proved similarly as in Proposition 8.2 of \cite{jacodtodorov}. For later use, let us directly consider a collection of times where we consider estimates of the spot volatilities instead of only one fixed time. In particular, for our test, we shall focus on finitely many jumps of $X$ with absolute value larger than some constant. Consider a finite set $(S_p)_{1\le p\le P}$ with fix $P<\infty$ of ordered stopping times exhausting those jump arrivals of $X$ on $[0,1]$. The restriction of $\Omega$ to
\begin{align}\Omega_n=\left\{\omega\in\Omega|S_1>r_n^{-1}h_n,S_P<1-r_n^{-1}h_n,\forall p:(S_p-S_{p-1})>2r_n^{-1}h_n\right\}\end{align}
satisfies $\P(\Omega_n)\rightarrow 1$ as $n\rightarrow \infty$. Thus, we work on $\Omega_n$. We aim at establishing for
\begin{align}\hspace*{-.1cm}\alpha_n=n^{\beta/2}\Bigg(\sum_{k=\lfloor S_ph_n^{-1}\rfloor +1}^{\lfloor S_ph_n^{-1}\rfloor +r_n^{-1}}r_n\zeta_k(\tilde C^n+\epsilon)-\sigma^2_{S_p}\,,\,\sum_{k=\lfloor S_ph_n^{-1}\rfloor -r_n^{-1}}^{\lfloor S_ph_n^{-1}\rfloor -1}r_n\zeta_k(\tilde C^n+\epsilon)-\sigma^2_{S_p-}\Bigg)_{1\le p\le P}\end{align}
that $\E[Zg (\alpha_n)]\rightarrow \E[Z g(\alpha)]$ with $\alpha=\big(2\sqrt{2}\sigma_{S_p}^{3/2}\eta_{S_p}^{1/4} U_p,2\sqrt{2}\sigma_{S_p-}^{3/2}\eta_{S_p-}^{1/4} U'_p\big)_{1\le p\le P}$ for any $\mathcal{F}$-measur\-able bounded random variable $Z$ and continuous bounded function $g$ and for $(U_p,U_p')$ a sequence of standard normals defined on an exogenous space being independent of $\mathcal{F}$. This is the definition of the claimed $\mathcal{F}$-stable convergence.\\ 
The strategy is to exclude intervals on which the spot estimators are built and conditioning. Thereto, define
\[B_n=\bigcup_{p=1}^P[(S_p-(r_n^{-1}+1)h_n)\vee 0, (S_p+(r_n^{-1}+1)h_n)\wedge 1]\]
and $\tilde{\mathcal{G}}_t^n$ as the smallest filtration to which $\tilde C^n$ and $U$ are adapted and such that the $\sigma$-field generated by the Poisson measure which determines $S_1,\ldots, S_P$ lies in $\tilde{\mathcal{G}}_0^n$. Then each $\alpha_n$ is $\tilde{\mathcal{G}}_1^n$-measurable. The following decomposition of $\tilde C^n$ is well-defined:
\begin{align*}\tilde X(n)_t=\int_0^t\1_{B_n}(s)\sigma_{\lfloor sh_n^{-1}\rfloor h_n}\,dW_s\,,\,\bar X(n)_t=\tilde C^n_t-\tilde X(n)_t\,,\end{align*}
and analogously $(\tilde U_t)$ and $(\bar U_t)$. It is enough to consider $Z$ being $\tilde{\mathcal{G}}_1^n$-measurable, as we can simply substitute with $\E[Z|\tilde{\mathcal{G}}_1^n]$ otherwise. When $\mathcal{H}_n$ is the $\sigma$-field generated by $\tilde{\mathcal{G}}_0^n$, $\bar X(n)_t$ and $\bar U_t$, $\big(\mathcal{H}_n\big)_n$ is an isotonic sequence and $\bigvee_n \mathcal{H}_n=\tilde{\mathcal{G}}_1^n$. Since $\E[Z|\mathcal{H}_n]\rightarrow Z$ in $L^1(\P)$, it is enough to show
\begin{align}\E[Z\1_{\Omega_n} g (\alpha_n)]\rightarrow \E[Z\,g(\alpha)]= \E[Z]\E[g(\alpha)]\end{align}
for $Z$ $\mathcal{H}_q$-measurable for some $q$. We can use the approximation with constant $\mathcal{H}_q$-measurable squared volatilities $\sigma^2_{S_p}, \sigma^2_{S_p-}$ and with $\eta_{S_p}$ locally constant on the single intervals of $B_n$, where the errors have been bounded above. Restricted to $\Omega_n$ the vector $\alpha_n$ then includes only Brownian increments $\Delta_i^n W$ independent of the Brownian increments of $\bar X(n)_t$. Further, the noise is under Assumption \ref{eta} only short-term dependent on the past and in particular covariances of any such $Z$ and $\alpha_n$ tend to zero. Then for all $n\ge q$, conditional on $\mathcal{H}_q$, the vector $\alpha_n$ has a law asymptotically independent of $\bar X(n)_t$ and $\bar U_t$, such that the ordinary central limit theorem implies the claimed convergence. The above proof includes the stable convergence of the spot volatility estimator at one fixed time $s\in(0,1)$ as a special case. Thus, we have verified all conditions and infer the stable limit theorem \eqref{stableCLT}.\hfill\qed \\[.2cm]

To prove that the same limit theorem as \eqref{stableCLT} is valid for $n^{\beta/2}\Big(\hat \sigma^2_{s}-\sigma^2_{s}\Big)$, we show for the other addends above that they converge to zero in probability for all $s\in(0,1)$. We proceed with
\begin{align}\label{rem1}&{\bf{Step~2:}}&n^{\beta/2}\Bigg(\sum_{k=\lfloor s h_n^{-1}\rfloor +1}^{\lfloor s h_n^{-1}\rfloor+r_n^{-1}}r_n\Bigg(\zeta_k( C+\epsilon)-\zeta_k(\tilde C^n+\epsilon)\Bigg)\Bigg)=\KLEINO_{\P}(1)\,,\end{align}
under the same conditions as in Step 1. This remainder due to approximating $C$ by the simplified processes $\tilde C^n$ has exactly the same structure as the one for integrated squared volatility estimation examined in paragraph 6.3 of \cite{stable}. We just incorporate the additional jump component in the volatility using an estimate as \eqref{boundjumps}. Then, repeating the proof along the same lines, only changing the mean over all bins to the mean over local windows of size $r_n^{-1}h_n$, renders with $\beta<1/2$ the order:
\[\zeta_k(C+\epsilon)-\zeta_k(\tilde C^n+\epsilon)=\mathcal{O}_{\P}\big(h_n^{\alpha}\big)=\KLEINO_{\P}\big(n^{-\beta/2}\big)\,,\]
uniformly for all $k$. Analogously to \cite{stable}, we require here the mild condition \eqref{drift}.

\begin{align}\label{rem2}&{\bf{Step~3:}}&n^{\beta/2}\Bigg(\sum_{k=\lfloor s h_n^{-1}\rfloor +1}^{\lfloor s h_n^{-1}\rfloor+r_n^{-1}}r_n\big(\zeta_k(Y)\1_{\{h_n|\zeta_k(Y)|\le u_n\}}-\zeta_k(C+\epsilon)\big)\Bigg)=\KLEINO_{\P}(1)\,,\end{align}
when, additional to the assumptions for Steps 1 and 2, we have $\beta<\tau(1-r/2)$ and $\tau<1-\beta/(p-2)$ when $p<\infty$ moments of the noise exist in Assumption \ref{eta}.\\
\textit{Proof of Step 3:}
This part of the proof is related to Chapter 13 of \cite{JP} and the proofs in \cite{bibwink2015}. Our strategy here is related, but slightly different. We differentiate three cases. For some fixed $\rho\in(0,1)$, for instance $\rho=1/2$, consider
\begin{align*}&n^{\beta/2}\Bigg(\sum_{k=\lfloor s h_n^{-1}\rfloor +1}^{\lfloor s h_n^{-1}\rfloor+r_n^{-1}}r_n\big(\zeta_k(Y)\1_{\{h_n|\zeta_k(Y)|\le u_n\}}-\zeta_k(C+\epsilon)\big)\Bigg)=\\
&\quad n^{\beta/2}\Bigg(\sum_{k=\lfloor s h_n^{-1}\rfloor +1}^{\lfloor s h_n^{-1}\rfloor+r_n^{-1}}r_n\1_{\{h_n|\zeta_k(C+\epsilon)|>\rho u_n\}}\big(\zeta_k(Y)\1_{\{h_n|\zeta_k(Y)|\le u_n\}}-\zeta_k(C+\epsilon)\big)\Bigg)\\
&\quad\quad -n^{\beta/2}\Bigg(\sum_{k=\lfloor s h_n^{-1}\rfloor +1}^{\lfloor s h_n^{-1}\rfloor+r_n^{-1}}r_n\1_{\{h_n|\zeta_k(C+\epsilon)|\le \rho u_n\}}\1_{\{h_n|\zeta_k(Y)|> u_n\}}\zeta_k(C+\epsilon)\Bigg)\\
&\quad\quad + n^{\beta/2}\Bigg(\sum_{k=\lfloor s h_n^{-1}\rfloor +1}^{\lfloor s h_n^{-1}\rfloor+r_n^{-1}}r_n\1_{\{h_n|\zeta_k(C+\epsilon)|\le \rho u_n\}}\1_{\{h_n|\zeta_k(Y)|\le u_n\}}\big(\zeta_k(Y)-\zeta_k(C+\epsilon)\big)\Bigg)\,.\end{align*}
We prove that all three sums tend to zero in probability. For the first term, when $h_n|\zeta_k(C+\epsilon)|>\rho u_n=c\rho h_n^{\tau}$, it suffices to prove that uniformly for all $k$:
\[|\zeta_k(Y)\1_{\{h_n|\zeta_k(Y)|\le u_n\}}-\zeta_k(C+\epsilon)|=\KLEINO_{\P}\big(n^{-\beta/2}\big)\,.\]
We can choose $N_0\in\mathds{N}$, such that $h_n^{N_0(1-\tau)}=\KLEINO(n^{-\beta/2-\varepsilon})$ for some $\varepsilon>0$. 
Given that $h_n|\zeta_k(C+\epsilon)|>\rho u_n$, when we have enough moments of the noise such that $\tau<1-\beta/(p-2)$, we conclude with Lemma \ref{moments} that 
\begin{align*}&|\zeta_k(Y)\1_{\{h_n|\zeta_k(Y)|\le u_n\}}-\zeta_k(C+\epsilon)|\le \big(h_n^{-1}u_n+|\zeta_k(C+\epsilon)|\big)\Big|\frac{h_n\zeta_k(C+\epsilon)}{\rho u_n}\Big|^{N_0+1}\\
&\quad \le  \big(|\zeta_k(C+\epsilon)|^{N_0+1}+|\zeta_k(C+\epsilon)|^{N_0+2} h_n^{1-\tau}(c\rho)^{-1}\big)(c\,\rho)^{-N_0}h_n^{N_0(1-\tau)}\\ &\quad =\mathcal{O}_{\P}\Big(\log(n) h_n^{N_0(1-\tau}\Big)=\KLEINO_{\P}\big(n^{-\beta/2}\big)\,.
\end{align*}
This shows that the first sum above tends to zero in probability. Next, we prove that
\begin{align}\sum_{k=\lfloor s h_n^{-1}\rfloor +1}^{\lfloor s h_n^{-1}\rfloor+r_n^{-1}}r_n\1_{\{h_n|\zeta_k(C+\epsilon)|\le \rho u_n\}}\1_{\{h_n|\zeta_k(Y)|> u_n\}}\zeta_k(C+\epsilon)=\KLEINO_{\P}\big(n^{-\beta/2}\big)\,.\end{align}
We have the decomposition
\begin{align*}\zeta_k(Y)\hspace*{-.075cm}=\zeta_k(C+\epsilon)\hspace*{-.05cm}+\hspace*{-.05cm}\sum_{j=1}^{J_n}w_{jk}\|\Phi_{jk}\|_n^{-2}\hspace*{-.05cm}\Bigg(\hspace*{-.05cm}\Bigg(\hspace*{-.05cm}\sum_{i=1}^n\hspace*{-.05cm}\Delta_i^n J\Phi_{jk}\big(\tfrac{i}{n}\big)\hspace*{-.05cm}\Bigg)^2\hspace*{-.05cm} +\hspace*{-.05cm} 2\sum_{i=1}^n\hspace*{-.05cm}\Delta_i^n J\Phi_{jk}\big(\tfrac{i}{n}\big)\hspace*{-.1cm}\sum_{v=1}^n\hspace*{-.05cm}\Delta_v^n C\Phi_{jk}\big(\tfrac{v}{n}\big)\hspace*{-.05cm}\Bigg),\end{align*}
neglecting cross terms of jumps and noise. All cross terms can be bounded using the Cauchy-Schwarz inequality. Observe that 
\[\1_{\{h_n|\zeta_k(C+\epsilon)|\le \rho u_n\}}\1_{\{h_n|\zeta_k(Y)|> u_n\}}\le \1_{\big\{h_n\big|\sum_{j=1}^{J_n}w_{jk}\|\Phi_{jk}\|_n^{-2}(\sum_{i=1}^n\Delta_i^n J\Phi_{jk}(i/n))^2\big|> \tilde\rho \,u_n\big\}}\,,\]
for some fix $\tilde\rho\in(0,1)$ depending on $\rho$. This means that if the terms from the continuous part are not exceptionally large, the jumps need to be sufficiently large such that $h_n|\zeta_k(Y)|> u_n$. The simple uniform bound $\Phi_{jk}(t)\le \sqrt{2}h_n^{-1/2}\|\Phi_{jk}\|_n $ yields that %$\|\Phi_{jk}\|_n^{-2} \Phi_{jk}^2(t)\le 2h_n^{-1}$ yields that
\[h_n\Big|\sum_{j=1}^{J_n}w_{jk}\|\Phi_{jk}\|_n^{-2}\Big(\sum_{i=1}^n\Delta_i^n J\Phi_{jk}\big(\tfrac{i}{n}\big)\Big)^2\Big|\le 2\Big(\sum_{i=n\lfloor k h_n\rfloor +1}^{n\lfloor (k+1) h_n\rfloor}\Delta_i^n J\Big)^2\]
and we obtain that
\[\1_{\{h_n|\zeta_k(C+\epsilon)|\le \rho u_n\}}\1_{\{h_n|\zeta_k(Y)|> u_n\}}\le \1_{\{| J_{(k+1)h_n}-J_{kh_n}|> \rho^* \sqrt{u_n}\}}\,,\]
with $\rho^*=\tilde\rho/\sqrt{2}$. Therefore, it is sufficient to prove that
\begin{align}\label{jt2}\sum_{k=\lfloor s h_n^{-1}\rfloor +1}^{\lfloor s h_n^{-1}\rfloor+r_n^{-1}}\1_{\{| J_{(k+1)h_n}-J_{kh_n}|> \rho^* \sqrt{u_n}\}}=\KLEINO_{\P}\big(n^{\beta/2}\big)\,.\end{align}
Similar terms have been addressed several times in the literature, see, for instance, (13.1.14) in \cite{JP}. Applying \eqref{e5} with $\epsilon=\rho^* u_n^{1/2}$, we derive the condition
\begin{align}\label{const1}r_n^{-1} h_n u_n^{-r/2}=\KLEINO\big(n^{\beta/2}\big)~\Leftrightarrow 1-\frac{\tau r}{2}>\beta\,,\end{align}
to ensure \eqref{jt2}. When $h_n|\zeta_k(C+\epsilon)|\le \rho u_n$ and $h_n|\zeta_k(Y)|\le  u_n $, it follows that 
\[h_n\Big|\sum_{j=1}^{J_n}w_{jk}\|\Phi_{jk}\|_n^{-2}\big(\sum_{i=1}^n\Delta_i^n J\Phi_{jk}(i/n)\big)^2\Big|\le c \,u_n\] 
with some constant $c<4$.
%With the counting process
%\[N_t^n=\int_0^t\int_{\R}\1_{\{\gamma(x)>\rho* u_n^{1/2}\}}\,\mu(dx,ds)\]
In this case, we obtain by \eqref{e2}:
\begin{align*}\Big(| J_{(k+1)h_n}-J_{kh_n}|\wedge \sqrt{c\,u_n}\Big)=\mathcal{O}_{\P}\big(h_n^{1/2}\,u_n^{1-r/2}\big)\,,\end{align*}
and hence, if we can ensure that $h_n^{\tau(1-r/2)}=\KLEINO(n^{-\beta/2})$, using again $\Phi_{jk}(t)\le \sqrt{2}h_n^{-1/2}\|\Phi_{jk}\|_n $,
\begin{align*}|\zeta_k(Y)-\zeta_k(C+\epsilon)|&\le c \,\Big(\Big|\sum_{j=1}^{J_n}w_{jk}\|\Phi_{jk}\|_n^{-2}\sum_{i=1}^n\big(\Delta_i^n J\wedge \sqrt{u_n}\big)^2\Phi_{jk}^2\big(\tfrac{i}{n}\big)\Big|\wedge u_n\Big)\\
&\le 2c\,h_n^{-1}\Big(\Big(\sum_{i=n\lfloor k h_n\rfloor +1}^{n\lfloor (k+1) h_n\rfloor}\Delta_i^n J\Big)^2\wedge u_n\Big)\\
&=\mathcal{O}_{\P}\big(u_n^{1-r/2}\big)=\KLEINO_{\P}\big(n^{-\beta/2}\big)\,,\end{align*}
on the set where ${\{h_n|\zeta_k(C+\epsilon)|\le \rho u_n,h_n|\zeta_k(Y)|\le  u_n\}}$.
The condition $\beta<\tau(1-r/2)$ implies \eqref{const1} and is exactly what we need to complete the proof of \eqref{rem2}.
%%%%%%%%%%%%%%%%%%%%%%%%%%%%%%%%%%%%%%%%%%%%%%%%%%%%%%%%%%%%%%%%%%%%%%%%%%%%%%%%%%%%%
\begin{align*}%\label{rem3}
&{\bf{Step~4:}}&n^{\beta/2}\Bigg(\sum_{k=\lfloor s h_n^{-1}\rfloor +1}^{\lfloor s h_n^{-1}\rfloor+r_n^{-1}}r_n\big(\zeta_k^{ad}(Y)\ind-\zeta_k(Y)\1_{\{h_n|\zeta_p(Y)|\le u_n\}}\big)\Bigg)=\KLEINO_{\P}(1)\,.\end{align*}
%when, additional to the assumptions of Steps 1-3, Assumption \ref{eta2} holds true.\\
\textit{Proof of Step 4:}
In Step 3 we have not used the specific form of the oracle weights \eqref{orweights} and the proof analogously extends to
\begin{align}n^{\beta/2}\Bigg(\sum_{k=\lfloor s h_n^{-1}\rfloor +1}^{\lfloor s h_n^{-1}\rfloor+r_n^{-1}}r_n\big(\zeta_k^{ad}(Y)\1_{\{h_n|\zeta_k^{ad}(Y)|\le u_n\}}-\zeta_k^{ad}(C+\epsilon)\big)\Bigg)=\KLEINO_{\P}(1)\,.\end{align}
Thus it suffices to prove that
\begin{align}n^{\beta/2}\Bigg(\sum_{k=\lfloor s h_n^{-1}\rfloor +1}^{\lfloor s h_n^{-1}\rfloor+r_n^{-1}}r_n\big(\zeta_k^{ad}(C+\epsilon)-\zeta_k(C+\epsilon)\big)\Bigg)=\KLEINO_{\P}(1)\,.\label{rem3}\end{align}
We decompose this remainder as follows. Since both, oracle weights $w_j(\sigma^2_{kh_n},\eta_{kh_n})$ and estimated weights $w_j(\hat \sigma^2_{kh_n},\hat\eta_{kh_n})$ sum up to one, we can replace $(S_{jk}^2-\|\Phi_{jk}\|_n^{-2}\hat\eta_{kh_n}/n)$ by $(S_{jk}^2-\E[S_{jk}^2])$. First, consider the difference of pre-estimated and oracle weights, when the pilot estimator is the same for the whole window. When $\max\big(\hat\eta_{\lfloor sh_n^{-1}\rfloor h_n}-\eta_{\lfloor sh_n^{-1}\rfloor h_n},\hat\sigma^{2,pilot}_{\lfloor sh_n^{-1}\rfloor h_n}-\sigma^2_{\lfloor sh_n^{-1}\rfloor h_n}\big)=\mathcal{O}_{\P}(\delta_n)$ with $\delta_n\rightarrow 0$ as $n\rightarrow \infty$, we derive that
\begin{align*}&\sum_{k=\lfloor s h_n^{-1}\rfloor +1}^{\lfloor s h_n^{-1}\rfloor+r_n^{-1}}r_n\sum_{j=1}^{J_n}\Big(w_j\Big(\hat \sigma^{2,pilot}_{\lfloor sh_n^{-1}\rfloor h_n},\hat\eta_{\lfloor sh_n^{-1}\rfloor h_n}\Big)-w_j\Big(\sigma^2_{\lfloor sh_n^{-1}\rfloor h_n},\eta_{\lfloor sh_n^{-1}\rfloor h_n}\Big)\Big)\big(S_{jk}^2-\E[S_{jk}^2]\big)\\
&\quad = r_n\sum_{j=1}^{J_n}\Big(w_j\Big(\hat \sigma^{2,pilot}_{\lfloor sh_n^{-1}\rfloor h_n},\hat\eta_{\lfloor sh_n^{-1}\rfloor h_n}\Big)-w_j\Big(\sigma^2_{\lfloor sh_n^{-1}\rfloor h_n},\eta_{\lfloor sh_n^{-1}\rfloor h_n}\Big)\Big)\,\sum_{k=\lfloor s h_n^{-1}\rfloor +1}^{\lfloor s h_n^{-1}\rfloor+r_n^{-1}} \big(S_{jk}^2-\E[S_{jk}^2]\big)\\
&\quad =\mathcal{O}_{\P}\left(r_n^{1/2}\sum_{j=1}^{J_n}\big(1+\|\Phi_{jk}\|_n^{-2}n^{-1}\big)w_j\Big(\sigma^2_{\lfloor sh_n^{-1}\rfloor h_n},\eta_{\lfloor sh_n^{-1}\rfloor h_n}\Big)\log{(n)}\delta_n\right)=\KLEINO_{\P}\big(n^{-\beta/2}\big)\,.\end{align*}
We have used that the expectation is zero and that the weights do not hinge on $k$. Then, we can bound the variance using the derivative bound \eqref{derweights}. Covariances of the $S_{jk}^2$ over different bins for $k\ne k'$ are negligible what is shown in Step 5 of the proof.
%Covariances of the $S_{jk}^2$ over different bins for $k\ne k'$ and for the same $j$ are only due to the noise parts (covariances of the discrete parts vanish by the martingale property). Not using Assumption \ref{eta2}, but only on Assumption \ref{eta} with \eqref{varrho}, we derive
%\begin{align}\label{specov}\cov\big(S_{jk},S_{jk'}\big)&=\cov\Big(\Big(\sum_{i=1}^{n-1}\epsilon_i\varphi_{jk}\Big(\frac{i+1/2}{n}\Big)n^{-1}\Big)^2,\Big(\sum_{l=1}^{n-1}\epsilon_l\varphi_{jk'}\Big(\frac{l+1/2}{n}\Big)n^{-1}\Big)^2\Big)\\
%&\notag=\sum_{i,u= k nh_n +1}^{ (k+1) nh_n}\;\sum_{l,v= k' nh_n +1}^{ (k'+1) nh_n}\E[\epsilon_i\epsilon_u\epsilon_l\epsilon_v]\varphi_{jk}^2\Big(\frac{i+1/2}{n}\Big)\varphi_{jk'}^2\Big(\frac{l+1/2}{n}\Big)n^{-4}+\KLEINO(1)\\
%&\notag=\mathcal{O}\big( ((k-k')n h_n)^{-2-\varrho}n^{-2}\big)\,,\end{align}
%where we use similar approximations as in the proof of Lemma \ref{lemmacornoise} and the Riemann approximation by $\int \varphi_{jk}^2(t)\,dt=1$ for all $(j,k)$. Summing over $k,k'$ yields that the covariances are asymptotically negligible.\\
Finally, since $r_n^{1/2}=n^{-\beta/2}\sqrt{\log{(n)}}$ some $\delta_n<n^{-\varepsilon}$ for any $\varepsilon>0$ is enough here, while we actually attain $\delta_n=n^{-\beta/2}$. It remains to bound
\begin{align*}&\sum_{k=\lfloor s h_n^{-1}\rfloor +1}^{\lfloor s h_n^{-1}\rfloor+r_n^{-1}}r_n^2\,\var\bigg(\sum_{j=1}^{J_n}\Big(w_j\Big(\hat \sigma^{2,pilot}_{kh_n},\hat\eta_{kh_n}\Big)-w_j\Big(\hat \sigma^{2,pilot}_{\lfloor sh_n^{-1}\rfloor h_n}\Big),\hat\eta_{\lfloor sh_n^{-1}\rfloor h_n}\Big)\Big(S_{jk}^2-\E[S_{jk}^2]\Big)\bigg)\\
&\quad +\sum_{k=\lfloor s h_n^{-1}\rfloor +1}^{\lfloor s h_n^{-1}\rfloor+r_n^{-1}}r_n^2\,\var\bigg(\sum_{j=1}^{J_n}\Big(w_j\big(\sigma^2_{kh_n},\eta_{kh_n}\big)-w_j\Big( \sigma^2_{\lfloor sh_n^{-1}\rfloor h_n},\eta_{\lfloor sh_n^{-1}\rfloor h_n}\Big)\Big)\Big(S_{jk}^2-\E[S_{jk}^2]\Big)\bigg)\\
& =\mathcal{O}\left(r_n\log^5{(n)}\big(n^{-\beta}\vee (r_n^{-1}h_n)^{2\alpha}\big)\right)=\KLEINO\big(n^{-\beta}\big)\,.\end{align*}
This proves \eqref{rem3}.% and completes the proof of Theorem \ref{cltspot}.
\begin{align*}
&{\bf{Step~5:}}&n^{{\beta}}\sum_{\substack{k,k'=\lfloor s h_n^{-1}\rfloor +1,\\ k\ne k'}}^{\lfloor s h_n^{-1}\rfloor+r_n^{-1}}r_n^2\cov\big(\zeta_k(\tilde C^n+\epsilon),\zeta_{k'}(\tilde C^n+\epsilon)\big)=\KLEINO(1)\,.\end{align*}
Moreover, it holds that 
\(\begin{aligned}\cov(\hat\sigma_s^2,\hat\sigma_{s-}^2)=\KLEINO(n^{-\beta})\,.\end{aligned}\)\\[.1cm]
\textit{Proof of Step 5:}
Covariances of $S_{jk}^2$ and $S_{uk'}^2$ for different bins $k\ne k'$ are only due to the noise parts and the endogeneity between noise and signal terms. All covariances of the signal parts vanish by the martingale property of $\tilde C_t^n$. Under \eqref{endo}, covariances of $S_{jk}^2$ and $S_{uk'}^2$ due to correlations between $(\epsilon_i)_{0\le i\le n}$ and $(\Delta_i^n X)_{1\le i\le n}$ are only non-zero when $|k-k'|=1$. Since there are only a finite number of indices with $|i-l|<\tilde Q$ on two neighboring bins, we obtain the bound
\begin{align*}&\cov\bigg(\Big(\sum_{i=1}^n\Delta_i^n\tilde C^{n}\Phi_{j(k-1)}\Big(\frac{i}{n}\Big)\Big)^2\,,\,\Big(\sum_{l=1}^{n-1}\epsilon_l\varphi_{uk}\Big(\frac{l+1/2}{n}\Big)\frac{1}{n}\Big)^2\bigg)\\
&=\bigg(\sum_{i=knh_n}^{knh_n+\tilde Q}\sum_{l=i-\tilde Q}^{i}\E[\epsilon_i\Delta_l^n\tilde C^n]\Phi_{j(k-1)}\Big(\frac{l}{n}\Big)\varphi_{uk}\Big(\frac{i+1/2}{n}\Big)\bigg)^2n^{-2}(1+\KLEINO(1))=\mathcal{O}(n^{-2})\,,
\end{align*}
uniformly for all $k,u,j$. We used the same approximations as for the variance terms in Step 1. For the two other covariance terms due to endogeneity, analogous estimates yield bounds of the same order. Under Assumption \ref{eta2}, a similar bound can be proved for the covariances due to serial correlation of the noise. Here, we provide a proof that does not use Assumption \ref{eta2}, but only Assumption \ref{eta} with \eqref{varrho}. We derive that
\begin{align}&\notag\cov\big(S_{jk}^2,S_{uk'}^2\big)=\|\Phi_{jk}\|_n^{-2}\|\Phi_{uk'}\|_n^{-2}\cov\Big(\Big(\sum_{i=1}^{n-1}\epsilon_i\varphi_{jk}\Big(\frac{i+1/2}{n}\Big)\frac{1}{n}\Big)^2,\Big(\sum_{l=1}^{n-1}\epsilon_l\varphi_{uk'}\Big(\frac{l+1/2}{n}\Big)\frac{1}{n}\Big)^2\Big)\\
&\notag =\sum_{i=knh_n}^{(k+1)nh_n-1}\sum_{p=knh_n-i}^{(k+1)nh_n-i-1}\sum_{l=k'nh_n}^{(k'+1)nh_n-1}\sum_{q=k'nh_n-p}^{(k'+1)nh_n-p-1}\big(\E[\epsilon_i\epsilon_{i+p}\epsilon_l\epsilon_{l+q}]-\E[\epsilon_i\epsilon_{i+p}]\E[\epsilon_l\epsilon_{l+q}]\big)\\
&\notag \quad\times \varphi_{jk}\Big(\frac{i+1/2}{n}\Big)\varphi_{jk}\Big(\frac{i+p+1/2}{n}\Big)\varphi_{uk'}\Big(\frac{l+1/2}{n}\Big)\varphi_{uk'}\Big(\frac{l+q+1/2}{n}\Big)\|\Phi_{jk}\|_n^{-2}\|\Phi_{uk'}\|_n^{-2}n^{-4}\\
%&\notag=\sum_{i,u= k nh_n +1}^{ (k+1) nh_n}\;\sum_{l,v= k' nh_n +1}^{ (k'+1) nh_n}\E[\epsilon_i\epsilon_u\epsilon_l\epsilon_v]\varphi_{jk}^2\Big(\frac{i+1/2}{n}\Big)\varphi_{jk'}^2\Big(\frac{l+1/2}{n}\Big)n^{-4}+\KLEINO(1)\\
&\notag \le K\sum_{i=knh_n}^{(k+1)nh_n-1}\sum_{l=k'nh_n}^{(k'+1)nh_n-1}((\E[\epsilon_i\epsilon_l])^2+\KLEINO(1))\|\Phi_{jk}\|_n^{-2}\|\Phi_{uk'}\|_n^{-2}n^{-4} \varphi_{jk}^2\Big(\frac{i+1/2}{n}\Big)\varphi_{uk'}^2\Big(\frac{l+1/2}{n}\Big)\\
&\notag \le K\sum_{i=knh_n}^{(k+1)nh_n-1}\sum_{l=k'nh_n}^{(k'+1)nh_n-1}((\E[\epsilon_i\epsilon_l])^2+\KLEINO(1))n^{-2} \varphi_{jk}^2\Big(\frac{i+1/2}{n}\Big)\varphi_{uk'}^2\Big(\frac{l+1/2}{n}\Big)\\
&\label{specov}=\mathcal{O}\big( ((k-k')n h_n)^{-2-\varrho}\big)\,,\end{align}
where we use similar approximations as in the proof of Lemma \ref{lemmacornoise} and that $\int \varphi_{jk}^2(t)\,dt=1$ for all $(j,k)$. Thereby, we obtain that
\begin{align*}&n^{{\beta}}\sum_{\substack{k,k'=\lfloor s h_n^{-1}\rfloor +1,\\ k\ne k'}}^{\lfloor s h_n^{-1}\rfloor+r_n^{-1}}r_n^2\cov\big(\zeta_k(\tilde C^n+\epsilon),\zeta_{k'}(\tilde C^n+\epsilon)\big)\\
&\quad =n^{\beta}r_n\sum_{\substack{k,k'=\lfloor s h_n^{-1}\rfloor +1\\ k\ne k'}}^{\lfloor s h_n^{-1}\rfloor+r_n^{-1}}r_n\sum_{u,j=1}^{J_n}w_{jk}w_{uk'}\cov(S_{jk}^2,S_{uk'}^2)\\
&\quad =\KLEINO(r_n^{-1}(\log(n))^2(nh_n)^{-2})=\KLEINO(r_n^{-1}n^{-1})=\KLEINO(1)\,.
\end{align*}
This completes the proof of the marginal central limit theorem. At the same time, we obtain analogously
\begin{align*}n^{{\beta}}\sum_{k=\lfloor s h_n^{-1}\rfloor +1}^{\lfloor s h_n^{-1}\rfloor+r_n^{-1}}\sum_{k'=\lfloor s h_n^{-1}\rfloor -r_n^{-1}}^{\lfloor s h_n^{-1}\rfloor-1}r_n^2\cov\big(\zeta_k(\tilde C^n+\epsilon),\zeta_{k'}(\tilde C^n+\epsilon)\big)=\KLEINO(1)\,.
\end{align*}
This yields that the covariances of $\hat \sigma_s^2$ and $\hat \sigma_{s-}^2$ are asymptotically negligible. We conclude the joint stable central limit theorem \eqref{cltspoteq}.
\hfill\qed

\subsection{Asymptotic theory for the test}
Denote by $\{S_1,\ldots,S_{N_1}\}$ the finite sequence of stopping times exhausting the jumps of $X$ on $[0,1]$ with $|\Delta X_{S_i}|>a$ for all $i$ and some $a\in\mathds{R}_+$ and the L\'{e}vy measure of $X$ does not have an atom in $\{a\}$. In the case of finite activity jumps, $r=0$ in Assumption \ref{H}, we can set $a=0$.
\begin{prop}\label{proptest}On the null hypothesis ${\mathds{H}}(a)_{[0,1]}$, when Assumptions \ref{eff}, \ref{eta2}, \ref{H} and \ref{sigma} are satisfied, the test statistic is asymptotically $\chi^2$-distributed,
\begin{align}n^{\beta}\,T_0(h_n,r_n,g)\stackrel{(st)}{\longrightarrow} \chi^2_{N_1}\,,\end{align}
with $N_1$ degrees of freedom .
\end{prop} 
\begin{cor}\label{corconsistent}Under the alternative hypothesis $\big(\Omega\setminus {\mathds{H}}(a)_{[0,1]}\big)$, when there exists at least one $s\in[0,1]$ with $|\Delta X_{s}|>a$ and $|\Delta\sigma^2_{s}|>0$, it holds as $n\to\infty$ that:
\begin{align}\P\Big(n^{\beta}T_0(h_n,r_n,g)>q_{1-\alpha}(\chi_{\hat N_1}^2)\Big)\to 1\,.\end{align}
\end{cor}
\begin{proof}[Proof of Proposition \ref{proptest}:]\hfill

\noindent
\textbf{1. Detection of (large) price jump arrivals}\\[.1cm]
Consider the set
\begin{align*}\tilde \Omega_n=\left\{\omega\in\Omega|S_1>r_n^{-1}h_n,S_{N_1}<1-r_n^{-1}h_n,S_i-S_{i-1}>2r_n^{-1}h_n~,i=1,\ldots,N_1-1\right\}\\
\cup\; \left\{\omega\in\Omega|S_i= k\cdot h_n~,i=1,\ldots,N_1,k=0,\ldots,h_n^{-1}\right\}^{\complement}\\
\cup\; \left\{\omega\in\Omega|\exists (s,i)~\mbox{s.t.}~|\Delta\sigma^2_s|>0~\mbox{and}~s\in [S_i-r_n^{-1}h_n,S_i+r_n^{-1}h_n]\setminus\{S_i\}\right\}^{\complement}\,.\end{align*}
Since $\P(\tilde\Omega_n)\rightarrow 1$ as $n\rightarrow \infty$ with \eqref{e5}, we work conditionally on $\tilde\Omega_n$.
The jump times $\{S_i,i=1,\ldots,N_1\}$ are estimated with thresholding by $\{\hat S_i,i=1,\ldots,\hat N_1\}$, where we set $\hat S_i=kh_n$ when $h_n|\zeta_k^{ad}(Y)|>u_n\vee a^2$. We prove that
\begin{align}\label{jumpsneg}\sum_{k=r_n^{-1}}^{h_n^{-1}-r_n^{-1}-1}g(\hat\sigma^2_{kh_n},\hat\sigma^2_{kh_n-})\1_{\{h_n|\zeta_k^{ad}(Y)|>u_n\vee a^2\}}-\sum_{s\le 1}g(\hat\sigma_s^2,\hat\sigma_{s-}^2)\1_{\{|\Delta X_s|>a\}}=\KLEINO_{\P}(n^{-\beta})\,.\end{align}
Denote $\mathcal{K}=\{0\le k\le h_n^{-1}-1|S_i\in (kh_n,(k+1)h_n)\}$ and $\mathcal{K}^{\complement}=\{0,\ldots, h_n^{-1}-1\}\setminus \mathcal{K}$. 
First, we show that
\begin{align}\label{jumps1}\sum_{k\in\mathcal{K}^{\complement}}g(\hat\sigma^2_{kh_n},\hat\sigma^2_{kh_n-})\1_{\{h_n|\zeta_k^{ad}(Y)|>u_n\vee a^2\}}=\KLEINO_{\P}(n^{-\beta})\,.\end{align}
With the Markov inequality, Lemma \ref{moments} and using that at least $p=8$ moments of the noise exist, we obtain that
\begin{align*}\P\Big(\sup_{k\in\mathcal{K}^{\complement}}|\zeta_k^{ad}(Y) |> h_n^{-1}\big(u_n\vee a^2\big)\Big)\le K h_n^{-1}\,\frac{\log(n)}{\big(u_n\vee a^2\big)^8} \,h_n^8\,,\end{align*}
for some constant $K$, and the same order without the factor $h_n^{-1}$ for $\P\big(|\zeta_k^{ad}(Y) |> h_n^{-1}(u_n\vee a^2)\big)$ and some $k\in\mathcal{K}^{\complement}$. Indicator functions $\1_{A_n}$, with $p_n=\P(A_n)\rightarrow 0$, satisfy $\1_{A_n}=\mathcal{O}_{\P}(p_n^{1/2})$, using that $\E[\1_{A_n}]=p_n$ and $\var(\1_{A_n})\le p_n$. Most factors $g(\hat\sigma^2_{kh_n},\hat\sigma^2_{kh_n-})$ in $\eqref{jumps1}$ tend to zero in probability. When $|\Delta\sigma_s^2|=0$ for all $s\in[(k-1)h_n,kh_n]$, we have that $g(\hat\sigma^2_{kh_n},\hat\sigma^2_{kh_n-})=\mathcal{O}_{\P}\big(n^{-\beta})$. However, for $k\in\mathcal{K}^{\complement}$, jumps in $(\sigma_s^2)_{0\le s\le 1 }$ can occur. From the summability of $\sum_{s\le 1}(\Delta\sigma_s^2)^2<\infty$, it follows that at most $n^{2v}$ volatility jumps of sizes bounded by $n^{-v},v\in[0,1/2)$, can occur. Since $g(x,y)=\mathcal{O}((x-y)^2)$ for \eqref{stang} as $(x-y)\to 0$, for $a=0$ we obtain that 
\begin{align*}\sum_{k\in\mathcal{K}^{\complement}}g(\hat\sigma^2_{kh_n},\hat\sigma^2_{kh_n-})\1_{\{h_n|\zeta_k^{ad}(Y)|>u_n\vee a^2\}}=\mathcal{O}_{\P}\big(n^{-\beta}\log(n)h_n^{3-4\tau}+\log(n)h_n^{4(1-\tau)}\big)\,,\end{align*}
respectively $n^{-\beta}\log(n)h_n^3+\log(n)h_n^{4}$ for $a> 0$. For $a>0$, \eqref{jumps1} is clearly satisfied, while for $a=0$ the condition 
\[3-4\tau>0 ~~\Rightarrow~~\tau<3/4\]
ensures \eqref{jumps1}. We have proven that the error due to false jump detections is asymptotically negligible.\\
It remains to prove that the error due to non-detection of one of the finitely many jump times $S_1,\ldots,S_{N_1}$ is also asymptotically negligible. This is ensured by
\begin{align}\label{jumps2}-\sum_{k\in\mathcal{K}}g(\hat\sigma^2_{kh_n},\hat\sigma^2_{kh_n-})\1_{\{h_n|\zeta_k^{ad}(Y)|\le u_n\vee a^2\}}=\KLEINO_{\P}(n^{-\beta})\,.\end{align}
By the results from Section 3.1.3 of \cite{bibwink2015},\footnote{See Proposition 3.2. of \cite{bibwink2015}.} for $S_i\in ((k-1)h_n,kh_n)$, it holds that
\[h_n|\zeta_k^{ad}(Y)|=(\Delta X_{S_i})^2+\xi_i=a^2+\epsilon+\xi_i~~\mbox{with}~~ \xi_i=\KLEINO_{\P}(1)~\mbox{and}~\epsilon>0\,.\footnote{Since the L\'{e}vy measure of $X$ does not have an atom in $\{a\}$.}\]
On the hypothesis, there are no simultaneous jumps in the volatility, i.e.\ $\sigma_{S_i}^2-\sigma_{S_i-}^2=0$ for all $i=1,\ldots,N_1$. On the finitely many bins with $k\in\mathcal{K}$, we thus have that
%\[g(\hat\sigma^2_{kh_n},\hat\sigma^2_{kh_n-})=\mathcal{O}_{\P}\big(n^{-\beta}+\1_{\{\sup_{s\in((k-r_n^{-1})h_n,(k+r_n^{-1})h_n)}|\Delta \sigma_s^2|>0\}}\big)\,.\]
%It holds that
\[\sup_{k\in\mathcal{K}} g(\hat\sigma^2_{kh_n},\hat\sigma^2_{kh_n-})=\mathcal{O}_{\P}(n^{-\beta})\,.\]
Hence, $\xi_i=\KLEINO_{\P}(1)$ suffices to ensure \eqref{jumps2}. \eqref{jumps1} and \eqref{jumps2} imply \eqref{jumpsneg}.\\[.4cm]
\textbf{2. Stable convergence of spot volatility estimates around detected (large) price-jump times}\\[.1cm]
The asymptotic distribution of the test statistic is derived with \eqref{jumpsneg} and the stable convergences of the spot volatility estimates:
\begin{align*}n^{\beta/2}\left(\begin{matrix}\hat \sigma^2_s-\sigma^2_s\\ \hat \sigma^2_{s-}-\sigma^2_{s-}\end{matrix}\right)\stackrel{(st)}{\longrightarrow} MN\left(0,\left(\begin{matrix}8\sigma_s^{3}\eta_s^{1/2}&0\\0&8\sigma_{s-}^{3}\eta_s^{1/2}\end{matrix}\right)\right)\,,\end{align*}
%\begin{align*}n^{\nicefrac{\beta}{2}}\Big(\hat \sigma^2_{S_i}-\sigma^2_{S_i}\Big)\stackrel{(st)}{\longrightarrow}MN\big(0,8\sigma^3_{S_i}\eta_{S_i}^{1/2}\big)\,,\end{align*}
%\begin{align*}n^{\nicefrac{\beta}{2}}\Big(\hat \sigma^2_{S_i -}-\sigma^2_{S_i-}\Big)\stackrel{(st)}{\longrightarrow}MN\big(0,8\sigma_{S_i-}^{3}\eta_{S_i-}^{1/2}\big)\,,\end{align*}
%\begin{align*}n^{\beta}\Bigg(\sum_{k=r_n^{-1}}^{h_n^{-1}-r_n^{-1}-1}\1_{\{h_n|\zeta_k^{ad}(Y)|>u_n\vee a^2\}}-\sum_{s\le 1}\1_{\{|\Delta X_s|>a\}}\Bigg)\stackrel{\P}{\longrightarrow} 0\,,\end{align*}
which hold jointly for all $i=1,\ldots,N_1$. The stable limit theorems of the spot volatility estimators are given in Theorem \ref{cltspot}. Concerning the convergence of the spot estimates at stopping times, observe that
\begin{itemize}
\item Thresholding and identification of a jump is based on $\zeta_k^{ad}(Y)$.
\item Given that $h_n|\zeta_k^{ad}(Y)|>u_n\vee a^2$, $\hat S_i=k h_n$ for some $i\in\{1,\ldots,\hat N_1\}$, $\hat \sigma^2_{\hat S_i}$ is computed from $\zeta_l^{ad}(Y),l=(k+1),\ldots,(k+r_n^{-1})$.
\item Given $\hat S_i=k h_n$, $\hat \sigma^2_{\hat S_i-}$ is computed from $\zeta_l^{ad}(Y),l=(k-r_n^{-1}),\ldots,(k-1)$.
\end{itemize}
We restrict to $\tilde\Omega_n$ again. For the stability of weak convergences, we have already considered a sequence of stopping times in Step 1 of the proof of Theorem \ref{cltspot}. Recall the definition of $\tilde{\mathcal{G}}_t^n$ from this paragraph. The $S_p,p=1,\ldots,N_1$, are $\tilde{\mathcal{G}}_0^n$-measurable random variables and denote $i_p$ integer-valued $\tilde{\mathcal{G}}_0^n$-measurable random variables such that $i_ph_n<S_p< (i_p+1) h_n$. The stable limit theorem in Theorem \ref{cltspot} is valid when replacing the fixed time $s$ by stopping times $S_p,p=1,\ldots,N_1$. Analogously as in Lemma 8.1 of \cite{jacodtodorov}, this readily follows with the points above by the asymptotic independence of the statistics in Step 1 of the proof of Theorem \ref{cltspot} with $s=S_p$ for $\hat\sigma^2_{S_p}$, or $s=i_p h_n$ for $\hat\sigma^2_{i_p h_n}$ respectively, from $\mathcal{F}_{S_p}$. Here, we exploit that the noise is under Assumption \ref{eta2} only weakly serially dependent over asymptotically decreasing intervals and only dependent on finitely many preceding increments of $X$, and the strong Markov property of Brownian motion.\\
On assumption \ref{sigma}, $\max_p|\sigma^2_{i_p h_n}-\sigma^2_{S_p}|=\mathcal{O}_{\P}(h_n^{\alpha})=\KLEINO_{\P}(n^{-\alpha/2})$, the latter being much smaller than $n^{-\beta/2}$. Therefore, a discretization of estimated jump arrivals is asymptotically negligible.\\
Moreover, on $\tilde \Omega_n$ all spot squared volatility estimates are computed from \emph{disjoint} data subsets. 
 % In particular this applies to left and right-hand estimates at a particular jump time $\tau_i$. 
Therefore, by \eqref{specov}, covariations between all estimates converge to zero in probability what implies joint weak convergence.\footnote{Note that by Step 4 in the proof of Theorem \ref{cltspot}, this is still true if the pre-estimated noise long-run variance was computed from all observations.} Stability of the convergence of the vector has been established above in Step 1 of the proof of Theorem \ref{cltspot}.\\[.1cm]
\textbf{3. Convergence of the test statistic}\\[.1cm]
For test functions which are twice continuously differentiable with bounded second derivatives, Taylor's formula yields
\begin{align*}g(x_1,x_2)-g(a_1,a_2)&=\frac{\partial g}{\partial x_1}(a_1,a_2)(x_1-a_1)+\frac{\partial g}{\partial x_2}(a_1,a_2)(x_2-a_2)+\frac{\partial^2 g}{2\,\partial x_1^2}(a_1,a_2)(x_1-a_1)^2\\
&+\frac{\partial^2 g}{2\,\partial x_2^2}(a_1,a_2)(x_2-a_2)^2+\frac{\partial^2 g}{\partial x_1\partial x_2}(a_1,a_2)(x_1-a_1)(x_2-a_2)\\
&+\KLEINO\big(\max\big((x_1-a_1)^2,(x_2-a_2)^2\big)\big)\,.\end{align*}
We apply the generalized $\Delta$-method for stable convergence and set $(a_1,a_2)=(\sigma^2_{S_i},\sigma^2_{S_i-})$ and the random vector $(x_1,x_2)=(\hat \sigma^2_{S_i},\hat \sigma^2_{S_i-})$ with estimators \eqref{rspot} and \eqref{lspot} at the times $S_i,S_{i-},i=1,\ldots, N_1$. When we focus on the test function \eqref{stang} in Theorem \ref{test}, it holds that
\[\frac{\partial g}{\partial x_1}(\sigma^2_{S_i},\sigma^2_{S_i})=\frac{\partial g}{\partial x_2}(\sigma^2_{S_i},\sigma^2_{S_i})= g(\sigma^2_{S_i},\sigma^2_{S_i})=0\,.\]
The second order term comes into play and the equalities
\begin{align}\label{secd}\frac{\partial^2 g}{\partial x_1^2}(\sigma^2_{S_i},\sigma^2_{S_i})=\frac{\partial^2 g}{\partial x_2^2}(\sigma^2_{S_i},\sigma^2_{S_i})= -\frac{\partial^2 g}{\partial x_1\partial x_2}(\sigma^2_{S_i},\sigma^2_{S_i})=\frac{1}{8}\sigma^{-3}_{S_i}\,.\end{align}
Under Assumption \ref{eta2} we have by Proposition \ref{cornoiseest} estimators $\hat\eta_{kh_n}^{1/2}=\eta_{kh_n}^{1/2}+\KLEINO_{\P}(n^{-\beta})$ for all $k$. This renders the estimation errors of $\hat\eta_{kh_n}^{-1/2}$ in \eqref{teststatistic} asymptotically negligible in \eqref{testlt}.\\
Cram\'{e}r-Wold's theorem gives equivalence of the weak convergence of the vector $(\hat\sigma_{S_i}^2,\hat\sigma_{S_i-}^2)_{1\le i\le N_1}$ to weak convergence of linear combinations. 
Under $\mathds{H}(a)_{[0,1]}$, when $\sigma^2_{S_i}=\sigma^2_{S_i-}$ for all $i$, the limit of $n^{\beta}T_0(h_n,r_n,g)$ can thus be described by a random variable
\[\sum_{i=1}^{N_1}\Bigg(\frac{\partial^2 g}{2\,\partial x_1^2}(\sigma^2_{S_i},\sigma^2_{S_i}
)Z_i^2+\frac{\partial^2 g}{2\,\partial x_2^2}(\sigma^2_{S_i},\sigma^2_{S_i}
)\tilde Z_i^2+\frac{\partial^2 g}{\partial x_1\partial x_2}(\sigma^2_{S_i},\sigma^2_{S_i})Z_i\tilde Z_i
\Bigg)\,8 \sigma_{S_i}^{3}\,,\]
where $Z_i$ and $\tilde Z_i$, $i=1,\ldots,N_1$, are two independent collections of i.i.d.\,standard normals defined on the orthogonal extension of $(\Omega,\mathcal{F},\P)$ in the product space that accommodates all random variables. Since $(1/\sqrt{2})(Z_i-\tilde Z_i)$ are i.i.d.\,standard normals, the $\chi^2$-distribution with $N_1$ degrees of freedom appears as limiting distribution. Proposition \ref{proptest} follows with the binomial formula and by the second derivatives of the test function \eqref{stang} from \eqref{secd}. Even though the limit above could depend on the particular choice of stopping times its $\mathcal{F}$-conditional law does not.
\end{proof}
\begin{proof}[Proof of Corollary \ref{corconsistent}:]\hfill\\
Under the alternative hypothesis, $\sigma_{S_i}^2\ne \sigma_{S_i-}^2$, for at least one $i\in\{1,\ldots,N_1\}$. In this case, we have that
\[n^{\beta}T_0(h_n,r_n,g)=\mathcal{O}_{\P}(1)+n^{\beta}\hat\eta^{-1/2}_{\lfloor S_ih_n^{-1}\rfloor h_n}g\big(\hat\sigma^2_{\lfloor S_ih_n^{-1}\rfloor h_n},\hat\sigma^2_{\lfloor S_ih_n^{-1}\rfloor h_n -}\big)\1_{\{h_n|\zeta_{\lfloor S_ih_n^{-1}\rfloor h_n}^{ad}(Y)|>(u_n\vee a^2)\}}\]
with Proposition \ref{proptest}. Since 
\[g\big(\hat\sigma^2_{\lfloor S_ih_n^{-1}\rfloor h_n},\hat\sigma^2_{\lfloor S_ih_n^{-1}\rfloor h_n -}\big)\ge c\,\Delta\sigma^2_{S_i}-\KLEINO_{\P}(1)\]
for some constant $c$ and since \(h_n|\zeta_k^{ad}(Y)|=(\Delta X_{S_i})^2+\xi_i=a^2+\epsilon+\xi_i\) with $\xi_i=\KLEINO_{\P}(1)$ and some $\epsilon>0$, we conclude with the reverse triangle inequality that 
\[\P\Big(n^{\beta}T_0(h_n,r_n,g)>q_{1-\alpha}(\chi_{\hat N_1}^2)\Big)\to 1\]
for any arbitrarily small $\alpha>0$. This proves Corollary \ref{corconsistent}.
\end{proof}
\section*{Acknowledgement}
\noindent
The authors are grateful to two anonymous referees whose valuable comments and suggestions helped improving the paper.
Lars Winkelmann acknowledges financial support from the Deutsche Forschungsgemeinschaft via CRC 649 `\"Okonomisches Risiko', Humboldt-Universität zu Berlin.
%\end{appendix}
%\thispagestyle{plain}
%\section{Introduction\label{sec:1}}
%\noindent
%\section*{\refname}
\bibliographystyle{chicago}
\bibliography{ref}

\begin{thebibliography}{}

\bibitem[\protect\citeauthoryear{A\"{\i}t-Sahalia, Fan, Laeven, Wang, and
  Yang}{A\"{\i}t-Sahalia et~al.}{2017}]{lev}
A\"{\i}t-Sahalia, Y., J.~Fan, R.~J.~A. Laeven, C.~D. Wang, and X.~Yang (2017).
\newblock Estimation of the continuous and discontinuous leverage effects.
\newblock {\em Journal of the American Statistical Association,\/}~{\em
  112\/}(520), 1744--1758.

\bibitem[\protect\citeauthoryear{A{\"\i}t-Sahalia and Jacod}{A{\"\i}t-Sahalia
  and Jacod}{2010}]{aitjac10}
A{\"\i}t-Sahalia, Y. and J.~Jacod (2010).
\newblock Is {B}rownian motion necessary to model high-frequency data?
\newblock {\em The Annals of Statistics\/}~{\em 38\/}(5), 3093--3128.

\bibitem[\protect\citeauthoryear{{A\"{\i}t-Sahalia} and
  {Jacod}}{{A\"{\i}t-Sahalia} and {Jacod}}{2014}]{sahaliajacod}
{A\"{\i}t-Sahalia}, Y. and J.~{Jacod} (2014).
\newblock {\em {High-frequency financial econometrics.}}
\newblock Princeton, NJ: Princeton University Press.

\bibitem[\protect\citeauthoryear{A\"{\i}t-Sahalia, Zhang, and
  Mykland}{A\"{\i}t-Sahalia et~al.}{2005}]{howoften}
A\"{\i}t-Sahalia, Y., L.~Zhang, and P.~A. Mykland (2005).
\newblock How often to sample a continuous-time process in the presence of
  market microstructure noise.
\newblock {\em Review of Financial Studies\/}~{\em 18}, 351--416.

\bibitem[\protect\citeauthoryear{Altmeyer and Bibinger}{Altmeyer and
  Bibinger}{2015}]{stable}
Altmeyer, R. and M.~Bibinger (2015).
\newblock Functional stable limit theorems for quasi-efficient spectral
  covolatility estimators.
\newblock {\em Stochastic Processes and their Applications\/}~{\em 125\/}(12),
  4556--4600.

\bibitem[\protect\citeauthoryear{Andersen and Bollerslev}{Andersen and
  Bollerslev}{1998}]{andersenbollerslev98}
Andersen, T.~G. and T.~Bollerslev (1998).
\newblock Answering the skeptics: Yes, standard volatility models do provide
  accurate forecasts.
\newblock {\em International Economic Review\/}~{\em 39\/}(4), 885--905.

\bibitem[\protect\citeauthoryear{Andersen, Bollerslev, Diebold, and
  Labys}{Andersen et~al.}{2001}]{abdl01}
Andersen, T.~G., T.~Bollerslev, F.~X. Diebold, and P.~Labys (2001).
\newblock The distribution of realized exchange rate volatility.
\newblock {\em Journal of the American Statistical Association\/}~{\em
  96\/}(453), 42--55.

\bibitem[\protect\citeauthoryear{Bandi and Ren\`{o}}{Bandi and
  Ren\`{o}}{2016}]{bandireno}
Bandi, F. and R.~Ren\`{o} (2016).
\newblock Price and volatility co-jumps.
\newblock {\em Journal of Financial Economics\/}~{\em 119\/}(1), 107--146.

\bibitem[\protect\citeauthoryear{Barndorff-Nielsen, Hansen, Lunde, and
  Shephard}{Barndorff-Nielsen et~al.}{2008}]{bn2}
Barndorff-Nielsen, O.~E., P.~R. Hansen, A.~Lunde, and N.~Shephard (2008).
\newblock Designing realised kernels to measure the ex-post variation of equity
  prices in the presence of noise.
\newblock {\em Econometrica\/}~{\em 76\/}(6), 1481--1536.

\bibitem[\protect\citeauthoryear{Barndorff-Nielsen and
  Shephard}{Barndorff-Nielsen and Shephard}{2002}]{bn3}
Barndorff-Nielsen, O.~E. and N.~Shephard (2002).
\newblock Econometric analysis of realized volatility and its use in estimating
  stochastic volatility models.
\newblock {\em Journal of the Royal Statistical Society\/}~{\em 64\/}(2),
  253--280.

\bibitem[\protect\citeauthoryear{Bibinger}{Bibinger}{2011}]{bibinger}
Bibinger, M. (2011).
\newblock Efficient covariance estimation for asynchronous noisy high-frequency
  data.
\newblock {\em Scandinavian Journal of Statistics\/}~{\em 38}, 23--45.

\bibitem[\protect\citeauthoryear{Bibinger, Hautsch, Malec, and
  Rei{\ss}}{Bibinger et~al.}{2014}]{BHMR}
Bibinger, M., N.~Hautsch, P.~Malec, and M.~Rei{\ss} (2014).
\newblock Estimating the quadratic covariation matrix from noisy observations:
  Local method of moments and efficiency.
\newblock {\em The Annals of Statistics\/}~{\em 42\/}(4), 1312--1346.

\bibitem[\protect\citeauthoryear{Bibinger, Hautsch, Malec, and
  Rei{\ss}}{Bibinger et~al.}{2017}]{BHMR2}
Bibinger, M., N.~Hautsch, P.~Malec, and M.~Rei{\ss} (2017).
\newblock Estimating the spot covariation of asset prices -- statistical theory
  and empirical evidence.
\newblock {\em Journal of Business \& Economic Statistics,\/}~{\em
  forthcoming}.

\bibitem[\protect\citeauthoryear{Bibinger, Jirak, and Vetter}{Bibinger
  et~al.}{2017}]{BJV}
Bibinger, M., M.~Jirak, and M.~Vetter (2017).
\newblock Nonparametric change-point analysis of volatility.
\newblock {\em The Annals of Statistics\/}~{\em 45\/}(4), 1542--1578.

\bibitem[\protect\citeauthoryear{Bibinger and Rei{\ss}}{Bibinger and
  Rei{\ss}}{2014}]{bibingerreiss}
Bibinger, M. and M.~Rei{\ss} (2014).
\newblock Spectral estimation of covolatility from noisy observations using
  local weights.
\newblock {\em Scandinavian Journal of Statistics\/}~{\em 41\/}(1), 23--50.

\bibitem[\protect\citeauthoryear{Bibinger and Winkelmann}{Bibinger and
  Winkelmann}{2015}]{bibwink2015}
Bibinger, M. and L.~Winkelmann (2015).
\newblock Econometrics of co-jumps in high-frequency data with noise.
\newblock {\em Journal of Econometrics\/}~{\em 184\/}(2), 361 -- 378.

\bibitem[\protect\citeauthoryear{Bloom}{Bloom}{2009}]{bloom}
Bloom, N. (2009).
\newblock The impact of uncertainty shocks.
\newblock {\em Econometrica\/}~{\em 77\/}(3), 623--685.

\bibitem[\protect\citeauthoryear{Clinet and Potiron}{Clinet and
  Potiron}{2017}]{clinet1}
Clinet, S. and Y.~Potiron (2017).
\newblock Efficient asymptotic variance reduction when estimating volatility in
  high frequency data.
\newblock {\em arxive:1701.01185\/}.

\bibitem[\protect\citeauthoryear{Comte and Renault}{Comte and
  Renault}{1998}]{comte}
Comte, F. and E.~Renault (1998).
\newblock Long memory in continuous-time stochastic volatility models.
\newblock {\em Mathematical Finance\/}~{\em 8\/}(4), 291--323.

\bibitem[\protect\citeauthoryear{Duffie, Pan, and Singleton}{Duffie
  et~al.}{2000}]{duffie}
Duffie, D., J.~Pan, and K.~Singleton (2000).
\newblock Transform analysis and asset pricing for affine jump-diffusions.
\newblock {\em Econometrica\/}~{\em 68\/}(6), 1343--1376.

\bibitem[\protect\citeauthoryear{Fan and Wang}{Fan and Wang}{2007}]{fanwang}
Fan, J. and Y.~Wang (2007).
\newblock Multi-scale jump and volatility analysis for high-frequency data.
\newblock {\em Journal of the American Statistical Association\/}~{\em
  102\/}(480), 1349--1362.

\bibitem[\protect\citeauthoryear{Hansen and Lunde}{Hansen and
  Lunde}{2006}]{hans06}
Hansen, P.~R. and A.~Lunde (2006).
\newblock Realized variance and market microstructure noise.
\newblock {\em Journal of Business \& Economic Statistics\/}~{\em 24\/}(2),
  127--161.

\bibitem[\protect\citeauthoryear{Hautsch and Podolskij}{Hautsch and
  Podolskij}{2013}]{haut13}
Hautsch, N. and M.~Podolskij (2013).
\newblock Preaveraging-based estimation of quadratic variation in the presence
  of noise and jumps: Theory, implementation, and empirical evidence.
\newblock {\em Journal of Business \& Economic Statistics\/}~{\em 31\/}(2),
  165--183.

\bibitem[\protect\citeauthoryear{Jacod}{Jacod}{2008}]{jacodjumps}
Jacod, J. (2008).
\newblock Asymptotic properties of realized power variations and related
  functionals of semimartingales.
\newblock {\em Stochastic Processes and their Applications\/}~{\em 118\/}(4),
  517--559.

\bibitem[\protect\citeauthoryear{Jacod, Kl\"uppelberg, and M\"uller}{Jacod
  et~al.}{2017}]{gernot}
Jacod, J., C.~Kl\"uppelberg, and G.~M\"uller (2017).
\newblock Testing for non-correlation between price and volatility jumps.
\newblock {\em Journal of Econometrics\/}~{\em 197\/}(2), 284--297.

\bibitem[\protect\citeauthoryear{Jacod, Li, Mykland, Podolskij, and
  Vetter}{Jacod et~al.}{2009}]{JLMPV}
Jacod, J., Y.~Li, P.~A. Mykland, M.~Podolskij, and M.~Vetter (2009).
\newblock Microstructure noise in the continous case: the pre-averaging
  approach.
\newblock {\em Stochastic Processes and their Applications\/}~{\em 119},
  2803--2831.

\bibitem[\protect\citeauthoryear{Jacod and Mykland}{Jacod and
  Mykland}{2015}]{jacodmykland}
Jacod, J. and P.~A. Mykland (2015).
\newblock Microstructure noise in the continuous case: Approximate efficiency
  of the adaptive pre-averaging method.
\newblock {\em Stochastic Processes and their Applications\/}~{\em 125}, 2910
  -- 2936.

\bibitem[\protect\citeauthoryear{Jacod and Protter}{Jacod and
  Protter}{2012}]{JP}
Jacod, J. and P.~Protter (2012).
\newblock {\em Discretization of processes}.
\newblock Springer.

\bibitem[\protect\citeauthoryear{Jacod and Todorov}{Jacod and
  Todorov}{2010}]{jacodtodorov}
Jacod, J. and V.~Todorov (2010).
\newblock Do price and volatility jump together?
\newblock {\em The Annals of Applied Probability\/}~{\em 20\/}(4), 1425--1469.

\bibitem[\protect\citeauthoryear{Kalnina and Xiu}{Kalnina and Xiu}{2017}]{xiu2}
Kalnina, I. and D.~Xiu (2017).
\newblock Nonparametric {E}stimation of the {L}everage {E}ffect: {A}
  {T}rade-{O}ff {B}etween {R}obustness and {E}fficiency.
\newblock {\em Journal of the American Statistical Association\/}~{\em
  112\/}(517), 384--396.

\bibitem[\protect\citeauthoryear{Koike}{Koike}{2016}]{koike2016}
Koike, Y. (2016).
\newblock Estimation of integrated covariances in the simultaneous presence of
  nonsynchronicity, microstructure noise and jumps.
\newblock {\em Econometric Theory\/}~{\em 32\/}(3), 533--611.

\bibitem[\protect\citeauthoryear{Lee and Mykland}{Lee and
  Mykland}{2008}]{leemykland}
Lee, S. and P.~A. Mykland (2008).
\newblock Jumps in finacial markets: A new nonparametric test and jump
  dynamics.
\newblock {\em Review of Financial Studies\/}~{\em 21\/}(6), 2535--2563.

\bibitem[\protect\citeauthoryear{Lee and Mykland}{Lee and
  Mykland}{2012}]{leemykland2}
Lee, S. and P.~A. Mykland (2012).
\newblock Jumps in equilibrium prices and market microstructure noise.
\newblock {\em Journal of Econometrics\/}~{\em 168\/}(2), 396--406.

\bibitem[\protect\citeauthoryear{Liu, Longstaff, and Pan}{Liu
  et~al.}{2003}]{liufinance}
Liu, J., F.~Longstaff, and J.~Pan (2003).
\newblock Dynamic asset allocation with event risk.
\newblock {\em Journal of Finance\/}~{\em 58\/}(1), 231--259.

\bibitem[\protect\citeauthoryear{Mancini}{Mancini}{2009}]{mancini}
Mancini, C. (2009).
\newblock Non-parametric threshold estimation for models with stochastic
  diffusion coefficient and jumps.
\newblock {\em Scandinavian Journal of Statistics\/}~{\em 36\/}(4), 270--296.

\bibitem[\protect\citeauthoryear{Mancini, Mattiussi, and Reno}{Mancini
  et~al.}{2015}]{manc15}
Mancini, C., V.~Mattiussi, and R.~Reno (2015).
\newblock Spot volatility estimation using delta sequences.
\newblock {\em Finance and Stochastics\/}~{\em 19\/}(2), 261--293.

\bibitem[\protect\citeauthoryear{Munk and Schmidt-Hieber}{Munk and
  Schmidt-Hieber}{2010a}]{munk2010}
Munk, A. and J.~Schmidt-Hieber (2010a).
\newblock Lower bounds for volatility estimation in microstructure noise
  models.
\newblock In J.~O. Berger, T.~T. Cai, and I.~M. Johnstone (Eds.), {\em
  Borrowing Strength: Theory Powering Applications -- A Festschrift for
  Lawrence D. Brown}, Volume~6 of {\em Collections}, pp.\  43--55. Beachwood,
  Ohio, USA: Institute of Mathematical Statistics.

\bibitem[\protect\citeauthoryear{Munk and Schmidt-Hieber}{Munk and
  Schmidt-Hieber}{2010b}]{munk2010b}
Munk, A. and J.~Schmidt-Hieber (2010b).
\newblock Nonparametric estimation of the volatility function in a
  high-frequency model corrupted by noise.
\newblock {\em Electronic Journal of Statistics\/}~{\em 4}, 781--821.

\bibitem[\protect\citeauthoryear{Pastor and Veronesi}{Pastor and
  Veronesi}{2012}]{pastor}
Pastor, L. and P.~Veronesi (2012).
\newblock Uncertainty about government policy and stock prices.
\newblock {\em Journal of Finance\/}~{\em 67\/}(4), 1219--1264.

\bibitem[\protect\citeauthoryear{Rei\ss}{Rei\ss}{2011}]{reiss}
Rei\ss, M. (2011).
\newblock Asymptotic equivalence for inference on the volatility from noisy
  observations.
\newblock {\em The Annals of Statistics\/}~{\em 39\/}(2), 772--802.

\bibitem[\protect\citeauthoryear{Tauchen and Todorov}{Tauchen and
  Todorov}{2011}]{voljumps}
Tauchen, G. and V.~Todorov (2011).
\newblock Volatility jumps.
\newblock {\em Journal of Business and Economic Statistics\/}~{\em 29\/}(3),
  356--371.

\bibitem[\protect\citeauthoryear{Todorov}{Todorov}{2010}]{premium}
Todorov, V. (2010).
\newblock Variance risk-premium dynamics: The role of jumps.
\newblock {\em Review of Financial Studies\/}~{\em 23\/}(1), 345--383.

\bibitem[\protect\citeauthoryear{Winkelmann, Bibinger, and Linzert}{Winkelmann
  et~al.}{2016}]{ecb}
Winkelmann, L., M.~Bibinger, and T.~Linzert (2016).
\newblock Ecb monetary policy surprises: Identification through cojumps in
  interest rates.
\newblock {\em Journal of Applied Econometrics\/}~{\em 31\/}(4), 613--629.

\bibitem[\protect\citeauthoryear{Zhang}{Zhang}{2006}]{zhang}
Zhang, L. (2006).
\newblock Efficient estimation of stochastic volatility using noisy
  observations: A multi-scale approach.
\newblock {\em Bernoulli\/}~{\em 12\/}(6), 1019--1043.

\bibitem[\protect\citeauthoryear{Zu and Boswijk}{Zu and Boswijk}{2014}]{zu14}
Zu, Y. and H.~P. Boswijk (2014).
\newblock Estimating spot volatility with high-frequency financial data.
\newblock {\em Journal of Econometrics\/}~{\em 181\/}(2), 117 -- 135.

\end{thebibliography}

\end{document}